\documentclass[11pt]{article}
\usepackage{amsmath,amssymb,color}

\definecolor{mycolorred}{rgb}{1, 0, 0}

\def\I{\mbox{\large \bf 1}}
\def\P{{\mathbb P}}
\def\R{{\mathbb R}}
\def\E{{\mathbb E}}

\def\N{{\mathbb N}}
\def\U{{\mathcal{U}}}

\def\<{\langle}
\def\>{\rangle}

\newtheorem{theorem}{Theorem}[section]

\newtheorem{assumption}[theorem]{Assumption}

\newtheorem{lemma}[theorem]{Lemma}

\newtheorem{proposition}[theorem]{Proposition}
\newtheorem{remark}[theorem]{Remark}

\numberwithin{equation}{section}

\begin{document}

\title{An invariance principle for stochastic series I. Gaussian limits}
\author{ \textsc{Vlad Bally}\thanks{%
Universit\'e Paris-Est, LAMA (UMR CNRS, UPEMLV, UPEC), INRIA, F-77454
Marne-la-Vall\'ee, France. Email: \texttt{bally@univ-mlv.fr}.} \smallskip \\
\textsc{Lucia Caramellino}\thanks{%
Dipartimento di Matematica, Universit\`a di Roma - Tor Vergata, Via della
Ricerca Scientifica 1, I-00133 Roma, Italy. Email: \texttt{%
caramell@mat.uniroma2.it}}\smallskip\\
}
\maketitle

\begin{abstract}
We study invariance principles and convergence to a Gaussian limit for
stochastic series of the form $S(c,Z)=\sum_{m=1}^{\infty }\sum_{\alpha
_{1}<...<\alpha _{m}}c(\alpha _{1},...,\alpha _{m})\prod_{i=1}^{m}Z_{\alpha
_{i}}$ where $Z_{k}$, $k\in \N$ is a sequence of centred independent random
variables of unit variance. In the case when the $Z_{k}$'s are Gaussian, $%
S(c,Z)$ is an element of the Wiener chaos and convergence to a Gaussian
limit (so the corresponding nonlinear CLT) has been intensively studied by  Nualart, Peccati, Nourdin and several other authors. The invariance principle
consists in taking $Z_{k}$ with a general law. It  has also been considered in
the literature, starting from the seminal papers of Jong, and a variety of
applications including $U$-statistics are of interest. Our main
contribution is to study the convergence in total variation distance and to
give estimates of the error.
\end{abstract}

\tableofcontents

\section{Introduction}

The aim of this paper is to provide invariance principles for stochastic
series of the form%
\begin{equation}
S_{N}(c,Z)=\sum_{m=1}^{N}\Phi _{m}(c,Z),\quad \mbox{with}\quad \Phi
_{m}(c,Z)=\sum_{\left\vert \alpha \right\vert =m}c(\alpha )Z^{\alpha }.
\label{Int1}
\end{equation}
Let us explain the notation: $\alpha =(\alpha _{1},...,\alpha _{m})\in \N^{m}$
is a multi-index with length $\left\vert \alpha \right\vert =m;$ $Z_{k}$, $k\in
\N$, is a sequence of independent random variables with $\E(Z_{k})=0$ and $%
\E(Z_{k}^{2})=1$ and $Z^{\alpha }=\prod_{i=1}^{m}Z_{\alpha _{i}}.$ Moreover $%
c(\alpha )\in \R$ are coefficients which are symmetric, null on the
diagonals and verify the normalization condition%
\begin{equation}
\E(S_{N}(c,Z)^2)=\sum_{m=1}^{N}m!\sum_{\left\vert \alpha \right\vert
=m}c^{2}(\alpha )=1.  \label{Int1a}
\end{equation}%
There are two types of results: first we consider infinite series, so $%
N=\infty $, and sequences of coefficients $c^{(n)}=(c^{(n)}(\alpha ))_{\alpha }$
which verify the above normalization condition and we give sufficient
conditions for the convergence of $S_{\infty}(c^{(n)},Z)$ to a standard Gaussian
law. This will be convergence in law on one hand and convergence in total
variation distance on the other hand. These are asymptotic results. In a
second stage we restrict ourself to finite series, so $N<\infty $ is fixed,
and we obtain non asymptotic estimates of the error. Here our first aim is to estimate
\begin{equation*}
\Delta _{Z,\overline{Z}}(c,f)=\big| \E(f(S_{N}(c,Z)))-\E(f(S_{N}(c,%
\overline{Z})))\big|
\end{equation*}%
where $\overline{Z}_{k}$, $k\in \N$, is a sequence of independent standard normal
random variables. The second aim is to estimate the distance between the law
of $S_{N}(c,Z)$ and the standard Gaussian law. In the case $N=1$ this is
just the CLT. Notice also that since $\overline{Z}_{k}$, $k\in \N$, are
standard normal random variables then $\Phi _{m}(c,\overline{Z})$ is (in law) a
multiple stochastic integral of order $m$ and, starting with the seminal
paper of Nualart and Peccati \cite{[PN]}, a lot of work has been done in
order to obtain the CLT for such multiple integrals. So, once we are able
to estimate $\Delta _{Z,\overline{Z}}(c,f)$ (this is the invariance
principle), we may use the above mentioned results concerning the Wiener
chaos, in order to obtain the distance to the Gaussian distribution.
However, the two problems have to be discussed separately because the
Gaussian law is not the single possible limit for such series: for example,
Nourdin and Peccati in \cite{[NP1]} give sufficient conditions in order that
such series converge to a chi-squared distribution. We address the
problem of non Gaussian limits in the working paper \cite{[BCNon]}.

This type of nonlinear invariance principle turns out to be of interest in
several very different fields of applications: Mossel, O'Donnell and  Oleszkiewicz in \cite{[MDO]}
provided interesting applications in theoretical computer science
and in social choice theory. And similar objects appear in the $U$-statistics theory see e.g. Koroljuk and Borovskich \cite{kb}.

The first results concerning the convergence in law of $S_{N}(c,Z)$ to the
Gaussian distribution has been obtained by Jong \cite{[dJ1]} and \cite%
{[dJ2]}. Afterwards, Mossel, O'Donnell and  Oleszkiewicz in \cite{[MDO]} obtained
an invariance principle in Kolmogorov distance. Finally, under a
supplementary regularity condition on the laws of $Z_{k}$ (that we discuss
below) Nourdin and Poly \cite{[NPy]} gave a convergence result in
total variation distance. Let us shortly present these results. The central
quantity which controls the convergence of the series $S_{N}(c,Z)$ is the so
called ``low influence factor'' defined by%
\begin{equation}\label{Int1b}
\overline{\delta }_{N}(c)=\sum_{m=1}^{N}\delta _{m}(c)\quad \mbox{with}\quad \delta
_{m}(c)=\max_{k}\sum_{\left\vert \alpha \right\vert =m-1}c^{2}(k,\alpha ).
\end{equation}%
Roughly speaking $\sum_{\left\vert \alpha \right\vert =m-1}c^{2}(k,\alpha )$
may be considered as the influence on the particle $k$ of all the other
particles. And if $\overline{\delta }_{N}(c)$ is small we say that we have
low influence. Consider now a sequence of coefficients $c^{(n)},n\in \N$ and
the corresponding series $S_{N}(c^{(n)},Z).$ In \cite{[MDO]} one proves that,
if $\lim_{n\to\infty}\overline{\delta }_{N}(c^{(n)})=0$ then
\begin{equation*}
\lim_{n\to\infty}\,\sup_{a\in\R}\Delta _{Z,\overline{Z}}(c^{(n)},1_{(a,\infty )})=0
\end{equation*}%
which means that the Kolmogorov distance between $S_{N}(c^{(n)},Z)$ and $%
S_{N}(c^{(n)},\overline{Z})$ converges to zero as $n\to\infty$. Actually the authors of that
paper look to a more particular problem, namely to a single level $\Phi
_{m}(c^{(n)},Z)$ and $\Phi _{m}(c^{(n)},\overline{Z})$, so in this sense our
problem is more general because it concerns series. Moreover, in
\cite{[NPy]} for $\Phi _{m}(c^{(n)},Z)$ and $\Phi _{m}(c^{(n)},\overline{Z})$
as well, under the hypothesis $\lim_{n\to\infty}\overline{\delta }_{N}(c^{(n)})=0$,
one proves convergence in total variation distance that is%
\begin{equation*}
\lim_{n\to\infty}\sup_{\left\Vert f\right\Vert _{\infty }\leq 1}\Delta _{Z,\overline{Z%
}}(c^{(n)},f)=0.
\end{equation*}%
But the authors are obliged to assume more regularity, namely that the law
of $Z_{k}$ is locally lower bounded by the Lebesgue measure: there exist $%
r,\varepsilon >0$ and $z_{k}\in \R$ such that for every measurable set $%
A\subset B_{r}(z_{k})$ one has%
\begin{equation}
\P(Z_{k}\in A)\geq \varepsilon \lambda (A)  \label{Int2}
\end{equation}%
where $\lambda $ is the Lebesgue measure. (\ref{Int2}) is analogous to what is known in the literature as the Doeblin condition. Then they use a splitting method
and the $\Gamma$-calculus settled in \cite{[BGL]} to obtain the
regularity which is needed in order to handle test functions $f$ which are
just measurable. This strategy is close to the method that we use ourselves in
this paper. Notice that the hypothesis (\ref{Int2}) is in fact very mild
(almost necessary): indeed, in the case of the classical CLT (which
corresponds to $N=1$), Prohorov proved in \cite{[PROH]} that in order
to obtain convergence in total variation distance one needs that the law of
the random variables has at least a piece of absolutely continuous component
(and it turns out that this is very close to (\ref{Int2}), see the discussion in Section 2 of \cite{[BC-CLT]}).

Let us now present the contributions of our paper. We first prove that if $%
f\in C_{b}^{3}$ then
\begin{equation}
\Delta _{Z,\overline{Z}}(c,f)\leq C\|f^{\prime \prime \prime }\|_{\infty }\,\overline{\delta }_{N}(c)
\label{Int3}
\end{equation}%
where $C$ is a constant which depends on $M_{p}=M_{p}(Z)=\max_{k}\E(\left%
\vert Z_{k}\right\vert ^{p})$ with $p=3$, see Theorem \ref{Conv-smooth} for
a precise statement. In this case the regularity condition (\ref{Int2})\ is
not required. The proof is a rather standard application of the Lindeberg
method.

We discuss now the convergence to a Gaussian law. First we have to introduce
the ``fourth cumulant'' defined for a random variable $X$ by $\kappa
_{4}(X)=\E(X^{4})-3\E(X^{2})^{2}$. This quantity is known to be a measure of
the distance between the law of $X$ and the standard Gaussian law in the
sense that, if $\lim_{n\to\infty}$ $\kappa _{4}(X_{n})=0$ then $X_{n}\to G$
in law, where $G$ is a standard normal random variable (for which $\kappa
_{4}(G)=0)$. The celebrated ``Fourth Moment Theorem'' proved in \cite{[PN]}
(and then refined in several other papers of Nualart, Nourdin, Peccati and
co-authors, see \texttt{https://sites.google.com/site/malliavinstein/home} for updated references on this subject) asserts that the convergence of the multiple
stochastic integrals $\Phi _{m}(c^{(n)},\overline{Z})$ to the normal law is
equivalent to the convergence of the fourth cumulant. So we define
\begin{equation}
\overline{\kappa }_{N}(c)=\sum_{m=1}^{N}\kappa _{4}^{1/4}(\Phi _{m}(c,%
\overline{Z})), \quad \mbox{with } \overline{Z}_k,k\in\N,\mbox{ i.i.d. standard normal}.  \label{Int4}
\end{equation}%
Notice that, having used independent standard normal random variables, $\overline{\kappa }_{N}(c)$ is in some sense an intrinsic quantity related to the coefficient $c$.

We present now our convergence results. Let $c^{(n)}=(c^{(n)}(\alpha ))_{\alpha
}$ be a sequence of coefficients which verify the normalization condition (%
\ref{Int1a}) and such that for every $N\in \N$%
\begin{equation}
\lim_{n\to\infty}\overline{\kappa }_{N}(c^{(n)})=0
\quad\mbox{and}\quad
\lim_{n\to\infty}\overline{\delta }_{N}(c^{(n)})=0,  \label{Int4b}
\end{equation}
$\overline{\delta}_N(c^{(n)})$ being given in (\ref{Int1b}). We will consider also the following ``uniformity'' assumption:%
\begin{equation}
\limsup_{N\rightarrow \infty }\,\limsup_{n\rightarrow \infty }\sum_{k\geq N}k^{q}\times k!\sum_{\left\vert
\alpha \right\vert =k}\vert c^{(n)}(\alpha )\vert ^{2}=0,\quad q\in\N.
\label{Int4a}
\end{equation}%
We prove that if (\ref{Int4b}) and (\ref{Int4a}) hold, with $q=0,$ then $%
S_{\infty }(c^{(n)},Z)\rightarrow G$ in law as $n\to\infty$ (so, a result for infinite series). In the case of multiple
stochastic integrals (that is when $Z_{k}$ are standard normal) this result has
already been proved in \cite{[HN]}, so, what is new here, is the invariance
principle we are going to introduce in (\ref{Int5}) (see Theorem \ref{CLT} for the precise statement, which needs some
more hypotheses on the moments). Moreover, if (\ref{Int4a}) holds with $q=1$, then
we prove that $S_{\infty }(c^{(n)},Z)\rightarrow G$ in total variation
distance (see Theorem \ref{CTV} for the precise statement). Notice that
\begin{equation*}
\| S_{\infty }(c^{(n)},Z)\|_{2}^{2}=\sum_{k\geq
1}k!\sum_{\left\vert \alpha \right\vert =k}|c^{(n)}(\alpha )| ^{2}
\end{equation*}%
so, in some sense, (\ref{Int4a}) with $q=0$ says that $S_{\infty
}(c^{(n)},Z)$, $n\in \N$, belongs to a ``uniform class'' in $L^{2}.$ And if (\ref%
{Int4a}) holds with $q=1$, one gets a stronger uniformity condition concerning the
Malliavin derivatives - which is morally coherent.

We come back now to our non asymptotic results. The challenging problem now is to
replace $\|f^{\prime \prime \prime }\|_{\infty }$ with $%
\|f\|_{\infty }$ in (\ref{Int3}), so to obtain the
distance in total variation between $S_{N}(c,Z)$ and $S_{N}(c,\overline{Z})$.
In Theorem \ref{TV} we prove that for each $p_{\ast }\geq 1$ one has
\begin{equation}
\Delta _{Z,\overline{Z}}(c,f)\leq C_{*N}\left\Vert
f\right\Vert _{\infty }(\overline{\delta }_{N}(c)+\overline{\kappa }%
_{N}^{p_*}(c))  \label{Int5}
\end{equation}%
where $C_{*N}$ is a constant which depends on $p_{\ast },$ on $N$ and on $%
M_{p}(Z) $ for some $p.$ Here we are obliged to take a finite $N.$

In the papers presented above, that is \cite{[MDO]} and \cite{[NPy]}, the only quantity which was supposed to be small is the low influence factor term $\overline{\delta}_N(c)$. So, the fact that  the 4th cumulant term $\overline{\kappa}_N(c)$ appears in (\ref{Int5}) may be seen as a weak point. However, as long as we deal with convergence to a Gaussian law, we know that, by the Fourth Moment Theorem, we need to ask $\overline{\kappa}_N(c^{(n)})\to 0$ as $n\to\infty$. For a general limit, in \cite{[BCNon]} we will prove that
\begin{equation}
\Delta _{Z,\overline{Z}}(c,f)\leq C_{N}\overline{\delta}_N^{1/N}(c).
\label{Int6}
\end{equation}%
The advantage of (\ref{Int6}) is that $\overline{\kappa}_N(c)$ does not appear in the right hand side, so (\ref{Int6}) works for general limits. But the interest of (\ref{Int5}) is that we get a more accurate estimate (because one has $\overline{\delta}_N(c)$ instead of $\overline{\delta}_N^{1/N}(c)$).

Let us now present our nonlinear CLT. We set $|c|^2_m=\sum_{|\alpha|=m}c(\alpha)^2$ and $\alpha_N(c)=\min_{m\leq N}|c|_m\,1_{|c|_m> 0}$. One may prove (see next (\ref{Law7'})) that $\overline{\delta}_N(c)\leq \alpha_N^{-1}(c)\,\overline{\kappa}_N(c)$, so (\ref{Int5}) with $p_*=1$ reads
\begin{equation}
\Delta _{Z,\overline{Z}}(c,f)\leq C_{*N}\left\Vert
f\right\Vert _{\infty }(1+\alpha_N^{-1}(c))\overline{\kappa}_N(c).
\label{Int5bis}
\end{equation}%
Moreover the Fourth Moment Theorem by Nourdin and Peccati
in \cite{[NP1]} says that if $\overline{Z}_k$, $k\in\N$, are standard normal then
\begin{equation*}
d_{TV}(S_N(c,\overline{Z}),G)\leq C_{N}\overline{\kappa}_N(c).
\end{equation*}%
Therefore, putting things together, in Theorem \ref{G} we prove that%
\begin{equation}
\big|\E(f(S_{N}(c,Z))-\E(f(G))\big| \leq C\| f\|
_{\infty }\,(1+\alpha_N^{-1}(c))\overline{\kappa}_N(c).  \label{Int7}
\end{equation}%

The proof of (\ref{Int5}) is based on integration by parts methodology,
inspired from Malliavin calculus and which has been settled in \cite{[BCl]}%
, \cite{[BC-EJP]} and has already been used in \cite{[BC-CLT]}. As usual the
difficult point which has to be handled is the non degeneracy condition. In
\cite{[NPy]} Nourdin and Poly use the Carbery-Wright inequality for small
balls probabilities in order to solve a similar problem - and we are doing
the same in \cite{[BCNon]}. This approach avoids to use the cumulants $%
\overline{\kappa }_{N}(c)$ but makes appear the power $1/N$ in $\overline{%
\delta }_{N}^{1/N}(c)$. So we give out this approach here and we use an
ad-hoc method based on martingale arguments and Burkholder inequality. A
serious technical difficulty comes from the fact that for general stochastic series
we do not have the product formula which is available for multiple
stochastic integrals (see Lemma \ref{lemma-app1} and Remark \ref{no-product} in Appendix \ref{app-B}).

The paper is organized as follows. In section \ref{notations} we introduce the rather
heavy notation and in Section \ref{conv1} we prove the convergence result for smooth
test functions, that is (\ref{Int3}). In Section \ref{splitting} we recall the variational
calculus that use here and in Section \ref{sect-MallS} we estimate the Sobolev norms of $S_{N}(c,Z)$ and give the non degeneracy estimate. In order to obtain this
last estimate a heavy calculus based on Burkholder inequalities and
a martingale method is needed - we postpone these calculations in
Appendix \ref{burk}. Finally in Section \ref{conv2} we prove the main results, namely (\ref{Int5}) and (\ref{Int7}).

\section{Notations}\label{notations}

The basic objects which appear in this paper are the following.

\begin{itemize}
\item[$\blacksquare$] We denote $\Gamma_m={\mathbb{N}}^m$, that is $\Gamma_m$
is the set of the multi-indexes $\beta =(\beta_{1},...,\beta _{m})$. When $%
m=0$, we define $\Gamma_0$ as the set containing only the null multi-index $%
\beta=\emptyset$. For $\beta\in \Gamma_m$, we say that $\beta$ has length $%
m $, and we define the length as $|\beta|=m$. We set $\Gamma=\cup_{m\geq 0}\Gamma_m$ the set of all multi-indexes. For a fixed $J\in{\mathbb{N}}$,
we set $\Gamma_m(J)$ as the multi-indexes whose components do not exceed $J$%
: $\Gamma_m(J)=\{\beta\in \Gamma_m:\beta_i\leq J\mbox{ for every } i\}$.
Finally we consider the set of ``ordered'' multi-indexes $\Gamma_m^o$ and $%
\Gamma^o_m(J)$: when considering the superscript $^o$ we mean that the multi-index $\beta$ has ordered components, that is $\beta_i<\beta_{i+1}$ for all $%
i$. For $z\in {\mathbb{R}}^{{\mathbb{N}}},z=(z_{k})_{k\in {\mathbb{N}}}$ and
for $\beta \in \Gamma_m$ we denote
\begin{equation}
z^{\beta }=\prod_{k=1}^{m}z_{\beta _{k}}.  \label{N1}
\end{equation}

\item[$\blacksquare$] We consider a sequence of independent (non necessarily identically distributed) random variables
$Z_{k},k\in {\mathbb{N}}$ which, for some $p\geq 1,$ verify
\begin{equation}  \label{Law1}
\mbox{$\E(Z_{k})=0$ and  $\E(Z_{k}^{2})=1\,\forall k$, $M_{p}(Z):=\max_{k}\E(|Z_{k}|^{p})<\infty$  and we set
$\overline{M}_{p}(Z)=b_{p}M_{p}(Z)$}
\end{equation}
where $b_{p}$ is the constant in the Burkholder inequality of order $p$ (see
Appendix \ref{burk}).

\item[$\blacksquare $] We consider a family of coefficients $c(\alpha )$, $%
\alpha \in \Gamma $, which are symmetric and null on all the diagonals: if $\alpha =(\alpha _{1},...,\alpha _{m})\in \Gamma _{m}$, then
for every permutation $\pi $ of $\{1,...,m\}$ one has $c(\alpha )=c(\alpha
_{\pi })$ with $\alpha _{\pi }=(\alpha _{\pi _{1}},...,\alpha _{\pi _{m}})$;
if $\alpha _{i}=\alpha _{j}$ for some $i\neq j$ then $c(\alpha )=0$.
\end{itemize}

We set

\begin{equation}
\left\vert c\right\vert _{m}=\Big(\sum_{\alpha \in \Gamma _{m}}c^{2}(\alpha )%
\Big)^{1/2},\qquad \left\Vert c\right\Vert _{m}=(\sum_{i=1}^{m}\left\vert
c\right\vert _{i}^{2})^{1/2},  \label{Law2}
\end{equation}%
and, for $q\in {\mathbb{N}}$ and $M>0$,
\begin{equation}
N_{q}(c,M)=\Big(\sum_{m=q}^{\infty }M^{m-q}\times \frac{m!}{(m-q)!}\times
m!\left\vert c\right\vert _{m}^{2}\Big)^{1/2}.  \label{Law3}
\end{equation}%
Moreover we denote
\begin{equation}
\delta _{1}(c)=\sup_{k}|c(k)|\quad \mbox{and}\quad \delta _{m}(c)=\sup_{k}%
\Big(\sum_{\alpha \in \Gamma _{m-1}}c^{2}(k,\alpha )\Big)^{1/2},\quad m\geq
2.  \label{Law3''}
\end{equation}%
We use the notation $(k,\alpha )=(k,\alpha _{1},...,\alpha _{m-1})$ for
$\alpha =(\alpha _{1},...,\alpha _{m-1})\in \Gamma _{m-1}$ (note that if $%
m=1 $ then $\Gamma _{m-1}$ contains only the void multi-index and $%
c^{2}(k,\emptyset )=c^{2}(k)$). Roughly speaking $\delta _{m}(c)$ quantifies
the maximum action of a single particle on the other ones, and, if $\delta
_{m}(c)$ is required to be small, we say that we have a \textquotedblleft
low influence\textquotedblright\ condition. Moreover we denote%
\begin{equation}
\varepsilon _{0}(c,M)=\sup_{k}\Big(\sum_{m=1}^{\infty }M^{2m}\times
m!(\sum_{\alpha \in \Gamma _{m}}c^{2}(k,\alpha ))\Big)^{1/2}\leq
\sum_{m=0}^{\infty }M^{2m}\times m!\times \delta _{m+1}(c).  \label{Law3'}
\end{equation}

Given $Z$ and $c$ as above we define
\begin{equation}
\Phi _{m}(c,Z)=\sum_{\alpha \in \Gamma _{m}}c(\alpha )Z^{\alpha }\quad %
\mbox{and}\quad \Phi _{m}^{o}(c,Z)=\sum_{\alpha \in \Gamma _{m}^{o}}c(\alpha
)Z^{\alpha }  \label{Law4}
\end{equation}%
Since $c$ is symmetric we have $\Phi _{m}(c,Z)=m!\Phi _{m}^{o}(c,Z).$

\begin{remark}\label{iterated}
We notice that $\Phi _{m}(c,Z)$ is (in law) a multiple stochastic integral when the r.v.'s $Z_k$, $k\in\N$, are i.i.d. standard normal. In fact, $W$ denoting a Brownian motion in $\R$, one has
$$
\Phi _{m}^{o}(c,Z)=\int_0^\infty\int_0^{t_m}\cdots\int_0^{t_{2}}f(t_1,\ldots,t_m)dW_{t_1}\cdots dW_{t_{m-1}}dW_{t_m}
$$
when we take $Z_k=W_{k+1}-W_{k}$, $k\in\N$ and
$$
f(t_1,\ldots,t_m)=\sum_{\alpha\in \Gamma_m}c(\alpha)\prod_{i=1}^m\I_{[\alpha_i,\alpha_i+1)}(t_i).
$$
\end{remark}

We
finally set
\begin{equation}
S(c,Z)=\sum_{m=1}^{\infty }\Phi _{m}(c,Z).  \label{Law7}
\end{equation}
For finite sums, we will use the notation
\begin{equation}
S_N(c,Z)=\sum_{m=1}^{N }\Phi _{m}(c,Z).  \label{Law7-bis}
\end{equation}

For a random variable $X$ we denote by $\kappa _{4}(X)$ the fourth cumulant,
that is
\begin{equation}
\kappa _{4}(X)={\mathbb{E}}(X^{4})-3{\mathbb{E}}(X^{2})^{2}.  \label{Law6''}
\end{equation}%
We will use the notation
\begin{equation}
\kappa _{4,m}(c)=\kappa _{4}(\Phi _{m}(c,\overline{Z}))\quad
\mbox{when
$\overline{Z}_1,\overline{Z}_2,\ldots$ are i.i.d. $\sim\mathcal{N}(0,1)$,}
\label{Law6'''}
\end{equation}%
where, from now on, $\mathcal{N}(\mu,\sigma^2)$ denotes the normal law with mean $\nu$ and variance $\sigma^2$. Recall that $%
\kappa _{4,1}(c)=0$, because $\Phi _{1}(c,\overline{Z})$ is centered
Gaussian. We also denote%
\begin{equation}
\overline{\kappa }_{4,N}(c)=\sum_{m=1}^{N}\kappa _{4,m}^{1/4}(c).
\label{Law6v}
\end{equation}%
Moreover it is known (see \cite{[NPRein]} or (\ref{acc10})) that%
\begin{equation}
\delta _{m}(c)\leq \frac{1}{m!m\left\vert c\right\vert _{m}}\sqrt{\kappa
_{4,m}(c)},\quad m\geq 2.  \label{Law7'}
\end{equation}%
In particular%
\begin{equation}
\varepsilon _{0}(c,M)\leq \Big(\sum_{m=1}^{\infty }\frac{M^{2m}}{m!m^{2}}%
\frac{\kappa _{4,m}(c)}{\left\vert c\right\vert _{m}^{2}}\Big)^{1/2}.
\label{Law7''}
\end{equation}

\section{Convergence of series for smooth test functions}\label{conv1}

The aim of this section is to discuss the convergence in law of the series
of the form (\ref{Law7}). The following estimates are immediate consequences
of the isometry property and of Burkholder inequality (but see also next
Lemma \ref{ES}):
\begin{align}
\left\Vert S(c,Z)\right\Vert _{2}& =N_{0}(c,1)=(\sum_{m=0}^{\infty
}m!\left\vert c\right\vert _{m}^{2})^{1/2}\qquad \mbox{and}  \label{Law8} \\
\left\Vert S(c,Z)\right\Vert _{p}& \leq N_{0}(c,\overline{M}%
_{p}^{2}(Z))=(\sum_{m=0}^{\infty }\overline{M}_{p}^{2m}(Z)\times
m!\left\vert c\right\vert _{m}^{2})^{1/2}.  \label{Law9}
\end{align}

Our first result consists in comparing $S(c,Z)$ and $S(c,\overline{Z})$ for
two different sequences $Z_{j},\overline{Z}_{j},j\in {\mathbb{N}}$ of random
variables.

\begin{theorem}
\label{Conv-smooth} Let $Z=(Z_{k})_{k\in N}$ and $\overline{Z}=(\overline{Z}%
_{k})_{k\in N}$ be two sequences of independent random variables such that ${%
\mathbb{E}}(Z_{k})={\mathbb{E}}(\overline{Z}_{k})=0$ and ${\mathbb{E}}%
(Z_{k}^{2})={\mathbb{E}}(\overline{Z}_{k}^{2})=1$. Recall (\ref{Law1}), (\ref{Law3}), (\ref{Law3'}) for the definition of $\overline{M}_{p}(Z)$ and $\overline{M}_{p}(\overline{Z})$, $N_{0}(c,M)$, $\varepsilon _{0}(c,M)$ respectively.
\begin{itemize}
\item[\textbf{A}.]
Let $M_{3}=\overline{M}_{3}^{2}(Z)\vee \overline{M}_{3}^{2}(%
\overline{Z})<\infty$. Then for every $f\in C_{b}^{3}({\mathbb{R}})$,
\begin{equation}
\left\vert {\mathbb{E}}(f(S(c,Z)))-{\mathbb{E}}(f(S(c,\overline{Z}%
))\right\vert \leq \frac{1}{3}\| f^{(3)}\| _{\infty
}M_{3}^{3}N_{0}(c,M_{3})^2\varepsilon _{0}(c,M_{3}).  \label{Law10}
\end{equation}
\item[\textbf{B}.]
Suppose moreover that ${\mathbb{E}}(Z_{k}^{3})={\mathbb{E}}(%
\overline{Z}_{k}^{3})=0$ and $M_{4}=\overline{M}_{4}^{2}(Z)\vee \overline{M}%
_{4}^{2}(\overline{Z})<\infty $. Then for every $f\in C_{b}^{4}({\mathbb{R}}%
) $%
\begin{equation}
\left\vert {\mathbb{E}}(f(S(c,Z)))-{\mathbb{E}}(f(S(c,\overline{Z}%
))\right\vert \leq \frac{1}{12}\Vert f^{(4)}\Vert _{\infty
}M_{4}^{4}N_{0}(c,M_{4})^{2}\varepsilon _{0}^{2}(c,M_{4}).  \label{Law10a}
\end{equation}%
\end{itemize}
\end{theorem}

\textbf{Proof}. \textbf{A}. Let $m,J\in {\mathbb{N}}$ and%
\begin{equation*}
S_{m,J}(c,Z)=\sum_{n=1}^{m}\sum_{\alpha \in \Gamma _{n}(J)}c(\alpha
)Z^{\alpha }.
\end{equation*}%
We will prove (\ref{Law10}) with $S_{m,J}(c,Z)$ instead of $S(c,Z).$ Since
the upper bound in the right hand side will not depend on $m$ and $J,$ the
inequality for $S(c,Z)$ is obtained by passing to the limit with $%
m,J\rightarrow \infty .$

For $k=1,...,J$ and $\theta \in \lbrack 0,1]$\ we define the vector $%
\widehat{Z}^{k}(\theta )$ by%
\begin{equation*}
\widehat{Z}^{k}(\theta )=(Z_{1},...,Z_{k-1},\theta Z_{k},\overline{Z}%
_{k+1},...,\overline{Z}_{J}).
\end{equation*}%
Then%
\begin{align*}
& {\mathbb{E}}(f(S_{m,J}(c,Z)))-{\mathbb{E}}(f(S_{m,J}(c,\overline{Z})) \\
& =\sum_{k=1}^{J}\big[{\mathbb{E}}(f(S_{m,J}(c,\widehat{Z}^{k}(1))))-{%
\mathbb{E}}(f(S_{m,J}(c,\widehat{Z}^{k}(0))))\big] \\
& \quad -\sum_{k=1}^{J}\big[{\mathbb{E}}(f(S_{m,J}(c,\widehat{Z}^{k}(0))))-{%
\mathbb{E}}(f(S_{m,J}(c,\widehat{Z}^{k-1}(1))))\big].
\end{align*}%
We write
\begin{equation*}
S_{m,J}(c,\widehat{Z}^{k}(1))=s_{k}+Z_{k}v_{k}
\end{equation*}%
with%
\begin{equation*}
s_{k}=S_{m,J}(c,\widehat{Z}^{k}(0))\qquad \mbox{and}\qquad
v_{k}=c(k)+\sum_{n=2}^{m}\sum_{\substack{ \beta \in \Gamma _{n-1}(J)  \\ %
k\notin \beta }}c(k,\beta )(\widehat{Z}^{k}(0))^{\beta }
\end{equation*}%
(recall that $(k,\beta )=(k,\beta _{1},...,\beta _{m-1})$ for $\beta =(\beta
_{1},...,\beta _{m-1})$). Then, by using a development in Taylor series of
order three,
\begin{align*}
f(S_{m,J}(c,\widehat{Z}^{k+1}(1))))-f(S_{m,J}(c,\widehat{Z}^{k+1}(0))&
=f^{\prime }(s_{k})Z_{k}v_{k}+\frac{1}{2}f^{\prime \prime
}(s_{k})Z_{k}^{2}v_{k}^{2}+ \\
& \quad +\frac{1}{6}\int_{0}^{1}f^{\prime \prime \prime }(s_{k}+\theta
Z_{k}v_{k})Z_{k}^{3}v_{k}^{3}\theta d\theta .
\end{align*}%
By taking expectation, by using independence and ${\mathbb{E}}(Z_{k})=0$, ${%
\mathbb{E}}(Z_{k}^{2})=1$, we obtain%
\begin{align*}
{\mathbb{E}}(f(S_{m,J}(c,\widehat{Z}^{k}(1)))-{\mathbb{E}}(f(S_{m,J}(c,%
\widehat{Z}^{k}(0)))& =\frac{1}{2}{\mathbb{E}}(f^{\prime \prime
}(s_{k})v_{k}^{2})+ \\
& \quad +\frac{1}{6}\int_{0}^{1}{\mathbb{E}}(f^{\prime \prime \prime
}(s_{k}+\theta Z_{k}v_{k})Z_{k}^{3}v_{k}^{3})\theta d\theta .
\end{align*}%
In a similar way we get%
\begin{align*}
{\mathbb{E}}(f(S_{m,J}(c,\widehat{Z}^{k}(0))))-{\mathbb{E}}(f(S_{m,J}(c,%
\widehat{Z}^{k-1}(1)))& =-\frac{1}{2}{\mathbb{E}}(f^{\prime \prime
}(s_{k})v_{k}^{2})+ \\
& \quad -\frac{1}{6}\int_{0}^{1}{\mathbb{E}}(f^{\prime \prime \prime
}(s_{k}+\theta \overline{Z}_{k}v_{k})\overline{Z}_{k}^{3}v_{k}^{3})\theta
d\theta .
\end{align*}%
The term containing $f^{\prime \prime }$ is the same in the two cases, so it
cancels when taking sums. Then we obtain%
\begin{align*}
{\mathbb{E}}(f(S_{m,J}(c,Z)))-{\mathbb{E}}(f(S_{m,J}(c,\overline{Z}))& =%
\frac{1}{6}\sum_{k=1}^{J}\int_{0}^{1}{\mathbb{E}}(f^{\prime \prime \prime
}(s_{k}+\theta Z_{k}v_{k})Z_{k}^{3}v_{k}^{3}))\theta d\theta \\
& \quad -\frac{1}{6}\sum_{k=1}^{J}\int_{0}^{1}{\mathbb{E}}(f^{\prime \prime
\prime }(s_{k}+\theta \overline{Z}_{k}v_{k})\overline{Z}_{k}^{3}v_{k}^{3}))%
\theta d\theta .
\end{align*}%
Since $Z_{k}$ is independent of $v_{k}$ we have
\begin{equation*}
\left\vert {\mathbb{E}}(f^{\prime \prime \prime }(s_{k}+\theta
Z_{k}v_{k})Z_{k}^{3}v_{k}^{3}))\right\vert \leq \left\Vert f^{\prime \prime
\prime }\right\Vert _{\infty }{\mathbb{E}}(\left\vert Z_{k}\right\vert ^{3}){%
\mathbb{E}}(\left\vert v_{k}\right\vert ^{3})\leq \left\Vert f^{\prime
\prime \prime }\right\Vert _{\infty }M_{3}^{3}(Z){\mathbb{E}}(\left\vert
v_{k}\right\vert ^{3}).
\end{equation*}%
Now, using (\ref{Law9}) with $p=3,$ we obtain (recall that $M_{3}=\overline{M%
}_{3}^{2}(Z)\vee \overline{M}_{3}^{2}(\overline{Z}))$%
\begin{equation*}
{\mathbb{E}}(\left\vert v_{k}\right\vert ^{3})\leq
(\sum_{n=0}^{m-1}M_{3}^{n}\times n!\sum_{\substack{ \left\vert \beta
\right\vert =n  \\ k\notin \beta }}c^{2}(k,\beta ))^{3/2}
\end{equation*}%
and so%
\begin{equation*}
\left\vert {\mathbb{E}}(f^{\prime \prime \prime }(s_{k}+\theta
Z_{k}v_{k})Z_{k}^{3}v_{k}^{3}))\right\vert \leq \left\Vert f^{\prime \prime
\prime }\right\Vert _{\infty }M_{3}^{3}(\sum_{n=0}^{m-1}M_{3}^{n}\times
n!\sum_{\substack{ \left\vert \beta \right\vert =n  \\ k\notin \beta }}%
c^{2}(k,\beta ))^{3/2}.
\end{equation*}%
The same estimate holds if we take $\overline{Z}_{k}$ instead of $Z_{k}.$ We
conclude that%
\begin{align*}
\left\vert {\mathbb{E}}(f(S_{m,J}(c,Z)))-{\mathbb{E}}(f(S_{m,J}(c,\overline{Z%
}))\right\vert \leq & \frac{1}{3}\left\Vert f^{\prime \prime \prime
}\right\Vert _{\infty
}M_{3}^{3}\sum_{k=1}^{J}(\sum_{n=0}^{m-1}M_{3}^{n}\times n!\sum_{\left\vert
\beta \right\vert =n}c^{2}(k,\beta ))^{3/2} \\
\leq & \frac{1}{3}\left\Vert f^{\prime \prime \prime }\right\Vert _{\infty
}M_{3}^{3}\max_{k\leq J}(\sum_{n=0}^{m-1}M_{3}^{n}\times n!\sum_{\left\vert
\beta \right\vert =n}c^{2}(k,\beta ))^{1/2}\times \\
& \times \sum_{k=1}^{J}\sum_{n=0}^{m-1}M_{3}^{n}\times n!\sum_{\left\vert
\beta \right\vert =n}c^{2}(k,\beta ).
\end{align*}%
Notice that
\begin{equation*}
\sum_{k=1}^{J}\sum_{n=0}^{m-1}M_{3}^{n}\times n!\sum_{\left\vert \beta
\right\vert =n}c^{2}(k,\beta )=\sum_{n=0}^{m-1}M_{3}^{n}\times
n!\sum_{\left\vert \alpha \right\vert =n+1}c^{2}(\alpha )\leq
N_{0}(c,M_{3})^{2}
\end{equation*}%
and%
\begin{equation*}
\max_{k\leq J}(\sum_{n=0}^{m-1}M^{n}\times n!\sum_{\left\vert \beta
\right\vert =n}c^{2}(k,\beta ))^{1/2}\leq \varepsilon _{0}(c,M).
\end{equation*}%
We conclude that%
\begin{equation*}
\left\vert {\mathbb{E}}(f(S_{m,J}(c,Z)))-{\mathbb{E}}(f(S_{m,J}(c,\overline{Z%
}))\right\vert \leq \frac{1}{3}M_{3}^{3}\left\Vert f^{\prime \prime \prime
}\right\Vert _{\infty }\varepsilon _{0}(c,M_{3})N_{0}(c,M_{3})^{2}
\end{equation*}%
so the proof of (\ref{Law10})\ is completed. The proof of (\ref{Law10a}) is
identical: one just go further to order 4 in the Taylor expansion.
$\square $ \medskip

We discuss now the convergence to a Gaussian random variable. This
immediately follows from the previous result and from Theorem 3 in \cite%
{[HN]}. We consider a sequence $c^{(n)}=(c^{(n)}(\alpha ))_{\alpha \in
\Gamma },n\in \N$ of coefficients and the corresponding stochastic series $%
S(c^{(n)},Z).$ Our assumptions will be the following:%
\begin{equation}\label{Law11}
\begin{array}{l}
\displaystyle
i)\quad \limsup_{N\rightarrow \infty }\limsup_{n\rightarrow
\infty }\sum_{k\geq N}k!\sum_{\left\vert \alpha \right\vert =k}|c^{(n)}(\alpha )| ^{2} =0,   \smallskip\\
\displaystyle
ii)\quad \lim_{n\rightarrow \infty }k!\sum_{\left\vert \alpha \right\vert
=k}|c^{(n)}(\alpha )|^{2} =:\sigma _{k}^{2},  \smallskip\\
\displaystyle
iii)\quad \sum_{k=1}^{\infty }\sigma _{k}^{2} =\sigma ^{2}  \smallskip\\
\displaystyle
iv)\quad \lim_{n\rightarrow \infty }\kappa _{4,k}(c^{(n)}) =0,\quad
\forall k\in \N.
\end{array}
\end{equation}

\begin{theorem}
\label{CLT} Let $Z=(Z_{k})_{k\in N}$ be a sequence of independent random
variables such that ${\mathbb{E}}(Z_{k})=0,$ ${\mathbb{E}}(Z_{k}^{2})=1$ and
$M_{3}=\max_{k\in \N}\E(\left\vert Z_{k}\right\vert ^{3})<\infty .$ Consider
also a sequence $c^{(n)}=(c^{(n)}(\alpha ))_{\alpha \in \Gamma },n\in \N$ of
coefficients which satisfy (\ref{Law11}). Then $\lim_{n\rightarrow \infty
}S(c^{(n)},Z)=\mathcal{N}(0,\sigma^2 )$ in law.
\end{theorem}

\textbf{Proof}. Notice first that the sequence of laws of $S(c^{(n)},Z)$ is
tight, so the only thing to be proven is that any limit point is in fact $%
\mathcal{N}(0,\sigma^2 )$. Let $\overline{Z}_{k},k\in \N$ be a sequence of independent and standard normal random variables, so that $S(c^{(n)},\overline{Z})$ is an infinite
sum of multiple stochastic integrals (see Remark \ref{iterated}). Then, under the hypothesis (\ref{Law11}%
), Hu and Nualart proved (see Theorem 3 in \cite{[HN]}) that $%
\lim_{n\rightarrow \infty }S(c^{(n)},\overline{Z})=\mathcal{N}(0,\sigma^2 )$ in law. And
(\ref{Law10}) guarantees that $\lim_{n\rightarrow \infty }S(c^{(n)},Z)$ is
the same. $\square $

\section{Variational calculus using a splitting method}\label{splitting}

In order to study the convergence in total variation and some related invariance
principles, our specific point is to consider a class of random
variables which have a regularity property allowing one to extrapolate an ``absolutely continuous noise''.

\subsection{The splitting procedure}
We say that the law
of the random variable $Z\in {\mathbb{R}}$ is locally lower bounded by the
Lebesgue measure if there exists $z\in {\mathbb{R}}$ and $\varepsilon
,r>0$ such that for every non negative and measurable function $f:{\mathbb{R}%
}\rightarrow {\mathbb{R}}_{+}$
\begin{equation}
\begin{array}{ll}
\hskip-3cmA(z, r,\varepsilon ):\qquad \qquad & {\mathbb{E}}(f(Z))\geq
\varepsilon \int f(\xi-z)1_{B(0,r)}(\xi-z)d\xi.%
\end{array}
\label{I1}
\end{equation}%
We denote by $\mathcal{L}(z,r,\varepsilon )$ the class of the random variables
which verify $A(z,r,\varepsilon )$. Given $r>0$\ we consider the functions $%
\theta _{r},\psi _{r}:{\mathbb{R}}\rightarrow {\mathbb{R}}_{+}$ defined by
\begin{equation}
\theta _{r}(t)=1-\frac{1}{1-(\frac{t}{r}-1)^{2}}\qquad \psi
_{r}(t)=1_{\{\left\vert t\right\vert \leq r\}}+1_{\{r<\left\vert
t\right\vert \leq 2r\}}e^{\theta _{r}(\left\vert t\right\vert )}.  \label{N4}
\end{equation}%
If $Z\in \mathcal{L}(z,r,\varepsilon )$ then%
\begin{equation}
{\mathbb{E}}(f(Z))\geq \varepsilon \int f(\xi-z)\psi _{r}(\left\vert
\xi-z\right\vert ^{2})d\xi.  \label{I3}
\end{equation}%
The advantage of $\psi _{r}(\left\vert \xi-z\right\vert ^{2})$ is that it
is a smooth function (which replaces the indicator function of the ball) and
(it is easy to check) that for each $l\in {\mathbb{N}},p\geq 1$ there exists
a constant $C_{l,p}\geq 1$ such that
\begin{equation}
\sup_{t\in\R}\psi _{r}(t)|\theta _{r}^{(l)}(|t|)|^{p}\leq \frac{C_{l,p}}{r^{lp}}
\label{N5}
\end{equation}%
where $\theta _{r}^{(l)}$ denotes the derivative of order $l$ of $\theta
_{r}.$ Moreover, in Proposition 3.1 in  \cite{[BC-CLT]} it is proved that  if $Z\in \mathcal{L}%
(z,r,\varepsilon )$ then $Z$ admits the following decomposition (the equality
is understood as identity of laws):
\begin{equation}
Z=\chi V+(1-\chi )U  \label{I2}
\end{equation}%
where $\chi ,U,V$ are independent random variables with the following laws:
\begin{equation}
\begin{array}{c}
\displaystyle{\mathbb{P}}(\chi =1)=\varepsilon m(r)\quad \mbox{and}\quad {%
\mathbb{P}}(\chi =0)=1-\varepsilon m(r),\smallskip \\
\displaystyle{\mathbb{P}}(V\in d\xi)=\frac{1}{m(r)}\psi _{r}(\left\vert
\xi-z)\right\vert ^{2}d\xi\smallskip \\
\displaystyle{\mathbb{P}}(U\in d\xi)=\frac{1}{1-\varepsilon m(r)}({\mathbb{P}}%
(Z\in d\xi)-\varepsilon \psi _{r}(\left\vert \xi-z\right\vert ^{2})d\xi)%
\end{array}
\label{I4}
\end{equation}%
with%
\begin{equation}
m(r)=\int \psi _{r}(|\xi|^{2})d\xi.  \label{N6}
\end{equation}

\begin{assumption}
From now on, we consider functionals of a sequence of independent random
variables $Z_{k}\in {\mathbb{R}},k\in {\mathbb{N}}$, having all moments and
such that $Z_{k}\in \mathcal{L}(z_k,r,\varepsilon )$ for every $k\in {\mathbb{N}}$.
Remark that $Z_{k}$ are not identically distributed but we assume that $r$ and $\varepsilon$ are the same for all of them (on the contrary, $z_{k}$ may depend on $k$).
\end{assumption}

\subsection{Differential operators and Sobolev spaces}

\label{sect-Malliavin}

We use the stochastic differential calculus (an abstract finite dimensional
Malliavin type calculus) based on $V_{k},k\in {\mathbb{N}}$ settled in \cite%
{[BCl]} \cite{[BC-EJP]} and, for this kind of splitting, in \cite{[BC-CLT]}.
The crucial point is that the law of $V_{k}$ is absolutely continuous and
has the nice density $\psi _{r}(\left\vert z-z_{k})\right\vert ^{2}).$ We
recall the results we need in the following sections.

\smallskip

We denote by $\mathcal{P}$ the subspace of the
measurable functions $\Phi :{\mathbb{R}}^{{\mathbb{N}}}\rightarrow {\mathbb{R%
}}$ that are polynomials. So $\Phi \in \mathcal{P}$ with degree $n$ means
that there exists $n\in {\mathbb{N}}$ and $(c(\beta ))_{\beta \in \cup
_{m=1}^{n}\Gamma _{m}}$ such that
\begin{equation}
\Phi (z)=\sum_{n=1}^{m}\sum_{\beta \in \Gamma _{n}}c(\beta )z^{\beta }.
\label{N2}
\end{equation}

We define the space of simple functionals%
\begin{equation}
\mathcal{S}=\{F=\Phi (Z):\Phi \in \mathcal{P}\}  \label{S1}
\end{equation}%
where $\mathcal{P}$\ is the space of the polynomials defined above. For $%
F=\Phi (Z)\in \mathcal{S}$ we define the derivative operator%
\begin{equation}
D_{k}F=\frac{\partial }{\partial V_{k}}F=\chi _{k}\frac{\partial }{\partial Z_{k}}F=\chi _{k}\partial_k\Phi(Z).
\label{S2}
\end{equation}%
We look to $DF=(D_{k}F)_{k\in {\mathbb{N}}}$ as to a random element of the
Hilbert space
\begin{equation}  \label{H}
\mathcal{H}=\Big\{z\in {\mathbb{R}}^{{\mathbb{N}}}:\left\vert z\right\vert _{%
\mathcal{H}}^{2}:=\sum_{k=1}^{\infty }z_{k}^{2}<\infty \Big\}.
\end{equation}%
Moreover we define the higher order derivatives in the following way. Let $%
n\in {\mathbb{N}}$ be fixed and let $\alpha =(\alpha _{1},...,\alpha _{n}).$
For $F=\Phi(Z)\in\mathcal{S}$, we define
\begin{equation}  \label{S4}
D_{\alpha }^{(n)}F =D_{\alpha_n}\cdots D_{\alpha_1}F =\Big(\prod_{j=1}^{n}\chi
_{\alpha_{j}}\Big)(\partial_{\alpha_n}\cdots\partial_{\alpha_1}\Phi)(Z) =%
\Big(\prod_{j=1}^{n}\chi _{\alpha _{j}}\Big)\partial_\alpha\Phi(Z).
\end{equation}%
We look to $D^{(n)}F=(D_{\alpha }^{(n)}F)_{\alpha \in \Gamma_{n}}$ as to a
random element of $\mathcal{H}^{\otimes n}$. For $n=1$, we write $%
D^{(1)}F=DF $.

We define now%
\begin{equation}
LF=-\sum_{k=1}^{\infty }(D_{k}D_{k}F+D_{k}F\times \Theta _{k})\qquad %
\mbox{with}\qquad \Theta _{k}=2\theta _{r}^{\prime }(\left\vert
V_{k}-z_{k}\right\vert ^{2})(V_{k}-z_{k}).  \label{S6}
\end{equation}%
Elementary integration by parts gives the following duality relation: for
every $F,G\in \mathcal{P}$%
\begin{equation}
{\mathbb{E}}(\left\langle DF,DG\right\rangle _{\mathcal{H}})={\mathbb{E}}%
(FLG)={\mathbb{E}}(GLF).  \label{S7}
\end{equation}

We define now the Sobolev norms. For $q\geq 1$ we set
\begin{equation}
\left\vert F\right\vert _{1,q}=\sum_{n=1}^{q}| D^{(n)}F| _{\mathcal{H}%
^{\otimes n}}\quad \mbox{and}\quad \left\vert F\right\vert _{q}=\left\vert
F\right\vert +\left\vert F\right\vert _{1,q}.  \label{S8}
\end{equation}
Moreover we define%
\begin{equation}
\left\Vert F\right\Vert _{1,q,p}=\big({\mathbb{E}}(\left\vert F\right\vert
_{1,q}^{p})\big)^{1/p},\qquad \left\Vert F\right\Vert _{q,p}=\big({\mathbb{E}%
}(\left\vert F\right\vert _{q}^{p})\big)^{1/p}  \label{S9}
\end{equation}%
and%
\begin{equation}
\left\Vert \left\vert F\right\vert \right\Vert _{1,q,p}=\left\Vert
F\right\Vert _{1,q,p}+\left\Vert LF\right\Vert _{q-2,p},\qquad \left\Vert
\left\vert F\right\vert \right\Vert _{q,p}=\left\Vert F\right\Vert
_{p}+\left\Vert \left\vert F\right\vert \right\Vert _{1,q,p}.  \label{S3}
\end{equation}

Finally we define the Sobolev spaces%
\begin{equation}
{\mathbb{D}}^{q,p}=\overline{\mathcal{S}}^{\|\cdot \|_{q,p}},\qquad {\mathbb{%
D}}^{q,\infty }=\cap _{p=1}^{\infty }{\mathbb{D}}^{q,p}\qquad {\mathbb{D}}%
^{\infty }=\cap _{q=1}^{\infty }{\mathbb{D}}^{q,\infty }.  \label{S10}
\end{equation}%
Notice that the duality relation (\ref{S7}) implies that the operators $%
D^{(n)}$ and $L$ are closable so we may extend these operators to ${\mathbb{D%
}}^{q,p}$ in a standard way.

\subsection{Integration by parts formula}

In this section we recall some results from \cite{[BC-EJP]} and \cite{[BCl]}%
. All these results are stated in that papers for a functional $F$ which
depends on $Z_{1},...,Z_{J}$ only (a finite number of random variables). But
all of them extend in a trivial way to $F\in {\mathbb{D}}^{q,p}.$

We recall first the basic computational rules and the integration by parts
formula. For $\phi \in C^{1}({\mathbb{R}}^{M})\in ({\mathbb{D}}^{2,\infty })^M$ we have%
\begin{equation}
D\phi (F)=\sum_{j=1}^{M}\partial _{j}\phi (F)DF^{j},  \label{FD9}
\end{equation}%
and for $\phi \in C^{2}({\mathbb{R}}^{M})$
\begin{equation}
L\phi (F)=\sum_{j=1}^{M}\partial _{j}\phi (F)LF^{j}-\frac{1}{2}%
\sum_{i,j=1}^{M}\partial _{i}\partial _{j}\phi (F)\left\langle
DF^{i},DF^{j}\right\rangle_{\mathcal{H}} .  \label{FD11}
\end{equation}%
In particular for $F,G\in {\mathbb{D}}^{2,\infty }$
\begin{equation}
L(FG)=FLG+GLF-\left\langle DF,DG\right\rangle_{\mathcal{H}} .  \label{FD12}
\end{equation}

For a functional $F=(F^{1},...,F^{M})\in({\mathbb{D}}^{2,\infty})^M$
we define the Malliavin covariance matrix $\sigma _{F}$ by%
\begin{equation}
\sigma _{F}^{i,j}=\left\langle DF^{i},DF^{j}\right\rangle_{\mathcal{H}} ,\qquad
i,j=1,...,M.  \label{FD1}
\end{equation}%
The lower eigenvalue of $\sigma _{F}$ is
\begin{equation}
\lambda _{F}=\inf_{\left\vert \xi \right\vert =1}\left\langle \sigma _{F}\xi
,\xi \right\rangle =\inf_{\left\vert \xi \right\vert =1}\sum_{k=1}^{\infty
}\left\langle D_{(k)}F,\xi \right\rangle_{\mathcal{H}} ^{2}.  \label{FD2}
\end{equation}%
If $\sigma _{F}$ is invertible we denote $\gamma _{F}=\sigma _{F}^{-1}.$
Moreover we denote
\begin{equation}
\overline{\sigma }_{F}(p)=1\vee {\mathbb{E}}((\det \sigma _{F})^{-p}),\qquad
\overline{\lambda }_{F}(p)=1\vee {\mathbb{E}}((\lambda _{F})^{-p})
\label{FD3}
\end{equation}

We are now able to give the Malliavin integration by parts formulae. Here, $%
C_{p}^{\infty}({\mathbb{R}}^{M})$ denotes the set of the infinitely
differentiable functions whose derivatives, of any order, have polynomial
growth.

\begin{theorem}
\label{TH1}Let $F=(F^{1},...,F^{M}),F_{i}\in {\mathbb{D}}^{2,\infty }$ and $%
G\in {\mathbb{D}}^{1,\infty }$ be such that $\overline{\sigma }%
_{F}(p)<\infty $ for every $p\geq 1.$ Then for every $\phi \in C_{p}^{\infty
}({\mathbb{R}}^{M})$ and every $i=1,...,M$%
\begin{equation}
{\mathbb{E}}(\partial _{i}\phi (F)G)={\mathbb{E}}(\phi (F)H_{i}(F,G))
\label{FD13}
\end{equation}%
with
\begin{equation}
H_{i}(F,G)=G\gamma _{F}LF+\left\langle D(G\gamma _{F}),DF\right\rangle_{\mathcal{H}}
\label{FD14}
\end{equation}%
Moreover let $m\in {\mathbb{N}},m\geq 2$ and $\alpha =(\alpha
_{1},...,\alpha _{m})\in \{1,...,M\}^{m}.$ Suppose that $%
F=(F^{1},...,F^{M}),F_{i}\in {\mathbb{D}}^{m+1,\infty }$ and $G\in {\mathbb{D%
}}^{m,\infty }.$ Then%
\begin{equation}
{\mathbb{E}}(\partial _{\alpha }\phi (F)G)=\E(\phi (F)H_{\alpha }(F,G))
\label{FD15}
\end{equation}%
with $H_{\alpha }(F,G)$ defined by $H_{(\alpha _{1},...,\alpha
_{m})}(F,G):=H_{\alpha _{m}}(F,H_{(\alpha _{1},...,\alpha _{m-1})}(F,G)).$
\end{theorem}

\textbf{Proof. } We give here only a sketch of the proof, a detailed one can
be found e.g. in \cite{[BC-EJP]} and \cite{[BCl]}. Using the chain rule $%
D\phi (F)=\nabla \phi (F)DF$ so that
\begin{equation*}
\left\langle D\phi (F),DF\right\rangle_{\mathcal{H}} =\nabla \phi (F)\left\langle
DF,DF\right\rangle_{\mathcal{H}} =\nabla \phi (F)\sigma _{F}.
\end{equation*}%
It follows that $\nabla \phi (F)=\gamma _{F}\left\langle D\phi
(F),DF\right\rangle_{\mathcal{H}}$. Then, by using (\ref{FD12}) and the duality formula (\ref{S7}),
\begin{eqnarray*}
{\mathbb{E}}(G\nabla \phi (F)) &=&{\mathbb{E}}(G\gamma _{F}\left\langle
D\phi (F),DF\right\rangle_{\mathcal{H}} )={\mathbb{E}}(G\gamma _{F}(L(\phi (F)F)-\phi
(F)LF+FL\phi (F)) \\
&=&{\mathbb{E}}(\phi (F)(FL(G\gamma _{F})+G\gamma _{F}LF+L(G\gamma _{F}F)).
\end{eqnarray*}%
We use once again (\ref{FD12}) in order to obtain $H_{i}(F,G)$ in (\ref{FD14}%
). By iteration one obtains the higher order integration by parts formulae.
$\square $

\medskip

We give now useful estimates for the weights which appear in (\ref{FD15}):

\begin{lemma}
\label{L1} Let $F\in \mathcal{S}^{M}$ be such that $\overline{\sigma }%
_{F}(p)<\infty $ for every $p\geq 1$ and let $G\in \mathcal{S}$. Then for
each $m,q\in {\mathbb{N}}$ there exists a universal constant $C\geq 1$
(depending on $M,m,q$ only) such that for every multi-index $\alpha $ with $%
\left\vert \alpha \right\vert \leq q$ one has%
\begin{equation}
\left\vert H_{\alpha }(F,G)\right\vert _{m}\leq C(1\vee (\det \sigma
_{F})^{-1})^{q(m+1)}(1+\left\vert F\right\vert
_{1,m+q+2}^{2Mq(m+2)}+\left\vert LF\right\vert _{m+q}^{q})\left\vert
G\right\vert _{m+q}.  \label{FD16}
\end{equation}%
In particular we have%
\begin{equation}
\left\Vert H_{\alpha }(F,G)\right\Vert _{p}\leq C\overline{\sigma }%
_{F}(2pq)(1+\left\Vert \left\vert F\right\vert \right\Vert
_{1,q+2,4p}^{6qM})\left\Vert G\right\Vert _{q,4p}  \label{FD4}
\end{equation}
\end{lemma}

The proof is long so we skip it, details may be found in \cite{[BCl]} and in
\cite{[BC-EJP]} Theorem 3.4. We will also need the following:

\begin{lemma}
For every $l\in\N$ there exists a constant $%
C_{l}\geq 1$ such that for every $q\in {\mathbb{N}}$, $p\geq 1$ and  $G\in {\mathbb{D}}^{q,\infty }$,
\begin{equation}
\|G^{l}\|_{q,p}\leq C_{l}\left\Vert G\right\Vert _{q,2^lp}^{l}  \label{FD4'}
\end{equation}
\end{lemma}

The proof is straightforward so we skip it.

\subsection{Regularization and non degeneracy}

In this section we consider a functional $F\in ({\mathbb{D}}^{2,\infty
})^{M} $. As it is clear from (\ref{FD16}), a delicate point in using the
integration by parts formulae is to ensure that the functionals at hand are
non degenerate, that is $\det \sigma _{F}>0.$ And in fact this is never
true almost surely: this is because $\chi _{1}=\cdots=\chi _{m}=0$ with
strictly positive probability. In order to bypass this difficulty we use a
regularization argument involving the lowest eigenvalue $\lambda _{F}$ of $%
\sigma _{F}$. We denote%
\begin{equation}
\lambda _{\delta ,\eta ,q}(F)=\delta ^{-1}{\mathbb{P}}(\lambda _{F}\leq \eta
)^{1/q}+\eta ^{-1},\quad \delta ,\eta >0,q\in {\mathbb{N}}.  \label{FD17}
\end{equation}%
We also set
\begin{equation}
\gamma _{\delta }(z)=\frac{1}{m(1)\sqrt{\delta }}\psi _{1}(\delta
^{-1}\left\vert z\right\vert ^{2})  \label{FD17'}
\end{equation}%
where $\psi _{1}$ is the function defined in (\ref{N4}) and $m(1)$ is the
normalization constant from (\ref{N6}) (with $r=1)$. For $f:{\mathbb{R}}%
^{d}\rightarrow {\mathbb{R}}$ we denote%
\begin{equation}
f_{\delta }=f\ast \gamma _{\delta },  \label{FD17''}
\end{equation}%
the symbol $\ast $ denoting convolution. We also consider a supplementary
random variable
\begin{equation*}
Z_{0}\sim \frac{1}{m(1)}\psi _{1}(\left\vert z\right\vert ^{2})dz
\end{equation*}%
which we assume to be independent of $Z_{k},k\in {\mathbb{N}}$, and we
define
\begin{equation}
F_{\delta }=F+\sqrt{\delta }Z_{0}.  \label{FD18}
\end{equation}

\begin{lemma}
\label{L2} Let $q\in {\mathbb{N}}$, $M\in {\mathbb{N}}$ with $M\geq 1$ and $%
p_{1},p_{2},p_{3}>0$ be such that $p_{1}^{-1}+p_{2}^{-1}+p_{3}^{-1}=1$. Then
there exists a constant $C$ such that for every $\delta ,\eta >0,$ every
multi-index $\alpha $ with $\left\vert \alpha \right\vert =q$, every
measurable function $f:{\mathbb{R}}^{M}\rightarrow {\mathbb{R}}$ and every $%
F\in ({\mathbb{D}}^{q+2,\infty })^{M}$ and $G\in {\mathbb{D}}^{q+1,\infty }$
one has%
\begin{equation}
\left\vert {\mathbb{E}}(\partial _{\alpha }f_{\delta }(F)G)\right\vert \leq
C\lambda _{\delta ,\eta ,Mqp_{1}}^{Mqp_{1}}(F)\left\Vert f\right\Vert
_{\infty }(1+\left\Vert \left\vert F\right\vert \right\Vert
_{1,q+2,4Mqp_{2}}^{4Mq})\left\Vert G\right\Vert _{q,qp_{3}}.  \label{FD26}
\end{equation}
\end{lemma}

\textbf{Proof}. We notice first that
\begin{equation*}
{\mathbb{E}}(\partial _{\alpha }f_{\delta }(F)G)={\mathbb{E}}(\partial
_{\alpha }f(F_{\delta })G).
\end{equation*}%
We work with the integration by parts formula based on $Z_{0}$ and on $%
Z_{k},k\in {\mathbb{N}}$. Then
\begin{equation*}
\mbox{$D_{k}F_{\delta } =D_{k}F$ for $k\geq 1$ and $D_0F_\delta=\sqrt{\delta
}$ for $k=0$.}
\end{equation*}%
and $D^{(n)}F_{\delta }=D^{(n)}F$ for $n\geq 2$. So $\left\Vert F_{\delta
}\right\Vert _{1,q+2,p}\leq \left\Vert F\right\Vert _{1,q+2,p}+\sqrt{\delta }%
.$ Moreover
\begin{equation*}
LF_{\delta }=LF+\sqrt{\delta }(\ln \psi _{1})^{\prime }(\left\vert
Z_{0}\right\vert ^{2})\times 2Z_{0}
\end{equation*}%
so that, by (\ref{N5}) we obtain $\left\Vert LF_{\delta }\right\Vert
_{q+2,p}\leq \left\Vert LF\right\Vert _{q+2,p}+C\sqrt{\delta }.$ We conclude
that
\begin{equation*}
\left\Vert \left\vert F_{\delta }\right\vert \right\Vert _{1,q+2,p}\leq
\left\Vert \left\vert F\right\vert \right\Vert _{1,q+2,p}+C\sqrt{\delta }.
\end{equation*}%
We look now to the covariance matrix:
\begin{eqnarray*}
\left\langle \sigma _{F_{\delta }},\xi ,\xi \right\rangle
&=&\sum_{k=0}^{\infty }\left\langle D_{k}F_{\delta },\xi \right\rangle
^{2}=\langle D_{0}(\sqrt{\delta }Z_{0}),\xi \rangle
^{2}+\sum_{k=0}^{\infty }\left\langle D_{k}F,\xi \right\rangle ^{2} \\
&\geq &\delta \left\vert \xi \right\vert ^{2}+\lambda _{F}\left\vert \xi
\right\vert ^{2},
\end{eqnarray*}%
so the lowest eigenvalue of $\sigma _{F_{\delta }}$ verifies $\lambda
_{F_{\delta }}\geq \delta +\lambda _{F}.$ Using the integration by parts
formulae (\ref{FD15}) for $F_{\delta }$ we obtain
\begin{equation*}
{\mathbb{E}}(\partial _{\alpha }f(F_{\delta })G)={\mathbb{E}}(f(F_{\delta
})H_{\alpha }(F_{\delta },G)).
\end{equation*}%
By (\ref{FD16}) with $m=0$%
\begin{equation*}
\left\vert H_{\alpha }(F_{\delta },G)\right\vert \leq C(1\vee (\lambda
_{F_{\delta }})^{-1})^{Mq}(1+\left\vert F\right\vert
_{1,q+2}^{4Mq}+\left\vert LF\right\vert _{q}^{q})\left\vert G\right\vert _{q}
\end{equation*}%
so that, using H\"{o}lder's inequality
\begin{equation*}
{\mathbb{E}}(\left\vert H_{\alpha }(F_{\delta },G)\right\vert )\leq C(1\vee {%
\mathbb{E}}((\lambda _{F_{\delta }})^{-Mqp_{1}})^{1/p_{1}}(1+\left\Vert
F\right\Vert _{1,q+2,Mqp_{2}}^{4Mq}+\left\Vert LF\right\Vert
_{q,qp_{2}}^{4Mq})\left\Vert G\right\Vert _{q,qp_{3}}.
\end{equation*}%
We write now%
\begin{eqnarray*}
{\mathbb{E}}((\lambda _{F_{\delta }})^{-Mqp_{1}}) &=&{\mathbb{E}}((\lambda
_{F_{\delta }})^{-Mqp_{1}}1_{\{\lambda _{F}\leq \eta \}})+{\mathbb{E}}%
((\lambda _{F_{\delta }})^{-Mqp_{1}}1_{\{\lambda _{F}>\eta \}}) \\
&\leq &\delta ^{-Mqp_{1}}{\mathbb{P}}(\lambda _{F}\leq \eta )+\eta
^{-Mqp_{1}} \\
&\leq &(\delta ^{-1}{\mathbb{P}}(\lambda _{F}\leq \eta )^{1/Mqp_{1}}+\eta
^{-1})^{Mqp_{1}}=\lambda _{\delta ,\eta ,Mqp_{1}}^{Mqp_{1}}(F).
\end{eqnarray*}%
We conclude that
\begin{eqnarray*}
{\mathbb{E}}(\partial _{\alpha }f_{\delta }(F)G) &=&\left\vert {\mathbb{E}}%
(\partial _{\alpha }f(F_{\delta })G)\right\vert \leq \left\Vert f\right\Vert
_{\infty }{\mathbb{E}}(\left\vert H_{\alpha }(F_{\delta },G)\right\vert ) \\
&\leq &C\lambda _{\delta ,\eta ,Mqp_{1}}^{Mq}(F)\left\Vert f\right\Vert
_{\infty }(1+\left\Vert F\right\Vert _{1,q+2,Mqp_{2}}^{4Mq}+\left\Vert
LF\right\Vert _{q,4qp_{2}}^{4Mq})\left\Vert G\right\Vert _{q,qp_{3}}
\end{eqnarray*}%
$\square $

\medskip

In order to pass from $f_{\delta }$ to $f$ we will use the following lemma:

\begin{lemma}
\label{L3}Let $M\in{\mathbb{N}}$, $M\geq 1$. There exist
constants $C,p,a\geq 1$ such that for every $\eta >0,\delta >0,$ every $F\in
({\mathbb{D}}^{3,p})^{M}$ and every bounded and measurable $f:{\mathbb{R}}%
^{M}\rightarrow {\mathbb{R}}$ one has%
\begin{equation}
\left\vert {\mathbb{E}}(f(F))-{\mathbb{E}}(f_{\delta }(F))\right\vert \leq
C\left\Vert f\right\Vert _{\infty }\Big({\mathbb{P}}(\lambda _{F}<\eta )+%
\frac{\sqrt{\delta }}{\eta ^{p}}(1+\left\Vert \left\vert F\right\vert
\right\Vert _{3,p})^{a}\Big)  \label{FD24'}
\end{equation}
\end{lemma}

The above Lemma is Lemma 2.5 in \cite{[BC-EJP]}. There $\gamma _{\delta }$
is the Gaussian density of covariance $\delta $ but the proof is exactly the
same with $\gamma _{\delta }$ defined in (\ref{FD17'}).

\section{Sobolev norms and non-degeneracy for stochastic series}
\label{sect-MallS}

The stochastic series $S(c,Z)=\sum_{m=1}^{\infty }\Phi _{m}(c,Z)$ are a
natural generalization of the decomposition in Wiener chaoses - indeed, if $%
Z_{k}$ are standard normal then the $\Phi _{m}(c,Z)$'s represent multiple
stochastic integrals of order $m$ (Remark \ref{iterated}). The aim of this section is to obtain
estimates of the Sobolev norms of $S(c,Z)$ and of $LS(c,Z)$ which are
analogous to the ones we have in the Gaussian case. To this purpose, it is
useful to introduce random variables taking values on a Hilbert space $%
\mathcal{U}$ that are derivable in Malliavin sense. In fact, $DS(c,Z)$ is a
r.v. in $\mathcal{H}$ (see (\ref{H})) and can be written again as a
stochastic series whose coefficients are in $\mathcal{H}$. So, in order to
handle properly our problem, we consider stochastic series whose
coefficients $c(\alpha )$ belongs to a separable Hilbert space ${\mathcal{U}}
$. We denote with $\<\cdot ,\cdot \>_{\mathcal{U}}$ and $|\cdot |_{\mathcal{U%
}}$ the associated inner product and norm, respectively. We set $L_{{%
\mathcal{U}}}^{p}=\{F\in {\mathcal{U}}:\left\Vert F\right\Vert _{\U,p}:=({%
\mathbb{E}}(\left\vert F\right\vert _{{\mathcal{U}}}^{p}))^{1/p}<\infty \}.$
Even if $F\in {\mathbb{R}}$ (as it is the case in our paper), the derivative
$DF$ takes values in ${\mathcal{U}}=\mathcal{H}$ which is a Hilbert space.

We set
\begin{equation}
\mathcal{H}(\mathcal{U})=\{x\in \mathcal{U}^{\mathbb{N}}\,:\,\sum_{k=1}^{%
\infty }|x_{k}|_{\mathcal{U}}^{2}<\infty \}.  \label{H-U}
\end{equation}%
$\mathcal{H}(\mathcal{U})$ is clearly a Hilbert space, the inner product and
the norm being given by
\begin{equation*}
\<x,y\>_{\mathcal{H}(\mathcal{U})}=\sum_{k=1}^{\infty }\<x_{k},y_{k}\>_{%
\mathcal{U}}\quad \mbox{and}\quad |x|_{\mathcal{H}(\mathcal{U})}=\Big(%
\sum_{k=1}^{\infty }|x_{k}|_{\mathcal{U}}^{2}\Big)^{1/2}
\end{equation*}%
respectively. Notice that $\mathcal{H}({\mathbb{R}})$ is the space $\mathcal{%
H}$ defined in (\ref{H}). Remark also that $\mathcal{H}(\mathcal{H}(\mathcal{%
U}))=\mathcal{H}(\mathcal{U})^{\otimes 2}$, and more generally $\mathcal{H}(%
\mathcal{H}(\mathcal{U})^{\otimes n})=\mathcal{H}(\mathcal{U})^{\otimes
(n+1)}$.

Let $A$ be a random variable taking values in $\mathcal{U}$ which is
measurable with respect to $\sigma(Z_i, i=1,2,\ldots)$ and take $p\geq 1$.
We set ${\mathbb{D}}^{0,p}_{\mathcal{U}}=L^p_{\mathcal{U}}$ and we set
\begin{equation*}
\|A\|_{\mathcal{U},0,p}\equiv \|A\|_{\mathcal{U},p}=\||A|_{\mathcal{U}}\|_p.
\end{equation*}
We say that $A\in {\mathbb{D}}^{1,p}_{\mathcal{U}}$ if $A\in L^p_{\mathcal{U}%
}$, $\<A,h\>_{\mathcal{U}}\in{\mathbb{D}}^{1,p}$ for every $h\in \mathcal{U}$
and there exists $DA\in\mathcal{H}(\mathcal{U})$ such that
\begin{equation*}
DA\in L^p_{\mathcal{H}(\mathcal{U})} \quad\mbox{and}\quad \<(DA)_k,h\>_{%
\mathcal{U}}=D_k\<A,h\>_{\mathcal{U}},\ k\in{\mathbb{N}}.
\end{equation*}
In the following, we use the notation $D_kA=(DA)_k$, $k\in{\mathbb{N}}$. For
$A\in{\mathbb{D}}^{1,p}_{\mathcal{U}}$ we define the Sobolev norm
\begin{equation*}
\|A\|_{\mathcal{U},1,p}=\|A\|_{{\mathcal{U}},p} +\|DA\||_{\mathcal{H}(%
\mathcal{U}),p}.
\end{equation*}
Note that if $\mathcal{U}={\mathbb{R}}$ then the above definition $F\in{%
\mathbb{D}}^{1,p}$ agrees with the standard definition: ${\mathbb{D}}^{1,p}_{%
\mathbb{R}}\equiv {\mathbb{D}}^{1,p}$.

This reasoning can be iterated in order to define a random variable $A$ in $%
\mathcal{U}$ which is $q\geq 2$ times differentiable: for $p\geq 1$, we say
that $A\in {\mathbb{D}}_{\mathcal{U}}^{q,p}$ if $A\in {\mathbb{D}}_{\mathcal{%
U}}^{q-1,p}$, $\<A,h\>_{\mathcal{U}}\in {\mathbb{D}}^{q,p}$ for every $h\in
\mathcal{U}$ and there exists $D^{(q)}A\in \mathcal{H}(\mathcal{U})^{\otimes
q}$ such that
\begin{equation*}
D^{(q)}A\in L_{\mathcal{H}(\mathcal{U})^{\otimes q}}^{p}\quad \mbox{and}%
\quad <(D^{(q)}A)_{\alpha },h\>_{\mathcal{U}}=D_{\alpha }^{(q)}<A,h\>_{%
\mathcal{U}},\ |\alpha |=q.
\end{equation*}%
We set $D_{\alpha }^{(q)}A=(D^{(q)}A)_{\alpha }$, $|\alpha |=q$, and
\begin{equation*}
\Vert A\Vert _{\mathcal{U},q,p}=\Vert A\Vert _{{\mathcal{U}},q-1,p}+\Vert
D^{(q)}A\Vert _{\mathcal{H}(\mathcal{U})^{\otimes q},p}=\Vert A\Vert _{{%
\mathcal{U}},p}+\sum_{j=1}^{q}\Vert D^{(j)}A\Vert _{\mathcal{H}(\mathcal{U}%
)^{\otimes j},p}.
\end{equation*}%
As an example, take $F$ a random variable in ${\mathbb{R}}$. Then $F\in {%
\mathbb{D}}^{2,p}$ if $F\in {\mathbb{D}}^{1,p}$ (standard definition) and $%
DF\in {\mathbb{D}}_{\mathcal{H}}^{1,p}$ following the above definition, that
is looking at $A=DF$ as a random variable taking values in the Hilbert space
$\mathcal{U}=\mathcal{H}$, and one has $\Vert F\Vert _{2,p}=\Vert F\Vert _{{%
\mathcal{U}},2,p}$ with ${\mathcal{U}}={\mathbb{R}}$. And in general, $F\in {%
\mathbb{D}}^{q,p}$ if for every $1\leq k\leq q-1$, $D^{(k)}F\in {\mathbb{D}}%
_{\mathcal{H}^{\otimes k}}^{1,p}$ and one has $\Vert F\Vert _{q,p}=\Vert
F\Vert _{{\mathcal{U}},q,p}$ with ${\mathcal{U}}={\mathbb{R}}$.

We consider now the random variable $S(c,Z)=\sum_{\alpha }c(\alpha
)Z^{\alpha }$ in (\ref{Law7}) with $c(\alpha )\in {\mathcal{U}}$ for every $%
\alpha $. We denote%
\begin{equation}
\left\vert c\right\vert _{{\mathcal{U}},m}=\Big(\sum_{\alpha \in \Gamma
_{m}}\left\vert c(\alpha )\right\vert _{{\mathcal{U}}}^{2}\Big)^{1/2}.
\label{nss2}
\end{equation}%
For $q,p\in {\mathbb{N}}$ and $M\in {\mathbb{R}}$ we set
\begin{equation}
N_{q}(c,M)=\Big(\sum_{m=q}^{\infty }M^{m-q}\times \frac{m!}{(m-q)!}\times
m!\left\vert c\right\vert _{{\mathcal{U}},m}^{2}\Big)^{1/2}  \label{nss3}
\end{equation}
(for a comparison, see (\ref{Law3}) for $\U=\R$).

We recall that we are assuming that $Z_{k},k\in {\mathbb{N}}$ are independent and%
\begin{equation}
{\mathbb{E}}(Z_{k})=0\qquad \mbox{and}\qquad {\mathbb{E}}(Z_{k}^{2})=1.
\label{nss3'}
\end{equation}

The following basic relations are immediate consequences of the isometry
property for Hilbert space valued discrete martingales. Let $\Phi
_{m}(c,Z)=\sum_{\left\vert \alpha \right\vert =m}c(\alpha )Z^{\alpha }$ with
$c(\alpha )\in {\mathcal{U}}.$ Then%
\begin{equation}
{\mathbb{E}}(\left\langle \Phi _{m}(c,Z),\Phi _{m^{\prime
}}(c,Z)\right\rangle _{{\mathcal{U}}})=\left\{
\begin{array}{ll}
0 & \qquad m\neq m^{\prime } \\
m!\left\vert c\right\vert _{{\mathcal{U}},m}^{2} & \qquad m=m^{\prime }.%
\end{array}%
\right.  \label{nss5}
\end{equation}

\begin{lemma}
\label{ES} Let $Z_{k},k\in {\mathbb{N}}$, be independent such that (\ref%
{nss3'}) holds. Let $S(c,Z)=\sum_{m=1}^{\infty }\Phi _{m}(Z)$ and $%
N_{q}(c,M) $ be defined in (\ref{nss3}).

\medskip

$(i)$ The series $S(c,Z)$ is convergent in $L_{{\mathcal{U}}}^{2}$ if and
only if $N_{0}(c,1)<\infty $ and in this case%
\begin{equation*}
\left\Vert S(c,Z)\right\Vert _{{\mathcal{U}},2}=N_{0}(c,1)=(\sum_{m=1}^{%
\infty }m!\left\vert c\right\vert _{{\mathcal{U}},m}^{2})^{1/2}.
\end{equation*}

$(ii)$ Let $b_{p}$ be the constant in the Burkholder's inequality (see (\ref%
{a1})) and%
\begin{equation}
M_{p}=\sqrt{2}b_{p}M_{p}(Z),\quad \mbox{with}\qquad M_{p}(Z)=\sup_{k\in {%
\mathbb{N}}}\left\Vert Z_{k}\right\Vert _{p}.  \label{nss4'}
\end{equation}%
Then%
\begin{equation}
\left\Vert S(c,Z)\right\Vert _{{\mathcal{U}},p}\leq N_{0}(c,M_{p}^{2})=\Big(%
\sum_{m=1}^{\infty }M_{p}^{2m}m!\left\vert c\right\vert _{{\mathcal{U}}%
,m}^{2}\Big)^{1/2}.  \label{nss4}
\end{equation}
\end{lemma}

\begin{remark}
\textrm{As an immediate consequence one has the following estimate for
multiple integrals in Wiener chaoses: $\left\Vert I_{m}(f)\right\Vert
_{p}\leq \sqrt{2}b_{p}M_{p}(G)\sqrt{m!}\left\Vert f\right\Vert
_{L^{2}(R_{+}^{m})}$ where $G\sim \mathcal{N}(0,1).$ More accurate estimates
concerning the Gaussian chaoses can be found in Latala \cite{L}. }
\end{remark}

\textbf{Proof of Lemma \ref{ES}.} $(i)$ immediately follows from (\ref{nss5}). As for $(ii)$,
we fix $J\in {\mathbb{N}}$ and we denote $\left\vert c\right\vert _{{%
\mathcal{U}},m,J}^{2}=\sum_{\alpha \in \Gamma _{m}(J)}\left\vert c(\alpha
)\right\vert _{{\mathcal{U}}}^{2}$ and%
\begin{equation*}
\Phi _{m,J}(c,Z)=\sum_{\alpha \in \Gamma _{m}(J)}c_{\alpha }Z^{\alpha
},\qquad \Phi _{m,J}^{o}(c,Z)=\sum_{\alpha \in \Gamma _{m}^{o}(J)}c(\alpha
)Z^{\alpha }.
\end{equation*}%
We set
\begin{equation*}
S_{n,J}(c,Z)=\sum_{m=1}^{n}\Phi _{m,J}(c,Z)=\sum_{m=1}^{n}m!\Phi
_{m,J}^{o}(c,Z)
\end{equation*}%
and we prove that for every $n,J\in {\mathbb{N}}$
\begin{equation}
\left\Vert S_{n,J}(c,Z)\right\Vert _{{\mathcal{U}},p}^{2}\leq
\sum_{m=1}^{n}M_{p}^{2m}m!\left\vert c\right\vert _{{\mathcal{U}}%
,m,J}^{2}\leq N_{0}^{2}(c,M_{p}).  \label{nss7}
\end{equation}%
Then (\ref{nss4}) follows by passing to the limit. We prove (\ref{nss7}) by
recurrence on $n$. For $n=1$, we use the Burkholder inequality (\ref{a2})
and we have
\begin{equation*}
\Vert S_{1,J}(c,Z)\Vert _{{\mathcal{U}},p}^{2}=\Big\|\sum_{j=1}^{J}c_{j}Z_{j}%
\Big\|_{{\mathcal{U}},p}^{2}\leq b_{p}^{2}\sum_{j=1}^{J}\Vert
c_{j}Z_{j}\Vert _{{\mathcal{U}},p}^{2}\leq
b_{p}^{2}M_{p}(Z)^{2}\sum_{j=1}^{J}|c_{j}|_{\mathcal{U}}^{2}\leq
M_{p}^{2}|c|_{{\mathcal{U}},1,J}^{2}.
\end{equation*}%
For $n>1$, we use the following basic decomposition:
\begin{equation*}
\Phi _{m,J}^{o}(c,Z)=\sum_{j=1}^{J}Z_{j}\sum_{\alpha \in \Gamma
_{m-1}^{o}(j-1)}c(\alpha ,j)Z^{\alpha }.
\end{equation*}%
This gives%
\begin{eqnarray*}
S_{n,J}(c,Z) &=&\sum_{m=1}^{n}m!\Phi
_{m,J}^{o}(c,Z)=\sum_{j=1}^{J}Z_{j}\sum_{m=1}^{n}m!\sum_{\alpha \in \Gamma
_{m-1}^{o}(j-1)}c(\alpha ,j)Z^{\alpha } \\
&=&\sum_{j=1}^{J}Z_{j}c_{j}+\sum_{j=1}^{J}Z_{j}\sum_{m=2}^{n}\sum_{\alpha
\in \Gamma _{m-1}(J)}(m1_{\Gamma (j-1)}(\alpha )c(\alpha ,j))Z^{\alpha } \\
&=&\sum_{j=1}^{J}Z_{j}(c_{j}+S_{n-1,J}(c^{j}))
\end{eqnarray*}%
with%
\begin{equation*}
c^{j}(\alpha )=(1+\left\vert \alpha \right\vert )\times 1_{\Gamma
(j-1)}(\alpha )c(\alpha ,j)\qquad \alpha \in \Gamma _{m-1}(J).
\end{equation*}%
Notice that $S_{n-1,J}(c^{j})$ is measurable with respect to $\sigma
(Z_{1},...,Z_{j-1})$ so that $S_{n,J}(c,Z)$ is a martingale. Then, by
Burkholder's inequality (\ref{a2}),%
\begin{eqnarray*}
\left\Vert S_{n,J}(c,Z)\right\Vert _{{\mathcal{U}},p}^{2} &\leq
&b_{p}^{2}\sum_{j=1}^{J}\left\Vert Z_{j}\right\Vert _{p}^{2}\left\Vert
c_{j}+S_{n-1,J}(c^{j},Z)\right\Vert _{{\mathcal{U}},p}^{2} \\
&\leq &2b_{p}^{2}M_{p}^{2}(Z)\Big(\sum_{j=1}^{J}\left\vert c_{j}\right\vert
_{{\mathcal{U}}}^{2}+\sum_{j=1}^{J}\left\Vert S_{n-1,J}(c^{j},Z)\right\Vert
_{{\mathcal{U}},p}^{2}\Big).
\end{eqnarray*}%
Using the recurrence hypothesis%
\begin{eqnarray*}
\sum_{j=1}^{J}\left\Vert S_{n-1,J}(c^{j},Z)\right\Vert _{{\mathcal{U}}%
,p}^{2} &\leq &\sum_{j=1}^{J}\sum_{m=1}^{n-1}M_{p}^{2m}m!\left\vert
c^{j}\right\vert _{{\mathcal{U}},m,J}^{2} \\
&=&\sum_{m=1}^{n-1}M_{p}^{2m}m!(m+1)^{2}\sum_{j=1}^{J}\sum_{\alpha \in
\Gamma _{m}(j-1)}\left\vert c(\alpha ,j)\right\vert _{{\mathcal{U}}}^{2}
\end{eqnarray*}%
Since%
\begin{equation*}
\left\vert c\right\vert _{{\mathcal{U}},m+1,J}^{2}=(m+1)\sum_{j=1}^{J}\sum_{%
\alpha \in \Gamma _{m}(j-1)}\left\vert c(\alpha ,j)\right\vert _{{\mathcal{U}%
}}^{2}
\end{equation*}%
we obtain%
\begin{equation*}
\sum_{j=1}^{J}\left\Vert S_{n-1,J}(c^{j},Z)\right\Vert _{{\mathcal{U}}%
,p}^{2}\leq \sum_{m=1}^{n-1}M_{p}^{2m}(m+1)!\left\vert c\right\vert _{{%
\mathcal{U}},m+1,J}^{2}.
\end{equation*}%
We conclude that%
\begin{equation*}
\left\Vert S_{n,J}(c,Z)\right\Vert _{{\mathcal{U}},p}^{2}\leq M_{p}^{2}\Big(%
\left\vert c\right\vert _{{\mathcal{U}},1}^{2}+%
\sum_{m=2}^{n}M_{p}^{2(m-1)}m!\left\vert c\right\vert _{{\mathcal{U}}%
,m,J}^{2}\Big)=\sum_{m=1}^{n}M_{p}^{2m}m!\left\vert c\right\vert _{{\mathcal{%
U}},m,J}^{2}.
\end{equation*}

$\square $

\medskip

We estimate now the derivatives of $S(c,Z).$

\begin{proposition}
\label{S-S} We assume that (\ref{nss3'}) holds and, for $p\geq 2$, set $M_p$
as in (\ref{nss4'}). %
%
%
%
For $q\in{\mathbb{N}}$ one has 
\begin{equation*}
\| D^{(q)}S(c,Z)\|_{\mathcal{H}({\mathcal{U}})^{\otimes q},p}\leq \sqrt{2}%
\Big(\sum_{n=q}^{\infty }n!\frac{n!}{(n-q)!}M_{p}^{2(n-q)} \left\vert
c\right\vert _{{\mathcal{U}},n}^{2}\Big)^{1/2}=\sqrt{2}N_{q}(c,M_{p}^{2}).
\end{equation*}
As a consequence,
\begin{equation}
\left\Vert S(c,Z)\right\Vert _{\mathcal{U}, q,p}\leq \sqrt{2}%
\sum_{k=0}^qN_{k}(c,M_{p}^{2}).  \label{nss8}
\end{equation}
\end{proposition}

\textbf{Proof}.
%
Let us denote $\partial _{j}\Phi _{m}(c,Z)=\partial _{Z_{j}}\Phi _{m}(c,Z)$
so that $D_{j}\Phi _{m}(c,Z)=\chi _{j}\partial _{j}\Phi _{m}(c,Z).$ We write%
\begin{equation*}
\partial _{j}\Phi _{m}(c,Z)=\sum_{\beta \in \Gamma _{m}}c(\beta )\partial
_{j}Z^{\beta }=m\sum_{\substack{ \alpha \in \Gamma _{m-1}  \\ j\notin \alpha
}}c(\alpha ,j)Z^{\alpha }=\Phi _{m-1}(\overline{c}^{j},Z)
\end{equation*}%
with
\begin{equation*}
\overline{c}^{j}(\alpha )=(1+\left\vert \alpha \right\vert )c(\alpha
,j)1_{\{j\notin \alpha \}}.
\end{equation*}%
So we may write%
\begin{equation*}
\partial _{j}S(c,Z)=\sum_{m=1}^{\infty }\partial _{j}\Phi
_{m,J}(c,Z)=c_{j}+\sum_{m=2}^{\infty }\Phi _{m-1,J}(\overline{c}^{j},Z).
\end{equation*}%
We set $\overline{c}_{0}=(c_{j})_{j\in {\mathbb{N}}}$ and $\overline{c}%
_{\alpha }=(\overline{c}_{\alpha }^{j})_{j\in {\mathbb{N}}}$. Notice that $%
\overline{c}_{0},\overline{c}_{\alpha }\in \mathcal{H}({\mathcal{U}})$, so
the above equality reads%
\begin{equation}
\nabla S(c,Z)=\overline{c}_{0}+\sum_{m=1}^{\infty }\Phi _{m,J}(\overline{c}%
,Z)=\overline{c}_{0}+S(\overline{c},Z)\in \mathcal{H}({\mathcal{U}})
\label{nss11}
\end{equation}%
where $\nabla S(c,Z)=(\partial _{j}S(c,Z))_{j\in {\mathbb{N}}}\in \mathcal{H}%
({\mathcal{U}})$. By using (\ref{nss4}),
\begin{equation}
\Vert \nabla S(c,Z)\Vert _{\mathcal{H}({\mathcal{U}}),p}^{2}\leq 2(|%
\overline{c}_{0}|_{\mathcal{H}({\mathcal{U}})}^{2}+\Vert S(\overline{c}%
,Z)\Vert _{\mathcal{H}({\mathcal{U}}),p}^{2})\leq 2\Big(|\overline{c}(0)|_{%
\mathcal{H}({\mathcal{U}})}^{2}+\sum_{m=1}^{\infty }M_{p}^{2m}m!|\overline{c}%
|_{\mathcal{H}({\mathcal{U}}),m}^{2}\Big)  \label{nss11bis}
\end{equation}%
We have $|\overline{c}_{0}|_{\mathcal{H}({\mathcal{U}})}^{2}=\sum_{j=1}^{%
\infty }|c_{j}|_{\mathcal{U}}^{2}=|c|_{{\mathcal{U}},1}^{2}$ and
\begin{eqnarray}
\left\vert \overline{c}\right\vert _{\mathcal{H}({\mathcal{U}}),m}^{2}
&=&\sum_{j=1}^{\infty }\sum_{\alpha \in \Gamma _{m}}\left\vert \overline{c}%
^{j}(\alpha )\right\vert _{{\mathcal{U}}}^{2}=\sum_{j=1}^{\infty
}\sum_{\alpha \in \Gamma _{m}}(1+\left\vert \alpha \right\vert
)^{2}\left\vert c(\alpha ,j)\right\vert _{{\mathcal{U}}}^{2}1_{\{j\notin
\alpha \}}  \label{nss12} \\
&=&(m+1)^{2}\sum_{\beta \in \Gamma _{m+1}}\left\vert c(\beta )\right\vert _{{%
\mathcal{U}}}^{2}=(m+1)^{2}\left\vert c\right\vert _{{\mathcal{U}},m+1}^{2}.
\notag
\end{eqnarray}%
By inserting in (\ref{nss11bis}), this gives%
\begin{align*}
\Vert \nabla S(c,Z)\Vert _{\mathcal{H}({\mathcal{U}}),p}^{2}& \leq 2\Big(%
|c|_{{\mathcal{U}},1}^{2}+\sum_{m=1}^{\infty
}M_{p}^{2m}m!(m+1)^{2}\left\vert c\right\vert _{{\mathcal{U}},m+1}^{2}\Big)
\\
& =2\Big(|c|_{{\mathcal{U}},1}^{2}+\sum_{m=1}^{\infty
}M_{p}^{2m}(m+1)!(m+1)\left\vert c\right\vert _{{\mathcal{U}},m+1}^{2}\Big)
\\
& \leq 2\sum_{m=1}^{\infty }M_{p}^{2m}m!m\left\vert c\right\vert _{{\mathcal{%
U}},m}^{2}=2N_{1}(c,M_{p}^{2})^{2},
\end{align*}%
that is $\Vert \nabla S(c,Z)\Vert _{\mathcal{H}({\mathcal{U}}),p}\leq \sqrt{2%
}\,N_{1}(c,M_{p}^{2})$. Finally, one has $D_kS(c,Z)=\chi_k\partial_jS(c,Z)$
and since $\chi_k\in[0,1]$ then $| DS(c,Z)| _{\mathcal{H}(\mathcal{U})}$ $%
\leq|\nabla S(c,Z)|_{\mathcal{H}({\mathcal{U}})}$. So, $\| DS(c,Z)\| _{%
\mathcal{H}(\mathcal{U}),p}$ $\leq \|\nabla S(c,Z)\|_{\mathcal{H}({\mathcal{U%
}}),p}$ and the proof is completed for the first order derivative.

Let us now estimate the second order derivatives. We have $DD\Phi
_{m,J}(c,Z)=D\Phi _{m,J}(\overline{c},Z).$ Using (\ref{nss12}) one checks
that
\begin{equation*}
\left\vert \overline{c}\right\vert _{\mathcal{H}({\mathcal{U}}),1}+N_{1}(%
\overline{c},M_{p}^{2})\leq \sqrt 2\,N_{2}(c,M_{p}^{2})
\end{equation*}%
and so we may use the result for the first order derivatives. For higher
order derivatives the argument is the same. $\square $

\medskip

We estimate now the Sobolev norms of $LS(c,Z).$

\begin{proposition}
\label{S-L} We assume that (\ref{nss3'}) holds and, for $p\geq 2$, set $M_p$
as in (\ref{nss4'}). For every $q\in {\mathbb{N}}$ there exists a constant $%
C $ (depending on $q$ and $p$ only) such that%
\begin{equation}
\left\Vert LS(c,Z)\right\Vert _{{\mathcal{U}},q,p}\leq \frac C{r^{q+1}}\Big(%
\sum_{n=1}^{q+1}|c|_{{\mathcal{U}},n}+N_{q+2}(c,M^2_p)\Big).  \label{nss14}
\end{equation}
\end{proposition}

\textbf{Proof}. Notice that $\left\langle DZ_{k},DZ_{j}\right\rangle =0$ for
$k\neq j.$ Using the computational rules one obtains for $\beta =(\beta
_{1},...,\beta _{m})$ with $m\geq 2$
\begin{equation*}
LZ^{\beta }=\sum_{k=1}^{m}LZ_{\beta _{k}}\prod_{j\neq k}Z_{\beta _{j}}
\end{equation*}%
and, using the symmetry of $c(\beta ),$ this gives%
\begin{eqnarray*}
L\Phi _{m}(c,Z) &=&m!L\Phi _{m}^{o}(c,Z)=m!m\sum_{j=1}^{\infty
}LZ_{j}\sum_{\alpha \in \Gamma _{m-1}^{o}(j-1)}c(\alpha ,j)Z^{\alpha } \\
&=&m^{2}\sum_{j=1}^{\infty }LZ_{j}\sum_{\alpha \in \Gamma
_{m-1}(j-1)}c(\alpha ,j)Z^{\alpha }=\sum_{j=1}^{\infty }LZ_{j}\Phi _{m-1}(%
\widehat{c}^{j},Z)
\end{eqnarray*}%
with%
\begin{equation*}
\widehat{c}^{j}(\alpha )=(1+\left\vert \alpha \right\vert )^{2}\times
1_{\max_k\alpha_k\leq j-1 }\,c(\alpha ,j).
\end{equation*}%
For $m=1$ we have $L\Phi _{1}(c,Z)=\sum_{j=1}^{\infty }c_{j}LZ_{j}.$ It
follows that, for $S_n(c,Z)=\sum_{m=1}^n\Phi_m(c,Z)$,
\begin{eqnarray*}
LS_n(c,Z) &=&\sum_{m=1}^{n}L\Phi _{m}(c,Z)=\sum_{j=1}^{\infty
}c_{j}LZ_{j}+\sum_{j=1}^{\infty }LZ_{j}\sum_{m=2}^{n}\Phi _{m-1}(\widehat{c}%
^{j},Z) \\
&=&\sum_{j=1}^{\infty }LZ_{j}(c_{j}+S_{n-1}(\widehat{c}^{j},Z)).
\end{eqnarray*}%
Notice that $S_{n-1}(\widehat{c}^{j},Z))$ is $\sigma (Z_{1},...,Z_{j-1})$
measurable. Since $LS_n(c)$ verifies (\ref{aa3}) with $B_{j}=c_{j}+S_{n-1}(%
\widehat{c}^{j},Z)$ and $\Lambda _{m}=0$, we will use (\ref{b2}) (actually,
we should use $S_{n,J}(c,Z)$ and then pass to the limit following the
standard technique). We have%
\begin{eqnarray*}
\sum_{j=1}^{\infty }\left\Vert B_{j}\right\Vert _{{\mathcal{U}},q,p}^{2}
&\leq &2\Big(\left\vert c\right\vert _{{\mathcal{U}},1}^{2}+\sum_{j=1}^{%
\infty }\left\Vert S_{n-1}(\widehat{c}^{j},Z))\right\Vert _{{\mathcal{U}}%
,q,p}^{2}\Big) \\
&\leq &C\Big(\left\vert c\right\vert _{{\mathcal{U}},1}^{2}+\sum_{j=1}^{%
\infty }\sum_{\ell =1}^{q}\left\Vert D^{(l)}S_{n-1}(\widehat{c}%
^{j},Z))\right\Vert _{\mathcal{H}({\mathcal{U}})^{\otimes \ell },p}^{2}\Big).
\end{eqnarray*}%
By using (\ref{nss8}) we obtain
\begin{equation*}
\sum_{j=1}^{\infty }\left\Vert B_{j}\right\Vert _{U,q,p}^{2}\leq C\Big(%
\left\vert c\right\vert _{{\mathcal{U}},1}^{2}+\sum_{j=1}^{\infty
}\sum_{\ell =1}^{q}\sum_{n=\ell }^{\infty }n!\frac{n!}{(n-\ell )!}%
M_{p}^{2(n-\ell )}\left\vert \widehat{c}^{j}\right\vert _{{\mathcal{U}}%
,n}^{2}\Big)=:C(\left\vert c\right\vert _{{\mathcal{U}},1}^{2}+I).
\end{equation*}%
In order to handle $I$, we compute
\begin{align*}
\sum_{j=1}^{\infty }\left\vert \widehat{c}^{j}\right\vert _{{\mathcal{U}}%
,n}^{2}& =(n+1)^{4}\sum_{j=1}^{\infty }\sum_{\alpha \in \Gamma
_{n}(j-1)}\left\vert c(\alpha ,j)\right\vert _{{\mathcal{U}}%
}^{2}=(n+1)^{4}n!\sum_{j=1}^{\infty }\sum_{\alpha \in \Gamma
_{n}^{o}(j-1)}\left\vert c(\alpha ,j)\right\vert _{{\mathcal{U}}}^{2} \\
& =(n+1)^{4}n!\sum_{\beta \in \Gamma _{n+1}^{o}}\left\vert c(\beta
)\right\vert _{{\mathcal{U}}}^{2}=(n+1)^{3}\sum_{\beta \in \Gamma
_{n+1}}\left\vert c(\beta )\right\vert _{{\mathcal{U}}}^{2}=(n+1)^{3}\left%
\vert c\right\vert _{{\mathcal{U}},n+1}^{2}.
\end{align*}%
Therefore,
\begin{align*}
I& =\sum_{\ell =1}^{q}\sum_{n=\ell }^{\infty }M_{p}^{2(n-\ell )}(n+1)!\frac{%
(n+1)!}{(n-\ell )!}(n+1)\left\vert c\right\vert _{{\mathcal{U}},n+1}^{2} \\
& =\sum_{\ell =1}^{q}\Big(((\ell +1)!)^{2}(\ell +1)\left\vert c\right\vert _{%
{\mathcal{U}},\ell +1}^{2}+\sum_{m=\ell +2}^{\infty }M_{p}^{2(m-1-\ell )}m!%
\frac{m!}{(m-1-\ell )!}m\left\vert c\right\vert _{{\mathcal{U}},m}^{2}\Big)
\\
& \leq C\sum_{\ell =1}^{q}\big(\left\vert c\right\vert _{{\mathcal{U}},\ell
+1}^{2}+N_{\ell +1}^{2}(c,M_{p}^{2})\big)
\end{align*}%
Now, for $k\leq q+2$ straightforward computations give
\begin{equation*}
N_{k}(c,M)^{2}\leq k!|c|_{{\mathcal{U}},k}^{2}+MN_{k+1}(c,M)^{2}\leq \cdots
\leq C\Big(\sum_{i=k}^{q+1}|c|_{{\mathcal{U}},i}^{2}+N_{q+2}(c,M)^{2}\Big)
\end{equation*}%
and by inserting we obtain
\begin{equation*}
I\leq C\Big(\sum_{n=2}^{q+1}\left\vert c\right\vert _{{\mathcal{U}}%
,n}^{2}+N_{q+2}(c,M_{p}^{2})\Big).
\end{equation*}%
By resuming,
\begin{equation*}
\sum_{j=1}^{\infty }\left\Vert B_{j}\right\Vert _{U,q,p}^{2}\leq C\Big(%
\sum_{n=1}^{q+1}|c|_{{\mathcal{U}},n}^{2}+N_{q+2}(c,M_{p}^{2})^{2}\Big)
\end{equation*}%
and by applying (\ref{b2}) we get
\begin{equation*}
\Vert LS(c,Z)\Vert _{{\mathcal{U}},q,p}\leq \frac{C}{r^{q+1}}\Big(%
\sum_{j=1}^{\infty }\left\Vert B_{j}\right\Vert _{U,q,p}^{2}\Big)^{1/2}\leq
\frac{C}{r^{q+1}}\Big(\sum_{n=1}^{q+1}|c|_{{\mathcal{U}}%
,n}+N_{q+2}(c,M_{p}^{2})\Big)
\end{equation*}%
and the statement is proved. $\square $

\medskip

As an immediate consequence of Proposition \ref{S-S} and \ref{S-L} we obtain:

\begin{proposition}
\label{S-S1} We assume that (\ref{nss3'}) holds and, for $p\geq 2$, set $M_p$
as in (\ref{nss4'}). For every $q\geq 2$ there exists a constant $C\geq 1$
depending on $p,q$ such that
\begin{equation}
\left\Vert \left\vert S(c,Z)\right\vert \right\Vert _{q,p}\equiv \Vert
S(c,Z)\Vert _{q,p}+\Vert LS(c,Z)\Vert _{q-2,p}\leq \frac{C}{r^{q-1}}\Big(%
\sum_{n=1}^{q-1}|c|_{{\mathcal{U}},n}+N_{q}(c,M_{p}^{2})\Big).  \label{nss15}
\end{equation}
\end{proposition}

We conclude this section with a result concerning the non degeneracy of the
Malliavin covariance matrix of $S(c,Z).$ Actually we are not able to obtain
such estimates for general series but for finite series only.
\begin{lemma}
\label{crucial-est}
For $N\in\N$, set
\begin{equation*}
S_{N}(c,Z)=\sum_{1\leq \left\vert \alpha \right\vert \leq N}c(\alpha
)Z^{\alpha }\quad\mbox{and}\quad i_{N}(c)=\sum_{m=1}^{N}\sum_{\alpha \in \Gamma
_{m}}m!c(\alpha )^{2}.
\end{equation*}
Then for every $p\geq 1$ such that $\sup_{k}\|Z_k\|_{2p}<\infty$ there exists a universal constant $C_{p}$ such that for every $\eta \leq \frac{1}{2}m(r)i_{N}(c)$ we have%
\begin{equation}
{\mathbb{P}}\big(\sigma _{S_{N}(c,Z)}\leq \eta \big)\leq \Big(\frac{C_p (1+i_N(c))}{%
m(r)i_{N}(c)}\, (N!)^32^NN^{-5/4}(\overline{\kappa }_{N}(c)+\overline{\delta}_N(c))\Big)^{p}  \label{Cov}
\end{equation}%
with
\begin{equation}
\overline{\kappa }_{N}(c)=\sum_{l=1}^{N}\kappa _{4,l}^{1/4}(c)
\quad \mbox{and}\quad
\overline{\delta }_{N}(c)=\sum_{l=1}^{N}\delta_{l}(c).  \label{Cov1}
\end{equation}
\end{lemma}

\textbf{Proof.} We write%
\begin{align*}
\sigma _{S_{N}(c,Z)}& =\sum_{j=1}^{\infty }\left\vert \partial
_{j}S_{N}(c,Z)\right\vert ^{2}\chi _{j}=\sum_{j=1}^{\infty }\left\vert
\partial _{j}S_{N}(c,Z)\right\vert ^{2}\widetilde{\chi }_{j}+m(r)%
\sum_{j=1}^{\infty }\left\vert \partial _{j}S_{N}(c,Z)\right\vert ^{2} \\
& =\sum_{j=1}^{\infty }\left\vert \partial _{j}S_{N}(c,Z)\right\vert ^{2}%
\widetilde{\chi }_{j}+m(r)\Big(\sum_{j=1}^{\infty }\left\vert \partial
_{j}S_{N}(c,Z)\right\vert ^{2}-i_{N}(c)\Big)+m(r)i_N(c).
\end{align*}%
We set
\begin{equation*}
\tilde{I}_{N}(c,Z)=\sum_{j=1}^{\infty }\left\vert \partial
_{j}S_{N}(c,Z)\right\vert ^{2}\widetilde{\chi }_{j}\quad \mbox{and}\quad
I_{N}(c,Z)=\sum_{j=1}^{\infty }\left\vert \partial _{j}S_{N}(c,Z)\right\vert^2.
\end{equation*}%
For $\eta \leq \frac{m(r)i(c)}{2}$, it follows that%
\begin{align*}
{\mathbb{P}}\big(\sigma _{S_{N}(c,Z)}\leq \eta )=& {\mathbb{P}}\Big(\tilde{I}%
_{N}(c,Z)+m(r)\big(I_{N}(c,Z)-i_N(c)\big)\leq -\frac{m(r)i_N(c)}{2}\Big) \\
\leq & {\mathbb{P}}\Big(\tilde{I}_{N}(c,Z)\leq -\frac{m(r)i_N(c)}{4}\Big)+{%
\mathbb{P}}\Big(m(r)\big(I_{N}(c,Z)-i_N(c)\big)\leq -\frac{m(r)i_N(c)}{4}\Big) \\
\leq & {\mathbb{P}}\Big(|\tilde{I}_{N}(c,Z)|\geq \frac{m(r)i_N(c)}{4}\Big)+{%
\mathbb{P}}\Big(|\big(I_{N}(c,Z)-i_N(c)|\geq \frac{i_N(c)}{4}\Big) \\
\leq & \Big(\frac{4}{m(r)i_N(c)}\Big)^{p}\big(\Vert \tilde{I}_{N}(c,Z)\Vert
_{p}^{p}+\Vert I_{N}(c,Z)-i_{N}(c)\Vert _{p}^{p})
\end{align*}%
Now, the real difficulty is to produce $L^{p}$ estimates for $\tilde{I}%
_{N}(c,Z)$ and $I_{N}(c,Z)-i_{N}(c)$. Section \ref{app-Lp} in Appendix \ref{app-B} is devoted to
such a problem, and the final result is given in Lemma \ref{Est-Ttilde}. So,
we use (\ref{NL0}) and the statement immediately follows. $\square $

\section{Convergence in total variation}\label{conv2}

In this section we study the convergence of stochastic series to the
Gaussian law in two situations: first, we consider finite series and we obtain estimates of the error which are not asymptotic; in a second stage we deal with infinite
series and we prove a convergence result, but in this case we are no more able to get the rate of convergence.

\subsection{Error estimates for finite series}

The aim of this section is to obtain non asymptotic estimates for the invariance principle in total variation distance. We stress that here $N$
is finite and fixed and that we consider a fixed set of coefficients $%
c=(c(\alpha ))_{\alpha }$ - so the results is not asymptotic. The estimates
will be given in terms of $\overline{\kappa }_{N}(c)$ defined in (\ref%
{Law6v}) and of \ $\overline{\delta }_{N}(c)=\sum_{l=1}^{N}l!\delta _{l}(c)$%
\ with $\delta _{N}(c)$ defined in (\ref{Law3''}). We will use the
normalization hypothesis%
\begin{equation}
i_{N}(c)=\sum_{m=1}^{N}\sum_{\alpha \in \Gamma _{m}}m!c(\alpha )^{2}=1
\label{A8}
\end{equation}%
and%
\begin{equation}
\sum_{\ell =1}^{N}\ell !\ell \delta _{\ell }(c)^{2}\leq \frac{i_{N}(c)}{4}=%
\frac{1}{4}  \label{A9}
\end{equation}%
Given two sequences $(Z_{k})_{k\in {\mathbb{N}}}$ and $(\overline{Z}%
_{k})_{k\in {\mathbb{N}}}$ we denote%
\begin{equation*}
M_{p}=M_{p}(Z,\overline{Z})=\max_{k}\left\Vert Z_{k}\right\Vert _{p}\vee
\left\Vert \overline{Z}_{k}\right\Vert _{p}.
\end{equation*}

The main result in this section is the following:

\begin{theorem}
\label{TV} Let $(Z_{k})_{k\in {\mathbb{N}}}$ and $(\overline{Z}_{k})_{k\in {%
\mathbb{N}}}$ satisfy (\ref{nss3'}) and such that $Z_{k}\in \mathcal{L}%
(r,\varepsilon )$, $\overline{Z}_{k}\in \mathcal{L}(r,\varepsilon )$ and $M_{p}=M_{p}(Z,\overline{Z})<\infty $ for every $p\geq 1$. We also assume that
(\ref{A8}) and (\ref{A9}) hold true. Then for every $p_*\geq 1$ there exist positive constant $C_*,d_*,c_*$ and $M_*$ such that for every $N\in \N$%
\begin{equation}\label{A10}
\left\vert {\mathbb{E}}(f(S_{N}(c,Z)))-{\mathbb{E}}(f(S_{N}(c,\overline{Z}%
)))\right\vert \leq
\frac{C_\ast}{(m(r)r)^{d_*}}\,M_*^NN^{c_*}(N!)^{3p_*}\left\Vert f\right\Vert _{\infty}
(\overline{\kappa }_{N}(c)^{p_*}+\overline{\delta}_{N}(c)).
\end{equation}
We stress that all constants depend on the random sequences $(Z_{k})_{k\in {\mathbb{N}}}$ and $(\overline{Z}_{k})_{k\in {%
\mathbb{N}}}$ only through $M_{p}=M_{p}(Z,\overline{Z})$ for a suitably large $p$.
\end{theorem}

\textbf{Proof.} We first give some estimates which
are specific to finite series (under the hypothesis that $c(\alpha )=0$ if $%
\left\vert \alpha \right\vert \geq N)$. First of all, for $q\geq 1,$
\begin{eqnarray}
N_{q}(c,M) &=&\Big(\sum_{m=q}^{N}M^{m-q}\times \frac{(m!)^{2}}{(m-q)!}%
\sum_{\alpha \in \Gamma _{m}}c(\alpha )^{2}\Big)^{1/2}\leq M^{(N-q)/2}N^{q/2}%
\Big(\sum_{m=q}^{N}m!\sum_{\alpha \in \Gamma _{m}}c(\alpha )^{2}\Big)^{1/2}
\label{A1} \\
&\leq &M^{(N-q)/2}N^{q/2}  \notag
\end{eqnarray}%
the last inequality being true if $i_{N}(c)\leq 1.$ As an immediate
consequence of this and of (\ref{nss15}), for every $p\geq 1$
\begin{equation}
\left\Vert \left\vert S_{N}(c,Z)\right\vert \right\Vert _{q,p}\leq \frac{%
C_{p}}{r^{q-1}}\times M_p^{2N}N^{q/2}.  \label{A2}
\end{equation}%
Moreover let $\tilde{c}^{k}(\alpha )=c(k,\alpha ).$ Then
\begin{eqnarray}
\sum_{k=1}^{\infty }N_{q}^{3}(\tilde{c}^{k},M) &\leq
&M^{3(N-q)/2}N^{3q/2}\sum_{k=1}^{\infty }\Big(\sum_{m=q}^{N}m!\sum_{\alpha
\in \Gamma _{m}}c(k,\alpha )^{2}\Big)^{3/2}  \label{A3} \\
&\leq &M^{3(N-q)/2}N^{3q/2}\overline{\delta }_{N}(c)\sum_{k=1}^{\infty }%
\Big(\sum_{m=q}^{N}m!\sum_{\alpha \in \Gamma _{m}}c(k,\alpha )^{2}\Big)  \notag
\\
&\leq &M^{3(N-q)/2}N^{3q/2}\overline{\delta }_{N}(c).  \notag
\end{eqnarray}%
We are now able to start the proof itself.

\smallskip

\textbf{Step 1}. Let $f\,:\,{\mathbb{R}}\rightarrow {\mathbb{R}}$ be such
that $\Vert f\Vert _{\infty }\leq 1$. For $\delta >0$, let $f_{\delta }$
denote its regularization, as in (\ref{FD17''}). Then we have
\begin{equation}
\left\vert {\mathbb{E}}(f(S_{N}(c,Z)))-{\mathbb{E}}(f(S_{N}(c,\overline{Z}%
)))\right\vert \leq a_{N}(\delta )+b_{N}(\delta )+\overline{b}_{N}(\delta )
\label{c1}
\end{equation}%
in which
\begin{equation*}
\begin{array}{c}
a_{N}(\delta )=|{\mathbb{E}}(f_{\delta }(S_{N}(c,Z)))-{\mathbb{E}}(f_{\delta
}(S_{N}(c,\overline{Z})))|,\smallskip \\
b_{N}(\delta )=|{\mathbb{E}}(f(S_{N}(c,Z)))-{\mathbb{E}}(f_{\delta
}(S_{N}(c,Z)))|,\qquad \overline{b}_{N}(\delta )=|{\mathbb{E}}(f(S_{N}(c,%
\overline{Z})))-{\mathbb{E}}(f_{\delta }(S_{N}(c,\overline{Z})))|.%
\end{array}%
\end{equation*}%
So, we study separately such contributions.

\smallskip

\textbf{Step 2: estimate of $a_{N}(\delta )$.} We use some facts developed in the proof of Theorem \ref{Conv-smooth}. Let $J\in {\mathbb{N}}$ and%
\begin{equation*}
S_{N,J}(c,Z)=\sum_{n=1}^{N}\sum_{\alpha \in \Gamma _{n}(J)}c(\alpha
)Z^{\alpha }.
\end{equation*}%
We also set
\begin{equation*}
a_{N,J}(\delta )=|{\mathbb{E}}(f_{\delta }(S_{N,J}(c,Z)))-{\mathbb{E}}%
(f_{\delta }(S_{N,J}(c,\overline{Z})))|.
\end{equation*}%
Since the estimate for $a_{N,J}$ will not depend on $J$, we will get the
result for $a_{N}(\delta )$ by passing to the limit. So, we recall the
following facts.

For $\theta \in (0,1)$, in Theorem \ref{Conv-smooth} we have denoted%
\begin{equation*}
\widehat{Z}^{k}(\theta )=(Z_{1},...,Z_{k-1},\theta Z_{k},\overline{Z}%
_{k+1},...,\overline{Z}_{J}),\quad k=0,1,\ldots ,J,
\end{equation*}%
with $\widehat{Z}^{0}(\theta )=(\overline{Z}_{1},\ldots ,\overline{Z}_{J})$
and $\widehat{Z}^{J}(\theta )=(Z_{1},...,Z_{J-1},\theta Z_{J})$. Moreover we
have denoted%
\begin{equation*}
s_{k}=S_{N,J}(c,\widehat{Z}^{k}(0))\qquad \mbox{and}\qquad
v_{k}=c(k)+\sum_{n=2}^{N}\sum_{\substack{ \beta \in \Gamma _{n-1}(J)  \\ %
k\notin \beta }}c(k,\beta )(\widehat{Z}^{k}(0))^{\beta }
\end{equation*}%
so that%
\begin{equation*}
S_{N,J}(c,\widehat{Z}^{k}(\theta ))=s_{k}+\theta Z_{k}v_{k}
\end{equation*}%
and we have proved that
\begin{align*}
{\mathbb{E}}(f_{\delta }(S_{N,J}(c,Z)))-{\mathbb{E}}(f_{\delta }(S_{N,J}(c,%
\overline{Z}))& =\frac{1}{6}\sum_{k=1}^{J}\int_{0}^{1}{\mathbb{E}}(f_{\delta
}^{\prime \prime \prime }(s_{k}+\theta Z_{k}v_{k})Z_{k}^{3}v_{k}^{3}))\theta
d\theta \\
& \quad -\frac{1}{6}\sum_{k=1}^{J}\int_{0}^{1}{\mathbb{E}}(f_{\delta
}^{\prime \prime \prime }(s_{k}+\theta \overline{Z}_{k}v_{k})\overline{Z}%
_{k}^{3}v_{k}^{3}))\theta d\theta .
\end{align*}%
In the original proof of Theorem \ref{Conv-smooth} we upper bounded $%
f_{\delta }^{\prime \prime \prime }(s_{k}+\theta Z_{k}v_{k})$ by $\left\Vert
f_{\delta }^{\prime \prime \prime }\right\Vert _{\infty }$ but now we will
use integration by parts in order to get rid of the derivatives. In order to
do it we will use Lemma \ref{L2} with $F_{k}=s_{k}+\theta Z_{k}v_{k}$ and $%
G_{k}=Z_{k}^{3}v_{k}^{3}$, $M=1$ and $q=3$. Then we apply (\ref{FD26}) with $%
p_{1}=p_{2}=p_{3}=3$: for every $\eta _{k}>0$ we obtain%
\begin{equation}
\left\vert {\mathbb{E}}(f_{\delta }^{\prime \prime \prime }(s_{k}+\theta
Z_{k}v_{k})Z_{k}^{3}v_{k}^{3}))\right\vert \leq C\left\Vert f\right\Vert
_{\infty }\lambda _{\delta ,\eta _{k},9}^{3}(s_{k}+\theta
Z_{k}v_{k})(1+\left\Vert \left\vert s_{k}+\theta Z_{k}v_{k}\right\vert
\right\Vert _{1,5,36}^{12})\left\Vert Z_{k}^{3}v_{k}^{3}\right\Vert _{3,9}
\label{A4}
\end{equation}%
with (see (\ref{FD17}))
\begin{equation*}
\lambda _{\delta ,\eta _{k},9}(s_{k}+\theta Z_{k}v_{k})=\frac{1}{\delta }{%
\mathbb{P}}(\lambda _{s_{k}+\theta Z_{k}v_{k}}\leq \eta _{k})^{\frac{1}{9}}+%
\frac{1}{\eta _{k}}.
\end{equation*}%
We denote
\begin{equation*}
c^{k,\theta }(\alpha )=c(\alpha )(1_{\{k\notin \alpha \}}+\theta 1_{\{k\in
\alpha \}}).
\end{equation*}%
Then
\begin{equation*}
s_{k}+\theta Z_{k}v_{k}=S_{N,J}(c,\widehat{Z}^{k}(\theta
))=S_{N,J}(c^{k,\theta },\widehat{Z}^{k}(1)).
\end{equation*}%
So, we take $\eta _{k}=m(r)i_{m,J}(c^{k,\theta })$, with $%
i_{N,J}(c^{k,\theta })=\sum_{\ell =1}^{N}\sum_{\alpha \in \Gamma _{\ell
}(J)}\ell !c^{k,\theta }(\alpha )^{2}$, and we use (\ref{Cov}) with $\bar{p}%
\geq 1$: we get
\begin{equation*}
{\mathbb{P}}(\lambda _{s_{k}+\theta Z_{k}v_{k}}\leq \eta _{k})\leq \Big(%
\frac{C_{p,N}}{m(r)i_{N,J}(c^{k,\theta })}\Big)^{\bar{p}}\overline{\kappa }%
_{N}(c^{k,\theta })^{\bar{p}}
\end{equation*}%
for every $\theta \in (0,1)$, where $C_{p,N}=C_p\mathcal{D}_N$ is given in Lemma \ref{crucial-est}:
$$
\mathcal{D}_N=(N!)^32^NN^{-5/4}.
$$
 Now, observe that
\begin{align*}
i_{N,J}(c^{k,\theta })& =\sum_{\ell =1}^{N}\sum_{\alpha \in \Gamma _{\ell
}(J)}\ell !c^{k,\theta }(\alpha )^{2}=\sum_{\ell =1}^{N}\sum_{\alpha \in
\Gamma _{\ell }(J)}\ell !c(\alpha )^{2}-(1-\theta )\sum_{\ell
=1}^{N}\sum_{\alpha \in \Gamma _{\ell }(J)}\ell !c(\alpha )^{2}%
\mbox{\large
\bf 1}_{k\in \alpha } \\
& \geq i_{N,J}(c)-(1-\theta )\sum_{\ell =1}^{N}\ell !\ell \delta _{\ell
}(c)^{2}.
\end{align*}%
Under (\ref{A9}),  for every $\theta \in (0,1)$ and for every $J$ large, one can write
\begin{equation*}
i_{N,J}(c^{k,\theta })\geq \frac{i_{N}(c)}{2}-\frac{i_{N}(c)}{4}=\frac{1}{4}.
\end{equation*}%
One also has
\begin{equation*}
\overline{\kappa }_{N}(c^{k,\theta })\leq 2\overline{\kappa }_{N}(c)\quad\mbox{and}\quad
\overline{\delta }_{N}(c^{k,\theta })\leq 2\overline{\delta }_{N}(c),
\end{equation*}%
so finally%
\begin{equation*}
{\mathbb{P}}(\lambda _{s_{k}+\theta Z_{k}v_{k}}\leq \eta _{k})\leq \Big(%
\frac{C_{\bar{p},N}}{m(r)}\Big)^{\bar{p}}(\overline{\kappa }_{N}(c)+\overline{\delta }_{N}(c))^{\bar{p}}.
\end{equation*}%
Moreover, for every $k\in {\mathbb{N}}$, $\theta \in (0,1)$ and for every $J$
large we have
\begin{equation*}
\frac{1}{\eta }_{k}\leq \frac{16}{m(r)}
\end{equation*}%
and then
\begin{equation*}
\lambda _{\delta ,\eta _{k},9}(s_{k}+\theta Z_{k}v_{k})=\frac{1}{\delta }%
\Big(\frac{C_{\bar{p},N}}{m(r)}\Big)^{\frac{\bar{p}}{9}}\overline{\kappa }_{N}(c)^{%
\frac{\bar{p}}{9}}+\frac{16}{m(r)}.
\end{equation*}%
Now we choose
\begin{equation}
\delta=\delta _{\overline{p}}=\Big(\frac{C_{\bar{p},N}}{m(r)}\Big)^{\frac{\bar{p}}{9}%
}(\overline{\kappa }_{N}(c)+\overline{\delta }_{N}(c))^{\frac{\bar{p}}{9}}  \label{A5}
\end{equation}%
and then%
\begin{equation*}
\lambda _{\delta _{\overline{p}},\eta _{k},9}(s_{k}+\theta Z_{k}v_{k})\leq 1+%
\frac{16}{m(r)}.
\end{equation*}%
Coming back to (\ref{A4}), we have
\begin{equation}
\left\vert {\mathbb{E}}(f_{\delta _{\overline{p}}}^{\prime \prime \prime
}(s_{k}+\theta Z_{k}v_{k})Z_{k}^{3}v_{k}^{3}))\right\vert \leq \frac{C}{%
m^{3}(r)}\left\Vert f\right\Vert _{\infty }(1+\left\Vert \left\vert
s_{k}+\theta Z_{k}v_{k}\right\vert \right\Vert
_{1,5,36}^{12})\left\Vert Z_{k}^{3}v_{k}^{3}\right\Vert _{3,9}.
\label{A6}
\end{equation}
We recall that $s_{k}+\theta Z_{k}v_{k}=S_{N,J}(c^{k,\theta },\widehat{Z}%
^{k}(1))$ and we use (\ref{A2}) with $c$ replaced by $c^{k,\theta }$ and we
obtain
\begin{equation*}
\left\Vert \left\vert s_{k}+\theta Z_{k}v_{k}\right\vert \right\Vert
_{1,5,36}\leq \frac{C}{r^{4}}M_{36}^{N}N^{5/2}.
\end{equation*}%
Since $Z_{k}$ and $v_{k}$ are independent we have
\begin{equation*}
\left\Vert Z_{k}^{3}v_{k}^{3}\right\Vert _{3,9}=\left\Vert
Z_{k}^{3}\right\Vert _{3,9}\left\Vert v_{k}^{3}\right\Vert _{3,9}\leq
CM_{27}^{3}\left\Vert v_{k}\right\Vert _{3,27}^{3}
\end{equation*}%
the last inequality being true because of (\ref{FD4'}). Summing on $k$ in (\ref%
{A6}) we get
\begin{equation*}
\sum_{k=1}^{J}\left\vert {\mathbb{E}}(f_{\delta _{\overline{p}}}^{\prime
\prime \prime }(s_{k}+\theta Z_{k}v_{k})Z_{k}^{3}v_{k}^{3}))\right\vert \leq
\frac{CM_{36}^{12N}N^{30}}{m^{3}(r)r^{48}}\left\Vert f\right\Vert _{\infty
}\sum_{k=1}^{J}\left\Vert v_{k}\right\Vert _{3,27}^{3}.
\end{equation*}%
Since $v_{k}=c(k)+S_{m-1,J}(\tilde{c}^{k},\widehat{Z}^{k}(0))$ with $\tilde{c%
}^{k}(\alpha )=c(k,\alpha )$, we use (\ref{nss8}) and (\ref{A3})\ and we
obtain
\begin{eqnarray*}
\sum_{k=1}^{J}\left\Vert v_{k}\right\Vert _{3,27}^{3} &\leq
&C\sum_{k=1}^{J}(\sum_{l=0}^{3}N_{l}(\tilde{c}^{k},M_{27}^{2}))^{3} \\
&\leq &C\sum_{l=0}^{3}\sum_{k=1}^{J}N_{l}^{3}(\tilde{c}^{k},M_{27}^{2}) \\
&\leq &CM_{27}^{3N}N^{9}\overline{\delta }_{N}(c).
\end{eqnarray*}%
Inserting in the previous inequality we get%
\begin{equation*}
\sum_{k=1}^{J}\left\vert {\mathbb{E}}(f_{\delta _{\overline{p}}}^{\prime
\prime \prime }(s_{k}+\theta Z_{k}v_{k})Z_{k}^{3}v_{k}^{3}))\right\vert \leq
\frac{CM_{27}^{15N}N^{39}}{m^{3}(r)r^{48}}\left\Vert f\right\Vert
_{\infty }\overline{\delta }_{N}(c).
\end{equation*}%
Since the same estimate holds with $Z$ replaced by $\overline{Z}$ we
conclude that%
\begin{equation}
a_{N}(\delta _{\overline{p}})\leq
\frac{CM_{27}^{15N}N^{39}}{m^{3}(r)r^{48}}\left\Vert f\right\Vert _{\infty }\overline{\delta }_{N}(c).
\label{A7}
\end{equation}%
\textbf{Step 3: estimate of $b(\delta _{\overline{p}})$ and $\overline{b}%
(\delta _{\overline{p}})$.} We use the regularization inequality (\ref{FD24'}%
) in Lemma \ref{L3} with $\eta =\frac{1}{2}m(r)i_{N}(c)=\frac{1}{2}m(r)$
and we obtain:
\begin{align*}
b_N(\delta)& \leq C_{\ast }\left\Vert f\right\Vert _{\infty }%
\Big({\mathbb{P}}\Big(\lambda _{S_{N}(c,Z)}\leq \frac{m(r)}{2}\Big)+
\frac{\sqrt{\delta}}{m(r)^{l_{\ast }}}\big(1+\left\Vert \left\vert
S(c,Z)\right\vert \right\Vert _{3,l_{\ast }}\big)^{a_{\ast }}\Big)
\end{align*}%
where $C_{\ast },l_{\ast },a_{\ast }$ are universal constants. We take now $\delta=\delta _{\overline{p}}$ as in (\ref{A5}) and we use  (\ref{Cov}), so that
\begin{align*}
b_N(\delta _{\overline{p}})
& \leq C_\ast\left\Vert f\right\Vert _{\infty }\Big(\Big(\frac{C_{\widehat{p}}\mathcal{D}_N}{m(r)}(\overline{\kappa }_{N}(c)+\overline{\delta}_{N}(c))\Big)^{%
\widehat{p}}
+\frac{(\frac{C_{\bar{p}}\mathcal{D}_N}{m(r)}(\overline{\kappa }_{N}(c)+\overline{\delta}_{N}(c)))^{\frac{\bar p}{18}}}{m(r)^{l_\ast}}
\big(1+\left\Vert \left\vert
S_N(c,Z)\right\vert \right\Vert _{3,l_{\ast }}\big)^{a_{\ast }}
\Big).
\end{align*}%
By applying  (\ref{A2}) and taking $\widehat{p}=p_{\ast }$ and $\overline{p}=18p_{\ast }$, with $p_*\geq 1$, we get
\begin{align*}
b_N(\delta _{\overline{p}})
& \leq C_\ast\left\Vert f\right\Vert _{\infty}
\Big(\frac{\mathcal{D}_N}{m(r)}(\overline{\kappa }_{N}(c)+\overline{\delta}_{N}(c))\Big)^{p_*}
\times \Big(
\frac{1}{r^{2}}\times M_{l_*}^{N}N^{3/2}\Big)^{a_*}\\
& \leq \frac{C_\ast}{(m(r)r)^{d_*}}\,M_*^NN^{c_*}\mathcal{D}_N^{p_*}\left\Vert f\right\Vert _{\infty}
(\overline{\kappa }_{N}(c)+\overline{\delta}_{N}(c))^{p_*}.
\end{align*}%
The same estimate holds for $\overline{b}(\delta _{\overline{p}}).$
This together with (\ref{A7}) and (\ref{c1}) yield (\ref{A10}).
$\square $

\subsection{\protect\bigskip Gaussian limit}

In this section we estimate the total variation distance between $S_{N}(c,Z)$
and a standard normal distributed random variable $G.$ This will be an
immediate consequence of the result from the previous section and the
following theorem due to Nourdin and Peccati \cite{[NP1]}.

\begin{theorem}
\label{G}Let $(\overline{Z}_{k})_{k\in {\mathbb{N}}}$ with $\overline{Z}_{k}$
standard normal random variables and let $G$ be another standard normal
random variable. Suppose that (\ref{A8}) holds. There exists an universal
constant $C$ such that for every $N$ and every measurable and bounded
function $f$%
\begin{equation}
\left\vert {\mathbb{E}}(f(S_{N}(c,\overline{Z})))-{\mathbb{E}}%
(f(G))\right\vert \leq 3\left\Vert f\right\Vert _{\infty
}N^{3}(2N)!N!^{3}\sum_{l=0}^{N}\kappa _{4,l}^{1/4}(c).  \label{G1}
\end{equation}%
Moreover, if $c(\alpha )=0$ for $\left\vert \alpha \right\vert =1,$ then
\begin{equation}
\left\vert {\mathbb{E}}(f(S_{N}(c,\overline{Z})))-{\mathbb{E}}%
(f(G))\right\vert \leq 3\left\Vert f\right\Vert _{\infty
}N^{3}(2N)!N!^{3}\sum_{l=0}^{N}\kappa _{4,l}^{1/2}(c).  \label{G2}
\end{equation}
\end{theorem}

\textbf{Proof}. The proof of (\ref{G2}) is an immediate consequence of the results in \cite{[NP1]}, see (3.38) in Theorem 3.1  and Proposition 3.7 therein. But in order to
obtain (\ref{G1}) we have to complete the argument from \cite{[NP1]}. Since
the argument is essentially the same we just sketch the proof and in
particular we explain why $\kappa _{4,l}^{1/4}(c)$ appears instead of $%
\kappa _{4,l}^{1/2}(c)$. Let us briefly recall the notations from \cite%
{[NP1]}. For a symmetric kernel $\phi _{n}\in L^1(\R_{+}^{n})$ one denotes by $%
I_{n}(\phi _{n})$ the multiple stochastic integral with kernel $\phi _{n}$.
This is an element of the Wiener space and the Malliavin derivative and the
Ornstein operator for it are defined as%
\begin{equation}
D_{s}I_{n}(\phi _{n})=nI_{n-1}(\phi _{n}(\circ ,s),\quad LI_{n}(\phi
_{n})=-nI_{n}(\phi _{n}).  \label{G3}
\end{equation}%
Consider now a functional $F_{N}=\sum_{n=1}^{N}I_{n}(\phi _{n}).$ The
operators $DF_{N}$ and $LF_{N}$ extend by linearity. Now, (3.38) in \cite{[NP1]} says that
\begin{equation}
\left\vert {\mathbb{E}}(f(F_{N}))-{\mathbb{E}}(f(G))\right\vert \leq
2\left\Vert f\right\Vert _{\infty }(\E((1-\left\langle
DF_{N},-DL^{-1}F_{N}\right\rangle )^{2}))^{1/2}.  \label{G4}
\end{equation}%
So our aim now is to estimate the quantity in the right hand side of (\ref%
{G3}). This is done in Proposition 3.7 from \cite{[NP1]} but there one
considers multiple integrals $I_{n}(\phi _{n})$ with $n\geq 2$ only. If $%
I_{1}(\phi _{1})$ comes in also, one more term appears and we explain this
now. Following \cite{[NP1]} we use (\ref{G3}) and we obtain%
\begin{eqnarray*}
\left\langle DF_{N},-DL^{-1}F_{N}\right\rangle
&=&n\int_{\R_{+}}I_{n-1}(\phi _{n}(\circ ,s))I_{m-1}(\phi
_{m}(\circ ,s))ds \\
&=&\sum_{n,m=1}^{N}n\int_{\R_{+}}I_{n-1}(\phi _{n}(\circ ,s))I_{m-1}(\phi
_{m}(\circ ,s))ds \\
&=&A+\left\Vert \phi _{1}\right\Vert _{L^{2}(\R_{+})}^{2}+B+B^{\prime }
\end{eqnarray*}%
with%
\begin{align*}
&A=\sum_{n,m=2}^{N}n\int_{\R_{+}}I_{n-1}(\phi _{n}(\circ ,s))I_{m-1}(\phi
_{m}(\circ ,s))ds,\\
&B =\sum_{m=2}^{N}\int_{\R_{+}}\phi _{1}(s)I_{m-1}(\phi _{m}(\circ
,s))ds=\sum_{m=2}^{N}I_{m}(\phi _{1}\otimes _{1}\phi _{m}), \\
&B^{\prime } =\sum_{n=2}^{N}n\int_{\R_{+}}I_{n-1}(\phi _{n}(\circ ,s))\phi
_{1}(s))ds=\sum_{n=2}^{N}nI_{n}(\phi _{1}\otimes _{1}\phi _{n}).
\end{align*}%
Using the product formula for multiple stochastic integrals (see (2.29) in
\cite{[NP1]} for this formula) one obtains
\begin{eqnarray*}
A &=&\sum_{n,m=2}^{N}n\sum_{r=0}^{n\wedge m-1}r!\left(
\begin{tabular}{c}
$n-1$ \\
$r$%
\end{tabular}%
\right) \left(
\begin{tabular}{c}
$m-1$ \\
$r$%
\end{tabular}%
\right) I_{n+m-2-2r}\Big(\int_{\R_{+}}\phi _{n}(\circ ,s)\widetilde{\otimes }%
_{r}\phi _{m}(\circ ,s)ds\Big) \\
&=&\sum_{n,m=2}^{N}n\sum_{r=0}^{n\wedge m-1}r!\left(
\begin{tabular}{c}
$n-1$ \\
$r$%
\end{tabular}%
\right) \left(
\begin{tabular}{c}
$m-1$ \\
$r$%
\end{tabular}%
\right) I_{n+m-2-2r}(\phi _{n}\widetilde{\otimes }_{r+1}\phi _{m}) \\
&=&\sum_{n,m=2}^{N}n\sum_{r=1}^{n\wedge m}(r-1)!\left(
\begin{tabular}{l}
$n-1$ \\
$r-1$%
\end{tabular}%
\right) \left(
\begin{tabular}{l}
$m-1$ \\
$r-1$%
\end{tabular}%
\right) I_{n+m-2r}(\phi _{n}\widetilde{\otimes }_{r}\phi _{m}) \\
&=&\sum_{n=2}^{N}n!\left\Vert \phi _{n}\right\Vert
_{L^{2}(\R_{+}^{n})}^{2}+A^{\prime }
\end{eqnarray*}%
with $A^{\prime }$ just defined by the above equality: so it represents the
sum over $(n,m,r)$ such that $(n,m,r)\neq (n,n,n).$ Notice that in this
case
\begin{eqnarray*}
\left\Vert I_{n+m-2r}(\phi _{n}\widetilde{\otimes }_{r}\phi _{m})\right\Vert
_{2} &=&(n+m)!\left\Vert \phi _{n}\widetilde{\otimes }_{r}\phi
_{m}\right\Vert _{2}\leq (n+m)!\kappa _{4,n}^{1/2}(c)\vee \kappa
_{4,m}^{1/2}(c) \\
&\leq &(n+m)!(\kappa _{4,n}^{1/2}(c)+\kappa _{4,m}^{1/2}(c))
\end{eqnarray*}
the last inequality being a consequence of well-known facts, which has been here collected in Appendix \ref{app-B}, see (\ref{acc8}) and (\ref{acc7}). So
\begin{equation*}
\left\Vert A^{\prime }\right\Vert _{2}\leq N^{3}(2N)!\times
N!^{3}\sum_{n=2}^{N}\kappa _{4,n}^{1/2}(c).
\end{equation*}%
And using (\ref{acc9-1}) we get $\left\Vert I_{m}(\phi _{1}\otimes _{1}\phi
_{m})\right\Vert _{2}=m!\left\Vert \phi _{1}\otimes _{1}\phi _{m}\right\Vert
_{2}\leq m!\left\Vert \phi _{1}\right\Vert _{2}\kappa _{4,n}^{1/4}(c)$ (we
stress that $\kappa _{4,n}^{1/4}(c)$ appears here instead of$\ \kappa
_{4,n}^{1/2}(c)$) so that%
\begin{equation*}
\left\Vert B\right\Vert _{2}+\left\Vert B^{\prime }\right\Vert _{2}\leq
2N\, N!\sum_{n=2}^{N}\kappa _{4,n}^{1/4}(c).
\end{equation*}%
We suppose now that $\sum_{n=1}^{N}n!\left\Vert \phi _{n}\right\Vert
_{L^{2}(\R_{+}^{n})}^{2}=1$ and we write $\left\langle
DF_{N},-DL^{-1}F_{N}\right\rangle =1+A^{\prime }+B+B^{\prime }$ so that
\begin{equation*}
\left\Vert 1-\left\langle DF_{N},DL^{-1}F_{N}\right\rangle \right\Vert
_{2}\leq \left\Vert A^{\prime }\right\Vert _{2}+\left\Vert B\right\Vert
_{2}+\left\Vert B^{\prime }\right\Vert _{2}\leq 3N^{3}(2N)!\times
N!^{3}\sum_{n=2}^{N}\kappa _{4,n}^{1/4}(c)
\end{equation*}%
and using (\ref{G4}) this gives
\begin{equation}
\left\vert {\mathbb{E}}(f(F_{N}))-{\mathbb{E}}(f(G))\right\vert \leq
3\left\Vert f\right\Vert _{\infty }N^{3}(2N)!\times
N!^{3}(\sum_{n=2}^{N}\kappa _{4,n}^{1/4}(c)).  \label{G5}
\end{equation}%
Of course, $\kappa _{4,n}^{1/4}(c)$ can be replaced by $\kappa _{4,n}^{1/2}(c)$ if $\phi_1=0$.

We come now back to stochastic series and we define $f_{c_{(m)}}$ to be the
kernels which correspond to the coefficients $c(\alpha ):$ $%
f_{c_{(m)}}(t_{1},...,t_{m})=c(\alpha )$ if $t_{i}\in \lbrack \alpha
_{i},\alpha _{i}+1),i=1,...,m$. Then $F_{N}=S_{N}(c,\overline{Z})$ and the
hypothesis (\ref{A8}) says that $\sum_{n=1}^{N}n!\Vert
f_{c_{(n)}}\Vert _{L^{2}(\R_{+}^{n})}^{2}=1.$ And using (\ref{G5}) we
obtain (\ref{G1}). $\square $

\medskip

The main result in this section is the following:

\begin{theorem}
\label{TV copy(2)} Let $(Z_{k})_{k\in {\mathbb{N}}}$ satisfy (\ref{nss3'})
and such that $Z_{k}\in \mathcal{L}(r,\varepsilon ).$ We also assume that $\sup_k\|Z_k\|_p<\infty $ for every $p\geq 1$ and we suppose that (\ref{A8})
and (\ref{A9}) hold true. There exist some constants $C,a\geq 1$ such that for every $N\in \N$ and every bounded and measurable function $f$ one
has%
\begin{equation}
\left\vert {\mathbb{E}}(f(S_{N}(c,Z)))-{\mathbb{E}}(f(G))\right\vert \leq C\frac{\left\Vert f\right\Vert _{\infty }}{(rm(r))^a}%
N^3(2N)!N!^{3}\big(\overline{\kappa} _{N}(c)+\overline{\delta}_N(c)\big).  \label{G6}
\end{equation}%
As a consequence, there exist $C,a\geq 1$ such that for every $N\in \N$ and every bounded and measurable function $f$ one has
\begin{equation}
\left\vert {\mathbb{E}}(f(S_{N}(c,Z)))-{\mathbb{E}}(f(G))\right\vert \leq C\frac{\left\Vert f\right\Vert _{\infty }}{(rm(r))^a}%
N^3(2N)!N!^{3}(1+\alpha_N^{-1}(c))\overline{\kappa}_N(c),  \label{G6bis}
\end{equation}%
in which $\alpha_N(c)=\min_{m\leq N}|c|_m\,1_{|c|_m> 0}$.
\end{theorem}

\textbf{Proof}. We take a sequence $\overline{Z}_k$, $k\in\N$, of standard normal r.v.'s and we write
$$
\left\vert {\mathbb{E}}(f(S_{N}(c,Z)))-{\mathbb{E}}(f(G))\right\vert
\leq
\left\vert {\mathbb{E}}(f(S_{N}(c,Z)))-{\mathbb{E}}(f(S_{N}(c,\overline{Z})))\right\vert +\left\vert {\mathbb{E}}(f(S_{N}(c,\overline{Z})))-{\mathbb{E}}(f(G))\right\vert.
$$
(\ref{G6}) now follows by applying Theorem \ref{TV} with $p_*=1$ and Theorem \ref{G}. Moreover, by using (\ref{Law7'}) one has $\overline{\delta}_N(c)\leq \alpha_N^{-1}(c)\,\overline{\kappa}_N(c)$, so (\ref{G6bis}) immediately follows from (\ref{G6}).
$\square$

\subsection{A convergence result for infinite series}

We consider a sequence $c^{(n)}=(c^{(n)}(\alpha ))_{\alpha }$ of coefficients
and the corresponding infinite series $S_{\infty }(c^{(n)},Z).$ Our aim is
to give sufficient conditions in order to obtain convergence to the Gaussian
law in total variation distance. Here are our hypotheses. First we assume
the normalization condition
\begin{equation}
i(c^{(n)})=\sum_{k=1}^{\infty }k!\sum_{\left\vert \alpha \right\vert
=k}c^{(n)}(\alpha )^{2}=1.  \label{CTV1}
\end{equation}%
We also assume that for every $p\geq 1$
\begin{equation}
\sup_{n}\sum_{k=0}^{3}N_{k}(c^{(n)},M_{p}^{2})<\infty  \label{CTV1a}
\end{equation}%
where $N_{k}(c^{(n)},M_{p}^{2})$ is defined in (\ref{Law3}). Moreover we
suppose that

\begin{equation}
\limsup_{N\rightarrow \infty }\limsup_{n\to\infty} %
\sum_{k\geq N}k\times k!\sum_{\left\vert \alpha \right\vert =k}\vert
c^{(n)}(\alpha )\vert ^{2}=0.  \label{CTV2}
\end{equation}%
Notice that this is analogous to the ``uniformity condition'' in (\ref{Law11}),
which is used by Hu and Nualart \cite{[HN]} for getting convergence in law for
infinite series. Then we have the following convergence result:

\begin{theorem}
\label{CTV}Let $(Z_{k})_{k\in {\mathbb{N}}}$ be a sequence of independent
centred random variables with $\E(Z_{k}^{2})=1$ and which have finite moments
of any order. Let $c^{(n)}=(c^{(n)}(\alpha ))_{\alpha }$ be a sequence of
coefficients which verify (\ref{CTV1}), (\ref{CTV1a}),  (\ref{CTV2}) and such that
\begin{equation}
\lim_{n\to\infty}\kappa _{4,m}(c^{(n)})=0\quad\mbox{and}\quad
\lim_{n\to\infty}\delta _{m}(c^{(n)})=0, \label{CTV3}
  \end{equation}%
for each $m\in \N$. Then
\begin{equation}
\lim_{n\to\infty}d_{TV}(S_{\infty }(c^{(n)},Z),G)=0  \label{CTV3a}
\end{equation}%
where $G$ is a standard Gaussian random variable and
\begin{equation}
d_{TV}(S_{\infty }(c^{(n)},Z),G):=\sup_{\left\Vert f\right\Vert _{\infty }\leq
1}\vert \E(f(S_{\infty }(c^{(n)},Z)))-\E(f(G))\vert .  \label{CTV3b}
\end{equation}
\end{theorem}

\begin{remark}
In view of (\ref{Law7'}), a sufficient condition in order that (\ref{CTV3}) holds is the following:
$$
\lim_{n\to\infty}\kappa _{4,m}(c^{(n)})\Big(1+\frac 1{|c^{(n)}|_m^2}\,1_{\{|c^{(n)}|>0\}}\Big)=0,\quad\mbox{for every }m\in\N.
$$
\end{remark}
\textbf{Proof of Theorem \ref{CTV}}.
We set $S^N(c,Z)=S_\infty(c,Z)-S_N(c,Z)$. Then, we have
\begin{align*}
\sigma_{S_\infty(c^{(n)},Z)}-\sigma_{S_N(c^{(n)},Z)}
&=|DS_\infty(c^{(n)},Z)|^2-|DS_N(c^{(n)},Z)|^2\\
&=\<DS^N(c^{(n)},Z),DS_\infty(c^{(n)},Z)+DS_N(c^{(n)},Z)\>.
\end{align*}
So, by Cauchy-Schwartz inequality
$$
\E(|\sigma_{S_\infty(c^{(n)},Z)}-\sigma_{S_N(c^{(n)},Z)}|)
\leq \|DS^N(c^{(n)},Z)\|_2 (\|DS_\infty(c^{(n)},Z)\|_2+\|DS_N(c^{(n)},Z)\|_2).
$$
Setting
$$
\varepsilon _{N}(n)=\sum_{k\geq N}k\times k!\sum_{\left\vert \alpha \right\vert =k}\vert
c^{(n)}(\alpha )\vert ^{2},
$$
by Proposition \ref{S-S}, we have
$$
\E(|\sigma_{S(c^{(n)},Z)}-\sigma_{S_N(c^{(n)},Z)}|)
\leq \varepsilon _{N}(n)\cdot 2\varepsilon _{0}(n)
$$
We take $\eta =\frac{1}{2}m(r)$ and we use (\ref{Cov}) in order to get%
\begin{eqnarray*}
\P(\sigma _{S_{\infty }(c^{(n)},Z)} &\leq &\eta )\leq \P(\vert \sigma
_{S_{\infty }(c^{(n)},Z)}-\sigma _{S_{N}(c^{(n)},Z)}\vert \geq \eta )+\P(\sigma
_{S_{N}(c^{(n)},Z)}\leq 2\eta ) \\
&\leq &\frac{1}{\eta }2\varepsilon _{0}(n)\varepsilon _{N}(n)+C_{N}(\overline{\kappa }_{N}(c^{(n)})+\overline{\delta }_{N}(c^{(n)})),
\end{eqnarray*}%
where $C_{N}$ is a constant which depends on $N$ but not on $n.$
Now we use (%
\ref{CTV3}) and we obtain, for each fixed $N$
\begin{align*}
\limsup_{n\to\infty} \P(\sigma _{S_{\infty }(c^{(n)},Z)}\leq
\eta )
&\leq \frac{1}{\eta }\limsup_{n\to\infty} 2\varepsilon_0(n)\varepsilon
_{N}(n)+C_{N}\limsup_{n\to\infty} (\overline{\kappa }%
_{N}(c^{(n)})+\overline{\delta}_N(c))\\
&=\frac{1}{\eta }\limsup_{n\to\infty} 2\varepsilon_0(n)\varepsilon
_{N}(n).
\end{align*}%
Then by  (\ref{CTV2}),
\begin{equation}
\limsup_{n\to\infty} \P(\sigma _{S_{\infty }(c^{(n)},Z)}\leq
\eta )\leq \frac{2}{\eta }\limsup_{N\rightarrow \infty }\limsup_{n\rightarrow \infty }\varepsilon_0(n)\varepsilon _{N}(n)=0.  \label{CTV5}
\end{equation}%
Now we use the regularization Lemma \ref{L3}: for every $\delta >0$
\begin{eqnarray*}
\left\vert {\mathbb{E}}(f(S_\infty(c^{(n)},Z)))-{\mathbb{E}}(f_{\delta
}(S_\infty(c^{(n)},Z)))\right\vert &\leq &C\left\Vert f\right\Vert _{\infty }\Big(%
\P(\sigma _{S_{\infty }(c^{(n)},Z)}\leq \eta )+\frac{\sqrt{\delta }}{\eta ^{p}}%
(1+\Vert \vert S(c^{(n)},Z)\vert \Vert _{3,p})^{a}\Big)
\\
&\leq &C\left\Vert f\right\Vert _{\infty }\Big(\P(\sigma _{S_{\infty
}(c^{(n)},Z)}\leq \eta )+\frac{\sqrt{\delta }}{\eta ^{p}}C\Big)
\end{eqnarray*}%
the last inequality being a consequence of (\ref{nss8}) and (\ref{CTV1a}).
And a similar inequality holds for $G.$ So%
\begin{eqnarray*}
\left\vert {\mathbb{E}}(f(S_\infty(c^{(n)},Z)))-{\mathbb{E}}(f(G))\right\vert
&\leq&\left\vert {\mathbb{E}}(f_{\delta }(S_\infty(c^{(n)},Z)))-{\mathbb{E}}(f_{\delta
}(G))\right\vert \\
&&+C\left\Vert f\right\Vert _{\infty }\Big(\P(\sigma _{S_{\infty
}(c^{(n)},Z)}\leq \eta )+\P(\sigma _{G}\leq \eta )+\frac{\sqrt{\delta }}{\eta
^{p}}C\Big)\\
&\leq&\left\vert\E(
f_{\delta }(S_N(c^{(n)},Z)))-{\mathbb{E}}(f_{\delta
}(G))\right\vert\\
&&+\left\vert {\mathbb{E}}(f_{\delta }(S_\infty(c^{(n)},Z)))-\E(
f_{\delta }(S_N(c^{(n)},Z)))\right\vert\\
&&+C\left\Vert f\right\Vert _{\infty }\Big(\P(\sigma _{S_{\infty
}(c^{(n)},Z)}\leq \eta )+\frac{\sqrt{\delta }}{\eta
^{p}}C\Big)\\
&=:& A^\delta_N(n)+B^\delta_N(n)+C\left\Vert f\right\Vert _{\infty }\Big(\P(\sigma _{S_{\infty
}(c^{(n)},Z)}\leq \eta )+\frac{\sqrt{\delta }}{\eta
^{p}}C\Big)
\end{eqnarray*}%
in which we have used the fact that $\sigma _{G}=1>\eta $, so that $\P(\sigma _{G}\leq \eta )=0$. Now, by using Theorem \ref{TV copy(2)} and by recalling that $\|f_\delta\|_\infty\leq \|f\|_\infty$,
$$
\limsup_{n\to\infty} \sup_{\|f\|_\infty\leq 1}A^\delta_N(n)=0,
$$
for every fixed $N>0$. Moreover, since $\|f'_\delta\|_\infty\leq \|f\|_\infty/\delta$, we can write
$$
B^\delta_N(n)\leq
\frac{\|f\|_\infty}{\delta}
\E\big(|S^N(c^{(n)},Z)|\big)
\leq \frac{\|f\|_\infty}{\delta}
\varepsilon_N(n)^{1/2}
$$
so that
$$
\limsup_{N\to\infty}\limsup_{n\to\infty} \sup_{\|f\|_\infty\leq 1}B^\delta_N(n)=0
$$
because of (\ref{CTV2}), for every $\delta>0$.  By using (\ref{CTV5}) we finally get
$$
\limsup_{n\to\infty}
\sup_{\|f\|_\infty\leq 1}\left\vert {\mathbb{E}}(f(S_\infty(c^{(n)},Z)))-{\mathbb{E}}(f(G))\right\vert
\leq
C\frac{\sqrt{\delta }}{\eta^{p}}.
$$
Since $\delta >0$ is arbitrary the proof is complete. $\square $

\appendix

\section{Burkholder inequality for Hilbert valued discrete time martingales}
\label{burk}

We consider a Hilbert space $\mathcal{U}$ and we denote $\left\vert
\cdot\right\vert _{\mathcal{U}}$ and $\<\cdot,\cdot\>_\mathcal{U}$
respectively the norm and the inner product on $\mathcal{U}$. Recall $L^p_{%
\mathcal{U}}$ and ${\mathbb{D}}^{q,p}_{\mathcal{U}}$ defined at the
beginning of Section \ref{sect-MallS}.

We consider a martingale $M_{n}\in \mathcal{U},n\in \N$ and we recall
Burkholder's inequality in this framework: for each $p\geq 2$ there exists a
universal constant $b_{p}\geq 1$ such that%
\begin{equation}
\left\Vert M_{n}\right\Vert _{\mathcal{U},p}\leq b_{p}\Big({\mathbb{E}}\Big(%
\Big(\sum_{k=1}^{n}\left\vert \Delta _{k}\right\vert _{\mathcal{U}}^{2}\Big)%
^{p/2}\Big)\Big)^{1/p},\qquad \Delta _{k}=M_{k}-M_{k-1}.  \label{a1}
\end{equation}%
As an immediate consequence
\begin{equation}
\left\Vert M_{n}\right\Vert _{\mathcal{U},p}\leq b_{p}\Big(%
\sum_{k=1}^{n}\left\Vert \Delta _{k}\right\Vert _{\mathcal{U},p}^{2}\Big)%
^{1/2}.  \label{a2}
\end{equation}%
Indeed, by using (\ref{a1}),
\begin{equation*}
\left\Vert M_{n}\right\Vert _{\mathcal{U},p}^{2}\leq b_{p}^{2}\big\|%
\sum_{k=1}^{n}|\Delta _{k}|_{\mathcal{U}}^{2}\big\|_{p/2}\leq
b_{p}^{2}\sum_{k=1}^{n}\Vert |\Delta _{k}|_{\mathcal{U}}^{2}\Vert
_{p/2}=b_{p}^{2}\sum_{k=1}^{n}\Vert \Delta _{k}\Vert _{\mathcal{U},p}^{2}.
\end{equation*}

\medskip We give now estimates which are used in order to upper bound the
Sobolev norms of $LS(c,Z)$. Recall the definition of the space ${\mathbb{D}}%
_{\mathcal{U}}^{q,p}$ given in Section \ref{sect-MallS} and we set ${\mathbb{%
D}}_{\mathcal{U}}^{\infty }=\cap _{p\geq 1}\cap _{q\geq 0}{\mathbb{D}}_{%
\mathcal{U}}^{q,p}$.

\begin{proposition}
Suppose that $Z_{k}\in \mathcal{L}(z_k,r,\varepsilon )$, $k\in {\mathbb{N}}$.
Let $B_{k},\Lambda _{k}\in \mathcal{U}$ be random variables such that $%
B_{k},\Lambda _{k}\in {\mathbb{D}}_{\mathcal{U}}^{\infty }$ for every $k$
and $B_{k}$ is $\sigma (Z_{1},...,Z_{k})$ measurable. Consider the process
\begin{equation}
Y_{m}=\sum_{k=1}^{m-1}B_{k}LZ_{k+1}+\Lambda _{m}.  \label{aa3}
\end{equation}%
Then for every $q\in {\mathbb{N}}$ and $p\geq 2$ there exists a universal
constant $C\geq 1$ such that
\begin{equation}
\max_{m\leq n}\left\Vert Y_{m}\right\Vert _{\mathcal{U},q,p}\leq \frac{C}{%
r^{q+1}}\times C_{q,p}(B,\Lambda )  \label{b2}
\end{equation}%
with
\begin{equation}
C_{q,p}(B,\Lambda )=\Big(\sum_{k=1}^{n}\left\Vert B_{k}\right\Vert _{%
\mathcal{U},q,p}^{2}\Big)^{1/2}+\max_{m\leq n}\left\Vert \Lambda
_{m}\right\Vert _{\mathcal{U},q,p}.  \label{b1}
\end{equation}
\end{proposition}

\textbf{Proof.} We will use the following facts, proved in Lemma 3.2 in \cite%
{[BC-CLT]}: ${\mathbb{E}}(LZ_{k})=0$ and there exists a universal constant $%
C $ such that%
\begin{equation}
\left\Vert LZ_{k}\right\Vert _{q,p}\leq \frac{C}{r^{q+1}}.  \label{aa4}
\end{equation}

\textbf{Step 1}. Let $q=0,$ so that $\left\Vert Y_{m}\right\Vert _{\mathcal{U%
},q,p}=\left\Vert Y_{m}\right\Vert _{\mathcal{U},p}.$ We have to check that%
\begin{equation}
\max_{m\leq n}\left\Vert Y_{m}\right\Vert _{\mathcal{U},p}\leq \frac{C}{r}%
\times C_{0,p}(B,\Lambda).  \label{a3'}
\end{equation}%
Since $B_{k}$ is $\sigma (Z_{1},...,Z_{k})$ measurable and ${\mathbb{E}}%
(LZ_{k+1})=0,$ it follows that $M_{m}=\sum_{k=1}^{m-1}B_{k}LZ_{k+1}$ is a
martingale. By (\ref{a2})%
\begin{equation*}
\left\Vert M_{m}\right\Vert _{\mathcal{U},p}\leq b_{p}\Big(%
\sum_{k=1}^{m}\left\Vert LZ_{k+1}B_{k}\right\Vert _{\mathcal{U},p}^{2}\Big)%
^{1/2}.
\end{equation*}%
Since $LZ_{k+1}$ and $B_{k}$ are independent,
\begin{equation*}
\left\Vert LZ_{k+1}B_{k}\right\Vert _{\mathcal{U},p}^{2}=\left\Vert
LZ_{k+1}\right\Vert _{p}^{2}\left\Vert B_{k}\right\Vert _{\mathcal{U}%
,p}^{2}\leq \frac{C}{r^{2}}\left\Vert B_{k}\right\Vert _{\mathcal{U},p}^{2}.
\end{equation*}%
From $Y_{m}=M_{m}+\Lambda_{m}$, we conclude that%
\begin{equation*}
\Vert Y_{m}\Vert _{\mathcal{U},p}\leq \Vert M_{m}\Vert _{\mathcal{U}%
,p}+\Vert \Lambda_{m}\Vert _{\mathcal{U},p}\leq \frac{C}{r}\Big(\Big(%
\sum_{k=1}^{m}\Vert B_{k}\Vert _{\mathcal{U},p}^{2}\Big)^{1/2}+\Vert \Lambda
_{m}\Vert _{\mathcal{U},p}\Big)
\end{equation*}%
and the statement holds for $q=0$.

\textbf{Step 2}. We estimate the derivatives of $Y_{m}$. We have%
\begin{equation*}
\overline{Y}_{m}:=DY_{m}=\sum_{k=1}^{m-1}\overline{B}_{k}LZ_{k+1}+\overline{%
\Lambda}_{m}.
\end{equation*}%
with $\overline{B}_{k}=DB_{k}$ and $\overline{\Lambda}_{m}=%
\sum_{k=1}^{m-1}DLZ_{k+1}B_{k}+D\Lambda_{m}.$ Notice that $\overline{Y}_{m}$%
, $\overline{B}_{k}$ and $\overline{\Lambda}_{m}$ take values in $\mathcal{H}%
(\mathcal{U})$ (defined in (\ref{H-U})). So, by applying the step above, we
get
\begin{equation*}
\max_{m\leq n}\Vert DY_{m}\Vert _{\mathcal{H}(\mathcal{U}),p}\leq \frac{C}{r}%
C_{0,p}(\overline{B},\overline{\Lambda}),
\end{equation*}%
where
\begin{equation*}
C_{0,p}(\overline{B},\overline{\Lambda})=\max_{m\leq n}\Big(\Big(%
\sum_{k=1}^{m}\left\Vert \overline{B}_{k}\right\Vert _{\mathcal{H}(\mathcal{U%
}),p}^{2}\Big)^{1/2}+\left\Vert \overline{\Lambda}_{m}\right\Vert _{\mathcal{%
H}(\mathcal{U}),p}\Big).
\end{equation*}%
If we prove that
\begin{equation}
C_{0,p}(\overline{B},\overline{\Lambda})\leq \frac{C}{r}\times
C_{1,p}(B,\Lambda)  \label{b3}
\end{equation}%
(hereafter, $C>0$ denotes a constant that may vary) and recalling that $%
C_{0,p}(B,\Lambda)\leq C_{1,p}(B,\Lambda)$, then we obtain
\begin{equation*}
\max_{m\leq n}\Vert Y_{m}\Vert _{\mathcal{U},1,p}\leq \frac{C}{r^{2}}%
C_{1,p}(B,\Lambda).
\end{equation*}%
And by iteration, we get (\ref{b2}). So, let us prove (\ref{b3}).

We have $\Vert \overline{B}_{k}\Vert _{\mathcal{H}(\mathcal{U}),p}=\Vert
DB_{k}\Vert _{\mathcal{H}(\mathcal{U}),p}\leq \left\Vert B_{k}\right\Vert _{%
\mathcal{U},1,p}$. We analyze now $\overline{\Lambda}_{m}.$ First, $\Vert
D\Lambda_{k}\Vert _{\mathcal{H}(\mathcal{U}),p}\leq \Vert \Lambda_{m}\Vert _{%
\mathcal{U},1,p}$. Let $I_{m}:=\sum_{k=1}^{m-1}DLZ_{k+1}B_{k}$. Since $%
D_{p}LZ_{k+1}=0$ if $p\neq k+1$ we obtain%
\begin{equation*}
\left\vert I_{m}\right\vert _{\mathcal{H}(\mathcal{U})}^{2}=\sum_{k=1}^{m-1}%
\left\vert D_{k+1}LZ_{k+1}\right\vert ^{2}\left\vert B_{k}\right\vert _{%
\mathcal{U}}^{2}.
\end{equation*}%
Recalling that $D_{k+1}LZ_{k+1}$ and $B_{k}$ are independent and that $%
\left\Vert D_{k+1}LZ_{k+1}\right\Vert _{p}^{2}\leq Cr^{-2}$, we can write
\begin{align*}
\Vert I_{m}\Vert _{\mathcal{H}(\mathcal{U}),p}& =\Vert |I_{m}|_{\mathcal{H}(%
\mathcal{U})}^{2}\Vert _{p/2}^{1/2}\leq \Big(\sum_{k=1}^{m-1}\Vert
\left\vert D_{k+1}LZ_{k+1}\right\vert ^{2}\left\vert B_{k}\right\vert _{%
\mathcal{U}}^{2}\Vert _{p/2}\Big)^{1/2} \\
& =\Big(\sum_{k=1}^{m-1}\Vert D_{k+1}LZ_{k+1}\Vert _{p}^{2}\Vert B_{k}\Vert
_{\mathcal{U},p}^{2}\Big)^{1/2}\leq \frac{C}{r}\times \Big(%
\sum_{k=1}^{m-1}\Vert B_{k}\Vert _{\mathcal{U},p}^{2}\Big)^{1/2}.
\end{align*}%
By inserting all these estimates, we get (\ref{b3}). $\square $

\section{The $L^p$ estimates in Lemma \protect\ref{crucial-est}}

\label{app-B}

\subsection{Contractions and cumulants}

\label{cumulants}

We briefly recall some well known facts concerning contractions of kernels and
cumulants, and we give some easy consequences which are used in our paper.
The results in this section involve $\left\vert c\right\vert _{m}$ (see (\ref%
{Law2})), $\delta _{m}(c)$ (see (\ref{Law3''})) and $\kappa _{4,m}(c)$ (see (%
\ref{Law6''})). We denote by $c_{(m)}(\alpha )=1_{\{\left\vert \alpha
\right\vert =m\}}c(\alpha ),$ so $c_{(m)}\in \mathcal{H}^{\otimes m}$
represents the restriction of $c$ to $\Gamma _{m}.$ Then for $m,n\in {%
\mathbb{N}}$ and $0\leq r\leq m\wedge n$ we define the contraction $%
c_{(m)}\otimes _{r}c_{(n)}\in \mathcal{H}^{\otimes (m+n-2r)}$ as follows
\begin{equation}
c_{(m)}\otimes _{r}c_{(n)}(\alpha ,\beta )=\sum_{\left\vert \gamma
\right\vert =r}c_{(m)}((\alpha ,\gamma ))c_{(n)}((\beta ,\gamma
))=\sum_{\left\vert \gamma \right\vert =r}c((\alpha ,\gamma ))c((\beta
,\gamma ))  \label{C1}
\end{equation}%
where $\alpha =(\alpha _{1},...,\alpha _{m-r}),\beta =(\beta _{1},...,\beta
_{n-r})$. Since for $m\neq n,$ $c_{(m)}\otimes _{r}c_{(n)}$ is not
symmetric, we define $c_{(m)}\widetilde{\otimes }_{r}c_{(n)}$ to be the
symmetrization of $c_{(m)}\otimes _{r}c_{(n)}$: for $\eta \in \Gamma
_{m+n-2r}$,
\begin{equation}
c_{(m)}\widetilde{\otimes }_{r}c_{(n)}(\eta )=\frac{1}{(n+m-2r)!}\sum_{\pi
\in \Pi _{m+n-2r}}c_{(m)}\otimes _{r}c_{(n)}(\mathfrak{p}_{m-r}(\eta _{\pi
}),\mathfrak{q}_{n-r}(\eta _{\pi })),  \label{C1'}
\end{equation}%
in which $\Pi _{m+n-2r}$ denotes the permutations of $\{1,\ldots
,m+n-2r\}$, for $\pi \in \Pi _{m+n-2r}$ then $\eta _{\pi }=(\eta _{\pi
_{1}},\ldots ,\eta _{\pi _{m+n-2r}})$, $\mathfrak{p}_{m-r}$
is the projection on the first $m-r$ coordinates and $\mathfrak{q}_{m-r}$ is
the projection on the last $n-r$ coordinates, with the convention $\mathfrak{p}_{0}=\mathfrak{q}_{0}=\emptyset$.

Finally we recall Remark \ref{iterated}: that if $Z=(Z_{k})_{k\in {\mathbb{N}}}$ with $Z_{k}$
independent standard normal random variables then $\Phi
_{m}(c,Z)=I_{m}(f_{c_{(m})})$ is the multiple stochastic integral with the piecewise constant
kernel $f_{c_{(m)}}(t_{1},...,t_{m})=c(\alpha )$ if $t_{i}\in \lbrack \alpha
_{i},\alpha _{i}+1),i=1,...,m$. So we come back to the Wiener space (the results known in the literature usually concern multiple stochastic integrals) and we summarize all the needed results in the next lemma.

\begin{lemma}
\begin{itemize}
\item One has
\begin{equation}
\kappa_{4,m}(c)
=\sum_{r=1}^{m-1}m!^{2}\left(
\begin{tabular}{l}
$m$ \\
$r$%
\end{tabular}%
\right) ^{2}\{\left\Vert c_{(m)}\otimes _{r}c_{(m)}\right\Vert ^{2}+\left(
\begin{tabular}{l}
$2m-2r$ \\
$m-r$%
\end{tabular}%
\right) \left\Vert c_{(m)}\widetilde{\otimes }_{r}c_{(m)}\right\Vert ^{2}\}.
\label{acc7}
\end{equation}

\item For $0\leq r\leq m\wedge n$, one has
\begin{align}
&\left\Vert c_{(m)}\widetilde{\otimes }_{r}c_{(n)}\right\Vert ^{2}\leq \frac{%
1}{2}(\left\Vert c_{(m)}\otimes _{m-r}c_{(m)}\right\Vert ^{2}+\left\Vert
c_{(n)}\otimes _{n-r}c_{(n)}\right\Vert ^{2})  \label{acc8-1} \\
&\left\Vert c_{(m)}\otimes_{r}c_{(n)}\right\Vert ^{2}\leq
\|c_{(m)}\otimes_{m-r}c_{(m)}\|\,\|c_{(n)}\otimes_{n-r}c_{(n)}\|
\label{acc8}
\end{align}

\item For $0< r< m\wedge n$%
\begin{equation}  \label{acc9}
\begin{array}{l}
\left\Vert c_{(m)}\widetilde{\otimes }_{r}c_{(n)}\right\Vert \leq \max\Big(%
\frac{\sqrt{\kappa_{4,m}(c)}}{m!%
\begin{pmatrix}
m \\
r%
\end{pmatrix}%
}, \frac{\sqrt{\kappa_{4,n}(c)}}{n!%
\begin{pmatrix}
n \\
r%
\end{pmatrix}%
}\Big)\smallskip \\
\left\Vert c_{(m)}\otimes_{r}c_{(n)}\right\Vert \leq \max\Big(\frac{\sqrt{%
\kappa_{4,m}(c)}}{m!%
\begin{pmatrix}
m \\
r%
\end{pmatrix}%
}, \frac{\sqrt{\kappa_{4,n}(c)}}{n!%
\begin{pmatrix}
n \\
r%
\end{pmatrix}%
}\Big)%
\end{array}%
\end{equation}

\item For $1\leq m\leq n-1$%
\begin{equation}  \label{acc9-1}
\begin{array}{l}
\left\Vert c_{(m)}\otimes_{m}c_{(n)}\right\Vert \leq
\|c_{(m)}\|^2\,\|c_{(n)}\otimes_{n-m}c_{(n)}\| \leq \|c_{(m)}\|^2\Big(%
\frac{\sqrt{\kappa_{4,n}(c)}}{n!\mbox{\scriptsize$\begin{pmatrix}n\\n-m%
\end{pmatrix}$}}\Big)^{1/2}.%
\end{array}%
\end{equation}

\item The following estimate for the influence factor $\delta _{m}(c)$ holds:%
\begin{equation}
\delta _{m}(c)\leq \frac{1}{\left\vert c\right\vert _{m}}\left\vert
c_{(m)}\otimes _{m-1}c_{(m)}\right\vert _{2m-2}\leq \frac{\sqrt{\kappa
_{4,m}(c)}}{m!m\left\vert c\right\vert _{m}}  \label{acc10}
\end{equation}
\end{itemize}
\end{lemma}

\textbf{Proof}. We first recall that $\kappa_{4,m}(c)=\kappa _{4}(\Phi _{m}(c,Z))$ with $Z_k$, $k\in\N$, standard normal. So, the identity (\ref{acc7}) is proved in \cite{[PN]} for iterated integrals and remains true for stochastic series because $\|
f_{c_{(m)}}\| _{L({\mathbb{R}}_{+}^{m})}$ $=\| c_{(m)}\|$ and $%
f_{c_{(m)}}\otimes f_{c_{(n)}}=f_{c_{(m)}\otimes c_{(n)}}$. (\ref{acc8-1})
is straightforward (but see also formula (13) in \cite{[NO]}) and (\ref{acc8}%
) appears in \cite{[NO]} and \cite{[NPRev]}. (\ref{acc9}) is an immediate
consequence of (\ref{acc8-1})-(\ref{acc8}) and (\ref{acc7}). Concerning (\ref%
{acc9-1}), straightforward computations give
\begin{equation*}
\|c_{(m)}\otimes_mc_{(n)}\|^2
=\sum_{\rho,\bar\rho\in\Gamma_m}c_{(m)}(\rho)c_{(m)}(\bar\rho)
c_{(n)}\otimes_{n-m} c_{(n)}(\rho,\bar\rho).
\end{equation*}
By using the Cauchy-Schwarz inequality we get $\|c_{(m)}\otimes_mc_{(n)}%
\|^2\leq \|c_{(m)}\|^2\|c_{(n)}\otimes_{n-m} c_{(n)}\|$. Last inequality in (%
\ref{acc9-1}) follows from (\ref{acc9}). Finally, the inequality (\ref{acc10}%
) has been proved in \cite{[NPRein]} $\square $

\subsection{Some $L^{p}$ estimates for series}\label{app-Lp}

\subsubsection{The basic lemma}

We start with the basic definitions of this section.
\begin{assumption}\label{1-4}
\begin{enumerate}
\item We fix $m,n\geq 0$ integers and we consider a coefficient
\begin{equation*}
a\,:\,\Gamma_m\times \Gamma_n\to {\mathbb{R}},\quad (\alpha,\beta)\mapsto
a(\alpha,\beta)
\end{equation*}
that satisfies:

\smallskip

\ \ $\bullet$ for $n+m\geq 1$, as a function of $\gamma=(\alpha,\beta)\in%
\Gamma_{m+n}$, $a$ is null on the diagonals;

\smallskip

\ \ $\bullet$ for $m,n\geq 2$, $\Gamma_m\ni\alpha \mapsto a(\alpha ,\beta )$ and
$\Gamma_n\ni\beta \mapsto a(\alpha ,\beta )$ are both symmetric (so $a$ is
symmetric in each argument, separately).

We define
\begin{equation}
\left\vert a\right\vert _{m,n,J}=\Big(\sum_{\alpha \in \Gamma
_{m}(J)}\sum_{\beta \in \Gamma _{n}(J)}a^{2}(\alpha ,\beta )\Big)^{1/2}
\label{RA3j}
\end{equation}

\item
Take now $\overline{a}=(a_{j})_{j\in {\mathbb{N}}}$, with $a_{j}:\,\Gamma
_{m}\times \Gamma _{n}\rightarrow {\mathbb{R}}$ which satisfies the
hypotheses in 1. and furthermore with
$$
\mbox{$a_{j}(\gamma )=0$ if $j\in \gamma $.}
$$
In this case, we denote%
\begin{equation}
\left\vert \overline{a}\right\vert _{m,n,J}=\Big(\sum_{j=1}^{\infty
}\sum_{\alpha \in \Gamma _{m}(J)}\sum_{\beta \in \Gamma
_{n}(J)}a_{j}^{2}(\alpha ,\beta )\Big)^{1/2}.  \label{RA3}
\end{equation}
\item
We consider a sequence of centred and independent random variables $%
(Z_{k},Y_{k},\widetilde{\chi }_{k})$, $k\in {\mathbb{N}}$ (with $\widetilde{%
\chi }_{j}=\chi _{j}-{\mathbb{E}}(\chi _{j})$, $Z_k$ and $\chi_k$ having the usual meaning) and we denote%
\begin{equation}
M_{p}(Z,Y)=\sup_{k}\left\Vert Z_{k}\right\Vert _{p}\vee \left\Vert
Y_{k}\right\Vert _{p}.  \label{RA1}
\end{equation}%
Notice that we do not require that $Y_{k}$ is independent of $Z_{k}$ and/or of $%
\widetilde{\chi }_{k}$.
We are interested in the following double series: for a fixed $a$ as in 1.,
\begin{equation}
t_{m,n}(J,a)=\sum_{\alpha \in \Gamma _{m}(J)}\sum_{\beta \in \Gamma
_{n}(J)}Z^{\alpha }Y^{\beta }a(\alpha ,\beta )  \label{RA4'}
\end{equation}%
and for $\overline{a}=(a_{j})_{j\in {\mathbb{N}}}$ as in 2.,
\begin{equation}
T_{m,n}(J,\overline{a})=\sum_{\alpha \in \Gamma _{m}(J)}\sum_{\beta \in
\Gamma _{n}(J)}Z^{\alpha }Y^{\beta }\sum_{j=1}^{\infty }a_{j}(\alpha ,\beta )%
\widetilde{\chi }_{j}.  \label{RA4}
\end{equation}
\end{enumerate}
\end{assumption}
\begin{lemma}
\label{lemma-estTt} Under Assumption \ref{1-4} one has
\begin{eqnarray}
\left\Vert t_{m,n}(J,a)\right\Vert _{p} &\leq & \big((m+n)!\big)^{1/2}\,(
\sqrt 2\,b_{p}M_{p}(Z,Y))^{(n+m)}|a|_{m,n,J},  \label{RA9'} \\
\left\Vert T_{m,n}(J,\overline{a})\right\Vert _{p} &\leq & \Big(\frac{%
8b_p^2(4^{m+n}-1)(m+n)!}{3}\Big)^{1/2}\big(\sqrt 2\,b_pM_p(Z,Y)\big)^{m+n}
\left\vert \overline{a}\right\vert _{m,n,J}.  \label{RA9}
\end{eqnarray}
\end{lemma}

\textbf{Proof.} We recall that $a$ and $a_{j}$ are all null on the
diagonals. So, the sums in (\ref{RA4'}) and (\ref{RA4}) are really done on
the multi-indexes $\alpha$ and $\beta $ that do not have common components.
So, we consider such kind of indexes.

\smallskip

\textbf{Step 1}. We denote
\begin{eqnarray*}
\Lambda _{m,n}^{\prime }(J^{\prime }) &=&(\Gamma _{m}(J^{\prime })\setminus
\Gamma _{m}(J^{\prime }-1))\times \Gamma _{n}(J^{\prime }-1) \\
\Lambda _{m,n}^{\prime \prime }(J^{\prime }) &=&\Gamma _{m}(J^{\prime
})\times (\Gamma _{n}(J^{\prime })\setminus \Gamma _{n}(J^{\prime }-1)).
\end{eqnarray*}%
So $(\alpha ,\beta )\in $ $\Lambda _{m,n}^{\prime }(J^{\prime })$ if $%
\max_{j=1,...m}\alpha _{j}=J^{\prime }$ and $\max_{j=1,...n}\beta _{j}\leq
J^{\prime }-1.$ And the definition of $\Lambda _{m,n}^{\prime \prime
}(J^{\prime })$ is similar, with $\alpha $ replaced by $\beta .$ Finally we
put
\begin{equation*}
\Lambda _{m,n}(J^{\prime })=\Lambda _{m,n}^{\prime }(J^{\prime })\cup
\Lambda _{m,n}^{\prime \prime }(J^{\prime }).
\end{equation*}
So, we have
\begin{equation*}
t_{m,n}(J,a)=\sum_{\alpha \in \Gamma _{m}(J)}\sum_{\beta \in \Gamma
_{n}(J)}Z^{\alpha }Y^{\beta }a(\alpha ,\beta )=\sum_{J^{\prime
}=1}^{J}\sum_{(\alpha ,\beta )\in \Lambda _{m,n}(J^{\prime })}Z^{\alpha
}Y^{\beta }a(\alpha ,\beta ).
\end{equation*}%
In order to prove (\ref{RA9'}), the first step is to establish a recurrence
formula. We define%
\begin{equation}  \label{Qa}
(Q_{J}^{\prime }a)(\alpha ,\beta )=a((\alpha ,J),\beta )\qquad \mbox{and}%
\qquad (Q_{J}^{\prime \prime }a)(\alpha ,\beta )=a(\alpha ,(\beta ,J))
\end{equation}%
and we write
\begin{eqnarray*}
t_{m,n}(J,a) &=&\sum_{J^{\prime }=1}^{J}\sum_{(\alpha ,\beta )\in \Lambda
_{m,n}^{\prime }(J^{\prime })}Z^{\alpha }Y^{\beta }a(\alpha ,\beta
)+\sum_{J^{\prime }=1}^{J}\sum_{(\alpha ,\beta )\in \Lambda _{m,n}^{\prime
\prime }(J^{\prime })}Z^{\alpha }Y^{\beta }a(\alpha ,\beta ).
\end{eqnarray*}
But $(\alpha ,\beta )\in \Lambda _{m,n}^{\prime }(J^{\prime })$ if and only
if $\beta\in\Gamma_n(J^{\prime }-1)$ and $\alpha$ contains $J^{\prime }$,
the remaining entries forming a multi-index in $\Gamma_{m-1}(J^{\prime }-1)$%
. And similarly for $(\alpha ,\beta )\in \Lambda_{m,n}^{\prime \prime
}(J^{\prime })$, changing the role to $\alpha$ and $\beta$. So, by using the
symmetry of $\alpha\mapsto a(\alpha,\beta)$ and $\beta\mapsto
a(\alpha,\beta) $, we can write
\begin{eqnarray*}
t_{m,n}(J,a)&=&m\sum_{J^{\prime }=1}^{J}Z_{J^{\prime }}\sum_{\alpha \in
\Gamma _{m-1}(J^{\prime }-1)}\sum_{\beta \in \Gamma _{n}(J^{\prime
}-1)}Z^{\alpha }Y^{\beta }a((\alpha ,J^{\prime }),\beta ) \\
&&+n\sum_{J^{\prime }=1}^{J}Y_{J^{\prime }}\sum_{\alpha \in \Gamma
_{m}(J^{\prime }-1)}\sum_{\beta \in \Gamma _{n-1}(J^{\prime }-1)}Z^{\alpha
}Y^{\beta }a(\alpha ,(\beta ,J^{\prime })) \\
&=&\sum_{J^{\prime }=1}^{J}(Z_{J^{\prime }}mt_{m-1,n}(J^{\prime
}-1,Q_{J^{\prime }}^{\prime }a)+Y_{J^{\prime }}nt_{m,n-1}(J^{\prime
}-1,Q_{J^{\prime }}^{\prime \prime }a)).
\end{eqnarray*}

Let $\mathcal{G}_{n}=\sigma \{(Z_{k},Y_{k},\widetilde{\chi }_{k}):k\leq n\}.$%
\ Notice that $t_{m-1,n}(J^{\prime }-1,Q_{J^{\prime }}^{\prime }a)$ and $%
t_{m,n-1}(J^{\prime }-1,Q_{J^{\prime }}^{\prime \prime }a)$ are $\mathcal{G}%
_{J^{\prime }-1}$ measurable so the above sums are martingales with respect
to the filtration $\mathcal{G}_{n}$, $n\in {\mathbb{N}}$, and we may use
Burkholder's inequality. Using the above recurrence formula and the
recurrence hypotheses (\ref{RA9'}) we obtain
\begin{eqnarray*}
\Vert t_{m,n}(J,a)\Vert _{p}^{2} &\leq &(\sqrt{2}b_{p}M_{p}(Z,Y))^{2}\Big(%
\sum_{J^{\prime }=1}^{J}m^{2}\Vert t_{m-1,n}(J^{\prime },Q_{J^{\prime
}}^{\prime }a)\Vert _{p}^{2}+\sum_{J^{\prime }=1}^{J-1}n^{2}\Vert
t_{m,n-1}(J^{\prime },Q_{J^{\prime }}^{\prime \prime }a)\Vert _{p}^{2}\Big)
\\
&\leq &(m+n-1)!(\sqrt{2}b_{p}M_{p}(Z,Y))\Big)^{2(n+m)}(\sum_{J^{\prime
}=1}^{J}\sum_{\alpha \in \Gamma _{m-1}(J^{\prime }-1)}\sum_{\beta \in \Gamma
_{n}(J^{\prime }-1)}m^{2}a^{2}((\alpha ,J^{\prime }),\beta ) \\
&&+\sum_{J^{\prime }=1}^{J}\sum_{\alpha \in \Gamma _{m}(J^{\prime
}-1)}\sum_{\beta \in \Gamma _{n-1}(J^{\prime }-1)}n^{2}a^{2}((\alpha
,J^{\prime }),\beta )) \\
&\leq &(m+n-1)!(\sqrt{2}b_{p}M_{p}(Z,Y))\Big)^{2(n+m)}\times
(m+n)\sum_{\alpha \in \Gamma _{m}(J)}\sum_{\beta \in \Gamma
_{n}(J)}a^{2}(\alpha ,\beta ),
\end{eqnarray*}%
so (\ref{RA9'}) is proved.

\medskip

\textbf{Step 2}. We prove (\ref{RA9}). We write%
\begin{eqnarray*}
T_{m,n}(J,\overline{a}) &=&\sum_{J^{\prime }=1}^{J}\sum_{(\alpha ,\beta )\in
\Lambda _{m,n}(J^{\prime })}Z^{\alpha }Y^{\beta }\sum_{j=J^{\prime
}+1}^{\infty }a_{j}(\alpha ,\beta )\widetilde{\chi }_{j} +\sum_{J^{\prime
}=1}^{J}\sum_{(\alpha ,\beta )\in \Lambda _{m,n}(J^{\prime })}Z^{\alpha
}Y^{\beta }\sum_{j=1}^{J^{\prime }-1}a_{j}(\alpha ,\beta )\widetilde{\chi }%
_{j} \\
&=&A_{m,n}(J,\overline{a})+B_{m,n}(J,\overline{a}).
\end{eqnarray*}%
Notice that the term $j=J^{\prime }$ does not appear: for $%
(\alpha,\beta)\in\Lambda_{m,n}(J^{\prime })$ then $J^{\prime
}\in(\alpha,\beta)$, so $a_{J^{\prime }}(\alpha,\beta)=0$ by our assumption.
Now we write%
\begin{equation*}
A_{m,n}(J,\overline{a})=\sum_{j=2}^{\infty }\widetilde{\chi }%
_{j}\sum_{J^{\prime }=1}^{J\wedge (j-1)}\sum_{(\alpha ,\beta )\in \Lambda
_{m,n}(J^{\prime })}Z^{\alpha }Y^{\beta }a_{j}(\alpha ,\beta
)=\sum_{j=1}^{\infty }\widetilde{\chi }_{j}t_{m,n}(J\wedge (j-1),a_j).
\end{equation*}%
We use Burkholder's inequality and (\ref{RA9'}) in order to obtain%
\begin{eqnarray}
\left\Vert A_{m,n}(J,\overline{a})\right\Vert _{p} &\leq
&2b_{p}(\sum_{j=2}^{\infty }\left\Vert t_{m,n}(J\wedge (j-1),a_j)\right\Vert
_{p}^{2})^{1/2}  \label{RA8} \\
&\leq &2b_{p}(\sqrt{2}b_{p}M_{p}(Z,Y))^{(n+m)}((m+n)!)^{1/2}\times  \notag \\
&&\times \Big(\sum_{j=1}^{\infty }\sum_{\alpha \in \Gamma _{m}(J\wedge
(j-1))}\sum_{\beta \in \Gamma _{n}(J\wedge (j-1))}a_{j}^{2}(\alpha ,\beta )%
\Big)^{1/2}.  \notag
\end{eqnarray}

\textbf{Step 3}. We estimate now $\left\Vert B_{m,n}(J,\overline{a}%
)\right\Vert _{p}.$ We write%
\begin{eqnarray*}
B_{n,m}(J,\overline{a}) &=&\sum_{J^{\prime }=1}^{J}Z_{J^{\prime
}}m\sum_{\alpha \in \Gamma _{m-1}(J^{\prime }-1)}\sum_{\beta \in \Gamma
_{n}(J^{\prime }-1)}Z^{\alpha }Y^{\beta }\sum_{j=1}^{\infty
}1_{\{j<J^{\prime }\}}(Q_{J^{\prime }}^{\prime }a_{j})(\alpha ,\beta )%
\widetilde{\chi }_{j} \\
&&+\sum_{J^{\prime }=1}^{J}Y_{J^{\prime }}n\sum_{\alpha \in \Gamma
_{m}(J^{\prime }-1)}\sum_{\beta \in \Gamma _{n-1}(J^{\prime }-1)}Z^{\alpha
}Y^{\beta }\sum_{j=1}^{\infty }1_{\{j<J^{\prime }\}}(Q_{J^{\prime }}^{\prime
\prime }a_{j})(\alpha ,\beta )\widetilde{\chi }_{j} \\
&=&\sum_{J^{\prime }=1}^{J}(Z_{J^{\prime }}mT_{m-1,n}(J^{\prime
}-1,q_{J^{\prime }}^{\prime }\overline{a})+Y_{J^{\prime
}}nT_{m,n-1}(J^{\prime }-1,q_{J^{\prime }}^{\prime \prime }\overline{a})).
\end{eqnarray*}%
with
\begin{align*}
&(q_{J^{\prime }}^{\prime }\overline{a})_{j}(\alpha ,\beta
)=1_{\{j<J^{\prime }\}}(Q_{J^{\prime }}^{\prime }a_{j})(\alpha ,\beta
)=1_{\{j<J^{\prime }\}}a_{j}((\alpha ,J^{\prime }),\beta ) \\
&(q_{J^{\prime }}^{\prime \prime }\overline{a})_{j}(\alpha ,\beta
)=1_{\{j<J^{\prime }\}}(Q_{J^{\prime }}^{\prime \prime }a_{j})(\alpha ,\beta
)=1_{\{j<J^{\prime }\}}a_{j}(\alpha ,(\beta ,J^{\prime })),
\end{align*}
$Q_{J}^{\prime }a_j$ and $Q_{J}^{\prime\prime }a_j$ being defined in (\ref%
{Qa}). By using Burkholder's inequality,
\begin{equation*}
\left\Vert B_{m,n}(J,\overline{a})\right\Vert _{p}^{2} \leq
2b_{p}^{2}M_{p}^{2}(Z,Y)\sum_{J^{\prime }=1}^{J}(m^2\left\Vert
T_{m-1,n}(J^{\prime }-1,q_{J^{\prime }}^{\prime }\overline{a})\right\Vert
_{p}^{2}+n^2\left\Vert T_{m,n-1}(J^{\prime }-1,q_{J^{\prime }}^{\prime
\prime }\overline{a})\right\Vert _{p}^{2}).
\end{equation*}
So finally%
\begin{align*}
&\left\Vert T_{m,n}(J,\overline{a})\right\Vert _{p}^{2} \leq 2\left\Vert
A_{m,n}(J,\overline{a})\right\Vert _{p}^{2}+2\left\Vert B_{m,n}(J,\overline{a%
})\right\Vert _{p}^{2} \\
&\quad\leq 2\left\Vert A_{n,m}(J,\overline{a})\right\Vert _{p}^{2}
+2b_{p}^{2}M_{p}^{2}(Z,Y)\sum_{J^{\prime }=1}^{J}(m^2\left\Vert
T_{m-1,n}(J^{\prime }-1,q_{J^{\prime }}^{\prime }\overline{a})\right\Vert
_{p}^{2}+n^2\left\Vert T_{m,n-1}(J^{\prime }-1,q_{J^{\prime }}^{\prime
\prime }\overline{a})\right\Vert _{p}^{2})
\end{align*}%
Using the recurrence hypothesis and (\ref{RA9'}) we conclude the proof of (%
\ref{RA9}). $\square $

\subsubsection{The ``product formula''}

We have to deal with $|S(f,Z)|^{2}$ with
\begin{equation*}
S(f,Z)=\sum_{m\geq 0}\sum_{\alpha \in \Gamma _{m}}f_{(m)}(\alpha )Z^{\alpha
},
\end{equation*}%
where $f=\{f_{(m)}\}_{m}$ is a symmetric sequence of coefficients in $%
\mathcal{H}^{\otimes m}$, $m\in {\mathbb{N}}$, which are null on all
diagonals. Note that the case $m=0$ is allowed, by setting $f_{(0)}\in {%
\mathbb{R}}$, $\Gamma _{0}=\{\emptyset \}$ and $Z^{\emptyset }=1$. We then
study $|S(f,Z)|^{2}$. To this purpose, we recall that $\Pi _{n}$ denotes the
set of all permutations of $(1,\ldots ,n)$; for $\eta \in \Gamma _{n}$ and $%
\pi \in \Pi _{n}$, we set $\eta _{\pi }=(\eta _{\pi _{1}},\ldots ,\eta _{\pi
_{n}})$.

\begin{lemma}
\label{lemma-app1} We have
\begin{equation}
|S(f,Z)|^{2}=\sum_{m\geq 0}\sum_{n\geq 0}\sum_{\gamma \in \Gamma
_{m}}\sum_{\eta \in \Gamma _{n}}(Z^{\gamma })^{2}Z^{\eta }A_{n,m}[f](\eta
,\gamma )  \label{R1}
\end{equation}%
where, for $n,m\geq 0$, $\eta \in \Gamma _{n}$ and $\gamma \in \Gamma _{m}$,
\begin{equation}
A_{n,m}[f](\eta ,\gamma )=\frac{m!}{n!}\sum_{a=0}^{n}%
\begin{pmatrix}
a+m \\
m%
\end{pmatrix}%
\begin{pmatrix}
n-a+m \\
m%
\end{pmatrix}%
\sum_{\pi \in \Pi _{n}}f_{(a+m)}(\mathfrak{p}_{a}(\eta _{\pi }),\gamma
)f_{(n-a+m)}(\mathfrak{q}_{n-a}(\eta _{\pi }),\gamma )  \label{A}
\end{equation}%
in which, for $\eta \in \Gamma _{n}$ and $a=0,1,\ldots ,n$,
\begin{equation}
\mathfrak{p}_{a}(\eta )=\left\{
\begin{array}{ll}
\emptyset & \mbox{ for }a=0\smallskip \\
(\eta _{1},\ldots ,\eta _{a}) & \mbox{ for }1\leq a\leq n%
\end{array}%
\right. \quad \mbox{and}\quad \mathfrak{q}_{n-a}(\eta )=\left\{
\begin{array}{ll}
(\eta _{a+1},\ldots ,\eta _{n}) & \mbox{ for }0\leq a\leq n-1\smallskip \\
\emptyset & \mbox{ for }a=n%
\end{array}%
\right.  \label{ppqq}
\end{equation}%
Note that $A_{m,n}[f](\eta ,\gamma )=0$ if $\eta $ and $\gamma $ have common
components and the maps $\eta \mapsto A_{m,n}[f](\eta ,\gamma )$ and $\gamma
\mapsto A_{m,n}[f](\eta ,\gamma )$ are both symmetric.
\end{lemma}

\begin{remark}\label{no-product}
People working in Wiener chaos use the product formula for
multiple stochastic integrals in order to compute $|S(f,Z)|^{2}.$ But this
is not possible here. Suppose that we want to do it in the case where the $Z_{k}$'s
are standard normal - so $S(f,Z)$ is a sum of multiple stochastic integrals (Remark \ref{iterated}).
We stress that the kernels of these integrals are piecewise constant on
the intervals $[k,k+1)$ and if we use the product formula we get multiple
stochastic integrals with kernels which are no more piecewise constant on the
same grid $[k,k+1)$, $k\in \N$ - so we get out from our framework. Put it
otherwise: stochastic series with Gaussian random variables $Z_{k}$, $k\in \N$,
are functionals of the increments of the Brownian motion on $[k,k+1)$, $k\in \N.$
And if we use the product formula for such series we obtain functionals of
the whole Brownian path. Just as an example, if $W_{t}$ is a Brownian motion
and if $Z_{1}=W_{1}$ then $Z_{1}^{2}=2\int_{0}^{1}W_{s}dW_{s}+1.$ Moreover,
in the general case we have no It\^{o} formula which permits to get the above
representation of $Z_{1}^{2}.$ So we have to replace the product formula by (%
\ref{R1}). This leads to some algebraic difficulties but not only. In fact,
the product formula allows to eliminate squares - one comes down to linear
combinations of multiple stochastic integrals and this is very nice because
then one may use the standard Burkholder inequality for them. But here this
does not work and then we have to estimate double series as the ones defined
in (\ref{RA4'}) and (\ref{RA4}).
\end{remark}

\textbf{Proof of Lemma \ref{lemma-app1}.} For $\alpha ,\beta $ multi-indexes, let $\#\alpha \cap \beta $
denote the number of the components which are common to both $\alpha $ and $%
\beta $. For $m,n\geq 0$ and $r=0,\ldots ,m\wedge n$, we set $\Lambda
_{r}^{m,n}=\{(\alpha ,\beta )\in \Gamma _{m}\times \Gamma _{n}\,:\,\#\alpha
\cap \beta =r\}$. Then,
\begin{equation*}
S(f,Z)=\sum_{m\geq 0}\sum_{n\geq 0}\sum_{(\alpha ,\beta )\in \Gamma
_{m}\times \Gamma _{n}}f_{(m)}(\alpha )f_{(n)}(\beta )Z^{\alpha }Z^{\beta
}=\sum_{m\geq 0}\sum_{n\geq 0}\sum_{r=0}^{m\wedge n}\sum_{(\alpha ,\beta
)\in \Lambda _{r}^{m,n}}f_{(m)}(\alpha )f_{(n)}(\beta )Z^{\alpha }Z^{\beta }.
\end{equation*}%
We set $\tilde{\Gamma}_{m}$ as the set of the non-ordered multi-index, that
is the set of all subsets of ${\mathbb{N}}^{m}$, and $\Pi _{m}$ the set of
all permutations of $(1,\ldots ,m)$. For $\alpha =\{\alpha _{1},\ldots
,\alpha _{m}\}\in \tilde{\Gamma}_{m}$ and for $\pi \in \Pi _{m}$ we set $%
\alpha _{\pi }\in \Gamma _{m}$ by $\alpha _{\pi }=(\alpha _{\pi _{1}},\ldots
,\alpha _{\pi _{m}})$. Finally, for $\alpha \in \tilde{\Gamma}_{m}$ and $%
\beta \in \tilde{\Gamma}_{m}$ we set $\alpha \cup \beta $ and $\alpha \cap
\beta $ as the standard reunion and intersection respectively.

Now, $(\alpha ,\beta )\in \Lambda _{r}^{m,n}$ if and only if there exist $%
\gamma \in \tilde{\Gamma}_{r}$, $\bar{\alpha}\in \tilde{\Gamma}_{m-r}$, $%
\bar{\beta}\in \tilde{\Gamma}_{n-r}$, $\pi \in \Pi _{m}$ and $\sigma \in \Pi
_{n}$ such that $\gamma \cap \bar{\alpha}=\emptyset $, $\gamma \cap \bar{%
\beta}=\emptyset $, $\bar{\alpha}\cap \bar{\beta}=\emptyset $ and finally, $%
\alpha =(\bar{\alpha}\cup \gamma )_{\pi }$ and $\beta =(\bar{\beta}\cup
\gamma )_{\sigma }$. Therefore, by using the symmetry property for $f_{(m)}$%
,
\begin{align*}
\sum_{(\alpha ,\beta )\in \Lambda _{r}^{m,n}}f_{(m)}(\alpha )f_{(n)}(\beta
)Z^{\alpha }Z^{\beta }=& \!\!\!\!\sum_{(\pi ,\sigma )\in \Pi _{m}\times \Pi
_{n}}\sum_{%
\mbox{\scriptsize$\begin{array}{c}(\gamma,\bar\alpha,\bar\beta)\in\tilde \Gamma_r\times\tilde\Gamma_{m-r}\times\tilde\Gamma_{n-r}\\
\gamma\cap\bar\alpha=\emptyset, \gamma\cap\bar\beta=\emptyset, \bar\alpha\cap\bar\beta=\emptyset\end{array}$}%
}f_{(m)}(\bar{\alpha},\gamma )f_{(n)}(\bar{\beta},\gamma )Z^{\bar{\alpha}}Z^{%
\bar{\beta}}(Z^{\gamma })^{2} \\
=& m!n!\sum_{\gamma \in \tilde{\Gamma}_{r}}\!\!\!\sum_{%
\mbox{\scriptsize$\begin{array}{c}(\bar\alpha,\bar\beta)\in\tilde \tilde\Gamma_{m-r}\times\tilde\Gamma_{n-r}\\
\bar\alpha\cap\bar\beta=\emptyset\end{array}$}}\!\!\!\!\!f_{(m)}(\bar{\alpha}%
,\gamma )f_{(n)}(\bar{\beta},\gamma )Z^{\bar{\alpha}}Z^{\bar{\beta}%
}(Z^{\gamma })^{2}
\end{align*}%
Since $\gamma \mapsto f_{(m)}(\bar{\alpha},\gamma )f_{(n)}(\bar{\beta}%
,\gamma )Z^{\bar{\alpha}}Z^{\bar{\beta}}(Z^{\gamma })^{2}$ is symmetric, and
similarly for $\bar{\alpha}$ and $\bar{\beta}$, we get
\begin{align*}
\sum_{(\alpha ,\beta )\in \Lambda _{r}^{m,n}}f_{(m)}(\alpha )f_{(n)}(\beta
)Z^{\alpha }Z^{\beta }=& m!n!\frac{1}{r!(m-r)!(n-r)!}\times \\
& \quad \times \sum_{\gamma \in \Gamma_{r}}\sum_{\eta \in \Gamma
_{m-r+n-r}}f_{(m)}(\mathfrak{p}_{m-r}(\eta ),\gamma )f_{(n)}(\mathfrak{q}%
_{n-r}(\eta ),\gamma )Z^{\eta }(Z^{\gamma })^{2}.
\end{align*}%
Then, we have
\begin{align*}
S(f,Z)=& \sum_{m\geq 0}\sum_{n\geq 0}\sum_{r=0}^{m\wedge n}r!%
\begin{pmatrix}
m \\
r%
\end{pmatrix}%
\begin{pmatrix}
n \\
r%
\end{pmatrix}%
\sum_{\eta \in \Gamma _{m+n-2r}}\sum_{\gamma \in \Gamma _{r}}f_{(m)}(%
\mathfrak{p}_{m-r}(\eta ),\gamma )f_{(n)}(\mathfrak{q}_{n-r}(\eta ),\gamma
)(Z^{\gamma })^{2}Z^{\eta } \\
=& \sum_{r\geq 0}r!\sum_{m\geq r}\sum_{n\geq r}%
\begin{pmatrix}
m \\
r%
\end{pmatrix}%
\begin{pmatrix}
n \\
r%
\end{pmatrix}%
\sum_{\eta \in \Gamma _{m+n-2r}}\sum_{\gamma \in \Gamma _{r}}f_{(m)}(%
\mathfrak{p}_{m-r}(\eta ),\gamma )f_{(n)}(\mathfrak{q}_{n-r}(\eta ),\gamma
)(Z^{\gamma })^{2}Z^{\eta }
\end{align*}%
We consider the change of variable $a=m-r$ and $b=m-r+n-r=a+n-r$. We get
\begin{align*}
S(f,Z)=& \sum_{r\geq 0}\sum_{a\geq 0}\sum_{b\geq a}r!%
\begin{pmatrix}
a+r \\
r%
\end{pmatrix}%
\begin{pmatrix}
b-a+r \\
r%
\end{pmatrix}%
\sum_{\eta \in \Gamma _{b}}\sum_{\gamma \in \Gamma _{r}}f_{(a+r)}(\mathfrak{p%
}_{a}(\eta ),\gamma )f_{(b-a+r)}(\mathfrak{q}_{(b-a)}(\eta ),\gamma
)(Z^{\gamma })^{2}Z^{\eta } \\
=& \sum_{r\geq 0}\sum_{b\geq 0}\sum_{\eta \in \Gamma _{b}}\sum_{\gamma \in
\Gamma _{r}}\tilde{A}_{b,r}[f](\eta ,\gamma )(Z^{\gamma })^{2}Z^{\eta }
\end{align*}%
where, for $\eta \in \Gamma _{b}$ and $\gamma \in \Gamma _{r}$,
\begin{equation*}
\tilde{A}_{b,r}[f](\eta ,\gamma )=r!\sum_{a=0}^{b}%
\begin{pmatrix}
a+r \\
r%
\end{pmatrix}%
\begin{pmatrix}
b-a+r \\
r%
\end{pmatrix}%
f_{(a+r)}(\mathfrak{p}_{a}(\eta ),\gamma )f_{(b-a+r)}(\mathfrak{q}%
_{b-a}(\eta ),\gamma ).
\end{equation*}%
We notice that $\gamma \mapsto \tilde{A}_{b,r}[f](\eta ,\gamma )$ is
symmetric but $\eta \mapsto \tilde{A}_{b,r}[f](\eta ,\gamma )$ is not. So,
in order to work with a coefficient $A_{b,r}[f](\eta ,\gamma )$ which is
(separately) symmetric in both variables $\eta $ and $\gamma $, we use the
fact that
\begin{equation*}
\sum_{\eta \in \Gamma _{b}}\tilde{A}_{b,r}[f](\eta ,\gamma )Z^{\eta
}=\sum_{\eta \in \Gamma _{b}}\frac{1}{b!}\sum_{\pi \in \Pi _{b}}\tilde{A}%
_{b,r}[f](\eta _{\pi },\gamma )Z^{\eta }
\end{equation*}%
where $\Pi _{b}$ denotes all the permutations of $(1,\ldots ,b)$ and for $%
\pi \in \Pi _{b}$, $\eta _{\pi }=(\eta _{\pi _{1}},\ldots ,\eta _{\pi _{b}})$%
. Therefore,
\begin{equation*}
S(f,Z)=\sum_{r\geq 0}\sum_{b\geq 0}\sum_{\eta \in \Gamma _{b}}\sum_{\gamma
\in \Gamma _{r}}A_{b,r}[f](\eta ,\gamma )(Z^{\gamma })^{2}Z^{\eta }
\end{equation*}%
and $A_{b,r}[f]$ fulfils formula (\ref{A}). $\square $

\medskip

For $n,r\geq 0$, $m\geq r$, $\eta \in \Gamma _{n}$ and $\rho \in \Gamma _{r}$
we define

\begin{align}
& B_{n,r,m}[f](\eta ,\rho )=\sum_{\beta \in \Gamma _{m-r}}A_{n,m}[f](\eta
,(\rho ,\beta )) \notag \\
& \qquad =\frac{m!}{n!}\sum_{a=0}^{n}%
\begin{pmatrix}
a+m \\
m%
\end{pmatrix}%
\begin{pmatrix}
n-a+m \\
m%
\end{pmatrix}%
\sum_{\pi \in \Pi _{n}}\sum_{\beta \in \Gamma _{m-r}}f_{(a+m)}(\mathfrak{p}%
_{a}(\eta _{\pi }),\rho ,\beta )f_{(n-a+m)}(\mathfrak{q}_{n-a}(\eta _{\pi
}),\rho ,\beta )  \notag \\
& \qquad =\frac{m!}{n!}\sum_{a=0}^{n}%
\begin{pmatrix}
a+m \\
m%
\end{pmatrix}%
\begin{pmatrix}
n-a+m \\
m%
\end{pmatrix}%
\sum_{\pi \in \Pi _{n}}f_{(a+m)}\otimes _{m-r}f_{(n-a+m)}\big((\mathfrak{p}%
_{a}(\eta _{\pi }),\rho ),(\mathfrak{q}_{n-a}(\eta _{\pi }),\rho )\big),
 \label{B}
\end{align}%
$\mathfrak{p}_{a}$ and $\mathfrak{q}_{n-a}$ being defined in (\ref{ppqq}). As a consequence, $B_{n,r,m}[f](\eta ,\rho )=0$ if $\eta $ and $\rho $
have common components and the maps $\eta \mapsto B_{n,r,m}(\eta ,\rho )$
and $\rho \mapsto B_{n,r,m}(\eta ,\rho )$ are both symmetric.

\begin{lemma}
\label{lemma-app2} Let $Y_{i}=Z_{i}^{2}-1$, $i\geq 1$ and let $t_{n,r}(\cdot
,\cdot )$ be defined in (\ref{RA4'}) with respect to $Z=(Z_{i})_{i\in \N}$
and $Y=(Y_{i})_{i\in\N}$. We have
\begin{equation}
|S(f,Z)|^{2}=\sum_{r\geq 0}\sum_{m\geq r}\sum_{n\geq 0}%
\begin{pmatrix}
m \\
r%
\end{pmatrix}%
t_{n,r}(+\infty ,B_{n,r,m}[f]),  \label{Sr}
\end{equation}
$B_{n,r,m}[f]$ being given in (\ref{B}).
\end{lemma}

\textbf{Proof.} We start from the equality
\begin{equation*}
\prod_{i=1}^{m}x_{i}^{2}=1+\sum_{r=1}^{m}\sum_{\Lambda =\{\Lambda
_{1},\ldots ,\Lambda _{r}\}\subset \{1,\ldots
,m\}}\prod_{i=1}^{r}(x_{\Lambda _{i}}^{2}-1).
\end{equation*}%
Fix now $\gamma \in \Gamma _{m}$. For a fixed $r=1,\ldots ,m$ and $\Lambda
=\{\Lambda _{1},\ldots ,\Lambda _{r}\}\subset \{1,\ldots ,m\}$, we set $%
\gamma _{\Lambda }=(\gamma _{\Lambda _{1}},\ldots ,\gamma _{\Lambda _{r}})$.
Then
\begin{eqnarray*}
(Z^{\gamma })^{2} &=&\prod_{i=1}^{m}Z_{\gamma
_{i}}^{2}=1+\sum_{r=1}^{m}\sum_{\Lambda =\{\Lambda _{1},\ldots ,\Lambda
_{r}\}\subset \{1,\ldots ,m\}}\prod_{i=1}^{r}(Z_{\gamma _{\Lambda
_{i}}}^{2}-1) \\
&=&1+\sum_{r=1}^{m}\sum_{\Lambda =\{\Lambda _{1},\ldots ,\Lambda
_{r}\}\subset \{1,\ldots ,m\}}Y^{\gamma _{\Lambda }}.
\end{eqnarray*}%
Then, for $\eta \in \Gamma _{n}$, we can write
\begin{align*}
\sum_{\gamma \in \Gamma _{m}}(Z^{\gamma })^{2}A_{n,m}[f](\eta ,\gamma )=&
\sum_{\gamma \in \Gamma _{m}}A_{n,m}[f](\eta ,\gamma
)+\sum_{r=1}^{m}\sum_{\gamma \in \Gamma _{m}}\sum_{\Lambda =\{\Lambda
_{1},\ldots ,\Lambda _{r}\}\subset \{1,\ldots ,m\}}A_{n,m}[f](\eta ,\gamma
)Y^{\gamma _{\Lambda }} \\
=& \sum_{\gamma \in \Gamma _{m}}A_{n,m}[f](\eta ,\gamma
)+\sum_{r=1}^{m}\sum_{\rho \in \Gamma _{r}}Y^{\rho }%
\begin{pmatrix}
m \\
r%
\end{pmatrix}%
\sum_{\beta \in \Gamma _{m-r}}A_{n,m}[f](\eta ,(\rho ,\beta )),
\end{align*}%
the last inequality following from the fact that $\gamma \mapsto
A_{n,m}[f](\eta ,\gamma )$ is symmetric. We also notice that the case $r=0$
can be easily inserted in the sum of the above r.h.s. (as usual, $\Gamma
_{0}=\{\emptyset \}$ and $Y^{\emptyset }=1$). Then,
\begin{equation*}
\sum_{\gamma \in \Gamma _{m}}(Z^{\gamma })^{2}A_{n,m}[f](\eta ,\gamma
)=\sum_{r=0}^{m}\sum_{\rho \in \Gamma _{r}}Y^{\rho }%
\begin{pmatrix}
m \\
r%
\end{pmatrix}%
\sum_{\beta \in \Gamma _{m-r}}A_{n,m}[f](\eta ,(\rho ,\beta )).
\end{equation*}%
By inserting in (\ref{R1}), we obtain
\begin{align*}
|S(f,Z)|^{2}=& \sum_{m\geq 0}\sum_{n\geq 0}\sum_{\eta \in \Gamma
_{n}}Z^{\eta }\sum_{\gamma \in \Gamma _{m}}(Z^{\gamma })^{2}A_{n,m}[f](\eta
,(\rho ,\beta )) \\
=& \sum_{m\geq 0}\sum_{n\geq 0}\sum_{\eta \in \Gamma _{n}}Z^{\eta
}\sum_{r=0}^{m}\sum_{\rho \in \Gamma _{r}}Y^{\rho }%
\begin{pmatrix}
m \\
r%
\end{pmatrix}%
\sum_{\beta \in \Gamma _{m-r}}A_{n,m}[f](\eta ,(\rho ,\beta )) \\
=& \sum_{r\geq 0}\sum_{m\geq r}\sum_{n\geq 0}%
\begin{pmatrix}
m \\
r%
\end{pmatrix}%
\sum_{\eta \in \Gamma _{n}}\sum_{\rho \in \Gamma _{r}}Z^{\eta }Y^{\rho
}\sum_{\beta \in \Gamma _{m-r}}A_{n,m}[f](\eta ,(\rho ,\beta )) \\
& =\sum_{r\geq 0}\sum_{m\geq r}\sum_{n\geq 0}%
\begin{pmatrix}
m \\
r%
\end{pmatrix}%
t_{n,r}(+\infty ,B_{n,r,m}[f]).
\end{align*}
$\square $

\medskip We take now $c=(c_{(m)})_{m\geq 1}$, with $c_{(m)}\in \mathcal{H}%
^{\otimes m}$, $c_{(m)}$ symmetric and null on all diagonals. We set $%
|c|_{m}=|c_{(m)}|_{m}=\Vert c_{(m)}\Vert _{\mathcal{H}^{\otimes m}}$. Recall
that
\begin{equation*}
\Phi _{m}(m,Z)=\sum_{\alpha \in \Gamma _{m}}c_{(m)}(\alpha )Z^{\alpha }\quad %
\mbox{and}\quad S(c,Z)=\sum_{m\geq 1}\Phi _{m}(c_{(m)},Z).
\end{equation*}%
We also set $\partial _{j}S(c,Z)$ as the derivative of $S(c,Z)$ w.r.t. $%
Z_{j} $ and $D_{j}S(c,Z)$ as the Malliavin derivative in the $j$th
direction. Thus,
\begin{equation*}
D_{j}S(c,Z)=\chi _{j}\partial _{j}S(c,Z).
\end{equation*}%
We recall that we have to deal with $\Vert \tilde{I}(c,Z)\Vert _{p}$ and $%
\Vert I(c,Z)-i(c)\Vert _{p}$, where
\begin{equation*}
\tilde{I}(c,Z)=\sum_{j}\tilde{\chi}_{j}|\partial _{j}S(c,Z)|^{2},\quad
I(c,Z)=\sum_{j}|\partial _{j}S(c,Z)|^{2},\quad i(c)=\sum_{m\geq
1}m!|c|_{m}^{2}=\Vert S(c,Z)\Vert _{2}^{2}.
\end{equation*}%
We also recall that
\begin{equation}
\partial _{j}S(c,Z)=\sum_{m\geq 0}\Phi _{m}(\hat{c}_{(m)}^{j},Z)=S(\hat{c}%
^{j},Z),\quad \mbox{with}\quad \hat{c}_{(m)}^{j}=(1+m)c(\alpha ,j)%
\mbox{\large \bf 1}_{j\notin \alpha },\ \alpha \in \Gamma _{m}.  \label{app2}
\end{equation}%
The case $m=1$ is allowed: just set $\Gamma _{0}=\{\emptyset \}$, $%
Z^{\emptyset }=1$ and $c_{(0)}^{j}=c(j)$. Therefore, we can write
\begin{equation}
\begin{array}{c}
\displaystyle\tilde{I}(c,Z)=\sum_{j}\tilde{\chi}_{j}|S(\hat{c}%
^{j},Z)|^{2}\quad \mbox{and}\quad I(c,Z)=\sum_{j}|S(\hat{c}%
^{j},Z)|^{2},\smallskip \\
\displaystyle\mbox{with}\quad \hat{c}_{(m)}^{j}=(1+m)c(\alpha ,j)%
\mbox{\large \bf 1}_{j\notin \alpha },\ \alpha \in \Gamma _{m}%
\end{array}
\label{Itilde}
\end{equation}

We can then write $\tilde I(c,Z)$ and $I(c,Z)$ as follows.

\begin{lemma}
\label{lemma3-4} $(i)$ Let $\tilde{I}(c,Z)$ and $T_{n,r}(\cdot ,\cdot )$
(associated to $Z_{i}$ and $Y_{i}=Z_{i}^{2}-1$, $i\in \N)$ be defined in (\ref%
{Itilde}) and (\ref{RA4}) respectively. Then,
\begin{equation}
\tilde{I}(c,Z)=\sum_{r\geq 0}\sum_{m\geq r}\sum_{n\geq 0}T_{n,r}(+\infty ,%
\tilde{e}_{n,r,m}[c])  \label{Itilde1}
\end{equation}%
in which $\tilde{e}_{n,r,m}[c]=(\tilde{e}_{n,r,m}^{j}[c])_{j\in {\mathbb{N}}}
$ and for $j\in {\mathbb{N}}$, $\eta\in\Gamma_n$, $\rho\in\Gamma_r$,
\begin{equation}
\begin{array}{ll}
\displaystyle\tilde{e}_{n,r,m}^{j}[c](\eta ,\rho )= & \displaystyle\frac{%
(m+1)(m+1)!}{n!}\sum_{a=0}^{n}%
\begin{pmatrix}
a+m+1 \\
m+1%
\end{pmatrix}%
\begin{pmatrix}
n-a+m+1 \\
m+1%
\end{pmatrix}%
\times \smallskip \\
& \displaystyle\quad \times \sum_{\pi \in \Pi _{n}}c_{(a+m+1)}\otimes
_{m-r}c_{(n-a+m+1)}((\mathfrak{p}_{a}(\eta _{\pi }),\rho ,j),(\mathfrak{q}%
_{n-a}(\eta _{\pi }),\rho ,j)).%
\end{array}
\label{etilde}
\end{equation}%
As a consequence, $\tilde{e}_{n,r,m}^{j}[c](\eta ,\rho )=0$ if $\eta $ and $%
\rho $ have common components or if $j\in (\eta ,\rho ).$ And the maps $\eta
\mapsto \tilde{e}_{n,r,m}^{j}[c](\eta ,\rho )$ and $\rho \mapsto \tilde{e}%
_{n,r,m}^{j}[c](\eta ,\rho )$ are both symmetric.

\medskip

$(ii)$ Let $I(c,Z)$ and $t_{n,r}(\cdot ,\cdot )$ (associated to $Z_{i}$ and $%
Y_{i}=Z_{i}^{2}-1$, $i\in \N$) be defined in (\ref{Itilde}) and (\ref{RA4'})
respectively. Then,
\begin{equation}
I(c,Z)=\sum_{r\geq 0}\sum_{m\geq r}\sum_{n\geq 0}t_{n,r}(+\infty ,e_{n,r,m}[c])
\label{Itilde2}
\end{equation}%
in which, for $\eta\in\Gamma_n$ and $\rho\in\Gamma_r$,
\begin{equation}
\begin{array}{ll}
\displaystyle e_{n,r,m}[c](\eta ,\rho )= & \displaystyle\frac{(m+1)(m+1)!}{n!(m-r+1)%
}\sum_{a=0}^{n}%
\begin{pmatrix}
a+m+1 \\
m+1%
\end{pmatrix}%
\begin{pmatrix}
n-a+m+1 \\
m+1%
\end{pmatrix}%
\times \smallskip \\
& \displaystyle\quad \times \sum_{\pi \in \Pi _{n}}c_{(a+m+1)}\otimes
_{m-r+1}c_{(n-a+m+1)}((\mathfrak{p}_{a}(\eta _{\pi }),\rho ),(\mathfrak{q}%
_{n-a}(\eta _{\pi }),\rho )).%
\end{array}
\label{etildebis}
\end{equation}%
As a consequence, $e_{n,r,m}[c](\eta ,\rho )=0$ if $\eta $ and $\rho $ have
common components and the maps $\eta \mapsto e_{n,r,m}[c](\eta ,\rho )$ and $%
\rho \mapsto e_{n,r,m}[c](\eta ,\rho )$ are both symmetric.
\end{lemma}

\textbf{Proof.} $(i)$ By (\ref{Itilde}), we use Lemma \ref{lemma-app2} with $%
f_{(m)}=\hat{c}_{(m)}^{j}$ and we obtain
\begin{eqnarray*}
\tilde{I}(c,Z) &=&\sum_{j}\tilde{\chi}_{j}|S(\hat{c}^{j},Z)|^{2} \\
&=&\sum_{j}\tilde{\chi}_{j}\sum_{r\geq 0}\sum_{m\geq r}\sum_{n\geq 0}%
\begin{pmatrix}
m \\
r%
\end{pmatrix}%
t_{n,r}(+\infty ,B_{n,r,m}[\hat{c}^{j}]).
\end{eqnarray*}%
We develop $t_{n,r}(+\infty ,B_{n,r,m}[\hat{c}^{j}])$ according to (\ref%
{RA4'}) and we find that the coefficient of $Z^{\eta }Y^{\rho }\tilde{\chi}%
_{j}$ is $B_{n,r,m}[\hat{c}^{j}]$ $(\eta ,\rho ).$ And by definition we have
denoted this quantity by $\tilde{e}_{n,r,m}^{j}[c](\eta ,\rho ).$ Using now
the expression given in (\ref{B}) we obtain
\begin{equation*}
\tilde{e}_{n,r,m}^{j}[c](\eta ,\rho )=\frac{m!}{n!}\sum_{a=0}^{n}%
\begin{pmatrix}
a+m \\
m%
\end{pmatrix}%
\begin{pmatrix}
n-a+m \\
m%
\end{pmatrix}%
\sum_{\pi \in \Pi _{n}}\hat{c}_{(a+m)}^{j}\otimes _{m-r}\hat{c}_{(n-a+m)}^{j}%
\big((\mathfrak{p}_{a}(\eta _{\pi }),\rho ),(\mathfrak{q}_{n-a}(\eta _{\pi
}),\rho )\big).
\end{equation*}%
Since $\hat{c}_{(m)}^{j}(\alpha )=(1+m)c_{(1+m)}(\alpha ,j)$, we get
\begin{equation}
\begin{array}{l}
\displaystyle\hat{c}_{(a+m)}^{j}\otimes _{m-r}\hat{c}_{(n-a+m)}^{j}\big((%
\mathfrak{p}_{a}(\eta _{\pi }),\rho ),(\mathfrak{q}_{n-a}(\eta _{\pi }),\rho
)\big)=\smallskip \\
\displaystyle\quad =(a+m+1)(n-a+m+1)c_{(a+m+1)}\otimes _{m-r}c_{(n-a+m+1)}((%
\mathfrak{p}_{a}(\eta _{\pi }),\rho ,j),(\mathfrak{q}_{n-a}(\eta _{\pi
}),\rho ,j)).%
\end{array}
\label{app}
\end{equation}%
So,
\begin{align*}
\tilde{e}_{n,r,m}^{j}[c](\eta ,\rho )=& \frac{m!}{n!}\sum_{a=0}^{n}%
\begin{pmatrix}
a+m \\
m%
\end{pmatrix}%
\begin{pmatrix}
n-a+m \\
m%
\end{pmatrix}%
(a+m+1)(n-a+m+1)\times \\
& \quad \times \sum_{\pi \in \Pi _{n}}c_{(a+m+1)}\otimes
_{m-r}c_{(n-a+m+1)}((\mathfrak{p}_{a}(\eta _{\pi }),\rho ,j),(\mathfrak{q}%
_{n-a}(\eta _{\pi }),\rho ,j)) \\
=& \frac{(m+1)(m+1)!}{n!}\sum_{a=0}^{n}%
\begin{pmatrix}
a+m+1 \\
m+1%
\end{pmatrix}%
\begin{pmatrix}
n-a+m+1 \\
m+1%
\end{pmatrix}%
\times \\
& \quad \times \sum_{\pi \in \Pi _{n}}c_{(a+m+1)}\otimes
_{m-r}c_{(n-a+m+1)}((\mathfrak{p}_{a}(\eta _{\pi }),\rho ,j),(\mathfrak{q}%
_{n-a}(\eta _{\pi }),\rho ,j))
\end{align*}%
and the proof of $(i)$ is completed.

\smallskip

$(ii)$ By (\ref{Itilde}), we use Lemma \ref{lemma-app2} with $f_{(m)}=\hat
c^j_{(m)}$ and we have the result with $e_{n,r,m}[c](\eta,\rho)= \sum_{j\geq
1}B_{n,r,m}[\hat c^j](\eta,\rho)$. By inserting formula (\ref{B}), we obtain
\begin{equation*}
e_{n,r,m}[c](\eta,\rho)=\sum_{j\geq 1} \frac{m!}{n!}\sum_{a=0}^n
\begin{pmatrix}
a+m \\
m%
\end{pmatrix}
\begin{pmatrix}
n-a+m \\
m%
\end{pmatrix}
\sum_{\pi\in\Pi_n}\hat c^j_{(a+m)}\otimes_{m-r}\hat c^j_{(n-a+m)}\big((%
\mathfrak{p}_a(\eta_\pi),\rho),(\mathfrak{q}_{n-a}(\eta_\pi),\rho)\big).
\end{equation*}
We use now (\ref{app}) and we get
\begin{align*}
&\sum_{j\geq 1} \hat c^j_{(a+m)}\otimes_{m-r}\hat c^j_{(n-a+m)}\big((%
\mathfrak{p}_a(\eta_\pi),\rho),(\mathfrak{q}_{n-a}(\eta_\pi),\rho)\big)= \\
\ &\quad =(a+m+1)(n-a+m+1)\sum_{j\geq
1}c_{(a+m+1)}\otimes_{m-r}c_{(n-a+m+1)}((\mathfrak{p}_a(\eta_\pi),\rho,j),(%
\mathfrak{q}_{n-a}(\eta_\pi),\rho,j)) \\
&\quad =\frac{(a+m+1)(n-a+m+1)}{m-r+1}c_{(a+m+1)} \otimes_{m-r+1}c_{(n-a+m+1)}((\mathfrak{p%
}_a(\eta_\pi),\rho),(\mathfrak{q}_{n-a}(\eta_\pi),\rho)).
\end{align*}
Then,
\begin{align*}
e_{n,r,m}[c](\eta,\rho)
=&\frac{m!}{n!(m-r+1)}\sum_{a=0}^n
\begin{pmatrix}
a+m \\
m%
\end{pmatrix}
\begin{pmatrix}
n-a+m \\
m%
\end{pmatrix}%
(a+m+1)(n-a+m+1)\times\\
&\ \times
\sum_{\pi\in\Pi_n} c_{(a+m+1)}\otimes_{m-r+1}c_{(n-a+m+1)}((%
\mathfrak{p}_a(\eta_\pi),\rho),(\mathfrak{q}_{n-a}(\eta_\pi),\rho)) \\
=&\frac {(m+1)(m+1)!}{n!(m-r+1)}\sum_{a=0}^{n}
\begin{pmatrix}
a+m+1 \\
m+1%
\end{pmatrix}%
\begin{pmatrix}
n-a+m+1 \\
m+1%
\end{pmatrix}\times\\
&\ \times\sum_{\pi\in\Pi_n}c_{(a+m+1)}\otimes_{m-r+1}c_{(n-a+m+1)}((%
\mathfrak{p}_a(\eta_\pi),\rho),(\mathfrak{q}_{n-a}(\eta_\pi),\rho))
\end{align*}
$\square$

\medskip

\subsection{$L^p$ estimates}

We can now use Lemma \ref{lemma-estTt} and finally state the estimate for
the $L^{p}$ norms of $\tilde{I}(c,Z)$ and $I(c,Z)-i(c)$. Here we are obliged
to restrict ourselves to finite series, so we fix $N$ and we suppose that
\begin{equation*}
c(\alpha )=0\quad \mbox{for}\quad \left\vert \alpha \right\vert >N.
\end{equation*}%
In this case we denote $\tilde{I}_{N}(c,Z)$, $I_{N}(c,Z)$ and $i_{N}(c)$
instead of $\tilde{I}(c,Z)$, $I(c,Z)$ and $i(c).$ We also recall that $%
\kappa _{4,\ell }(c)$ is the 4th cumulant given in (\ref{Law6''}) and $%
M_{p}(Z)=1\vee \sup_{k}\Vert Z_{k}\Vert _{p}$.

\begin{lemma}
\label{Est-Ttilde} Fore each $p\geq 1$ there exists a universal constant $%
C_{p}\geq 1$ such that
\begin{equation}
\begin{array}{rl}
\| \widetilde{I}_{N}(c,Z)\| _{p}+\|
I_{N}(c,Z)-i_{N}(c))\| _{p}
\leq &
C_p(1+i_N(c))^{1/2}(N!)^32^NN^{-5/4}\times\smallskip\\
&\displaystyle
\times\Big[\sum_{l=1}^{N}\kappa_{4, l}(c)^{1/4}
+
\sum_{l=1}^{N}\delta_{l}(c)\Big]\label{NL0},
\end{array}
\end{equation}
where $C_p>0$ is a constant depending on $p$ only.
\end{lemma}

\textbf{Proof.} \textbf{Step 1: estimate of $\tilde I_N(c,Z)$.} We recall the expression (\ref{Itilde1}) for $\widetilde{I}_{N}(c,Z)$, based on the coefficients $\tilde{e}%
_{n,r,m}^{j}[c](\eta ,\rho )=(\tilde{e}^j_{n,r,m}[c](\eta ,\rho ))_j$ given in (\ref{etilde}). We notice that $\tilde{e}^j_{n,r,m}[c]\equiv 0$ if, for every $a=0,1,\ldots,n$, one has $c_{(a+m+1)}\otimes_{m-r}c_{(n-a+m+1)}\equiv 0$, and this latter property holds if $a+m+1>N$ or $n-a+m+1>N$. Then, $\tilde{e}^j_{n,r,m}[c]\equiv 0$ if $m+1>N$.

Moreover, for $a=0,1,\ldots,n$ we have
\begin{align*}
\begin{pmatrix}
a+m+1 \\
m+1%
\end{pmatrix}%
\begin{pmatrix}
n-a+m+1 \\
m+1%
\end{pmatrix}%
&=\frac{(a+m+1)!}{(m+1)!}\,\frac{(n-a+m+1)!}{(m+1)!}
\begin{pmatrix}
n\\
a%
\end{pmatrix}%
\frac 1{n!}\\
&\leq \frac{2^n}{n!}\,\frac{(a+m+1)!}{(m+1)!}\,\frac{(n-a+m+1)!}{(m+1)!}.
\end{align*}
Recall that in $\tilde{e}^j_{n,r,m}[c](\eta ,\rho )$ this term multiplies
$c_{(a+m+1)}\otimes_{m-r}c_{(n-a+m+1)}$ which is null for $a+m+1>N$ or $n-a+m+1>N$. So, we consider only the terms with $a$ such that  $a+m+1\leq N$ and $n-a+m+1\leq N$, and notice that this gives the following request:
\begin{equation}\label{ass0}
m+1\leq N\quad\mbox{and}\quad n\leq 2(N-m-1).
\end{equation}
So, we can write
\begin{align}\label{est-coeff-bin}
\begin{pmatrix}
a+m+1 \\
m+1%
\end{pmatrix}%
\begin{pmatrix}
n-a+m+1 \\
m+1%
\end{pmatrix}%
&\leq \frac{2^n}{n!}\,\frac{(N!)^2}{((m+1)!)^2}.
\end{align}
Then,
$$
|\tilde{e}_{n,r,m}^{j}[c](\eta ,\rho )|
\leq
\displaystyle
\frac{2^n(N!)^2}{(n!)^2m!}
\sum_{a=0}^{n}\sum_{\pi \in \Pi _{n}}|c_{(a+m+1)}\otimes
_{m-r}c_{(n-a+m+1)}((\mathfrak{p}_{a}(\eta _{\pi }),\rho ,j),(\mathfrak{q}%
_{n-a}(\eta _{\pi }),\rho ,j))|,
$$
thereby
\begin{equation}\label{NL1}
\begin{array}{l}
\displaystyle
|\tilde{e}_{n,r,m}^{j}[c](\eta ,\rho )|^2 \smallskip\\
\leq
\displaystyle
\Big(\frac{2^n(N!)^2}{(n!)^2m!}\Big)^2\times n\times n!
\sum_{a=0}^{n}\sum_{\pi \in \Pi _{n}}|c_{(a+m+1)}\otimes
_{m-r}c_{(n-a+m+1)}((\mathfrak{p}_{a}(\eta _{\pi }),\rho ,j),(\mathfrak{q}%
_{n-a}(\eta _{\pi }),\rho ,j))|^2.
\end{array}
\end{equation}%
Now we have to distinguish two cases: $m\geq r+1$ and $m=r.$ We assume first
that $m\geq r+1$ and we use (\ref{acc9}) in order to obtain
$$
\begin{array}{ll}
\displaystyle
|\tilde{e}_{n,r,m}[c]|^2_{n,r,\infty}
&
\displaystyle
=\sum_{j\geq 1}\sum_{\eta\in\Gamma_n}\sum_{\rho\in\Gamma_r}|\tilde{e}^j_{n,r,m}[c](\eta,\rho)|^2
\smallskip\nonumber\\
&\leq
\displaystyle
\Big(\frac{2^n(N!)^2}{(n!)^2m!}\Big)^2\times n\times n!
\sum_{a=0}^{n}|c_{(a+m+1)}\otimes
_{m-r}c_{(n-a+m+1)}|^2\smallskip\nonumber\\
&\leq
\displaystyle
\Big(\frac{2^n(N!)^2}{(n!)^2m!}\Big)^2\times n\times n!
2\sum_{a=0}^{n}\frac{\kappa_{4, a+m+1}(c)}{((a+m+1)!)^2\begin{pmatrix}
a+m+1 \\
m-r%
\end{pmatrix}^2}%
\end{array}
$$
so that

\begin{equation}\label{est1}
\begin{array}{ll}
\displaystyle
|\tilde{e}_{n,r,m}[c]|^2_{n,r,\infty}
&
\leq
\displaystyle
\Big(\frac{2^n(N!)^2}{(n!)^2m!}\Big)^2\times n\times n!\times
\frac 2{((m+1)!)^2}\sum_{l=1}^{N}\kappa_{4, l}(c).
\end{array}
\end{equation}

If instead $m=r$, (\ref{NL1}) gives
$$
\begin{array}{ll}
\displaystyle
|\tilde{e}_{n,r,r}[c]|^2_{n,r,\infty}
&\leq
\displaystyle
\Big(\frac{2^n(N!)^2}{(n!)^2r!}\Big)^2\times n\times n!\times\smallskip\\
&\quad\displaystyle\times
\sum_{a=0}^{n}\sum_{\pi \in \Pi _{n}}\sum_{j\geq 1}\sum_{\eta\in\Gamma_n}\sum_{\rho\in\Gamma_r}|c_{(a+r+1)}(\mathfrak{p}_{a}(\eta _{\pi }),\rho ,j)
c_{(n-a+r+1)}(\mathfrak{q}%
_{n-a}(\eta _{\pi }),\rho ,j)
|^2\smallskip\\
&\leq
\displaystyle
\Big(\frac{2^n(N!)^2}{(n!)^2r!}\Big)^2\times n\times n!\times
\sum_{a=0}^{n}\delta_{a+r+1}^2(c)|c|_{n-a+r+1}^2
\end{array}
$$
and we obtain
\begin{equation}\label{est2}
\begin{array}{ll}
\displaystyle
|\tilde{e}_{n,r,r}[c]|^2_{n,r,\infty}
&\leq
\displaystyle
\Big(\frac{2^n(N!)^2}{(n!)^2r!}\Big)^2\times n\times n!\times i_N(c)
\sum_{l=1}^{N}\delta_{l}^2(c).
\end{array}
\end{equation}
Therefore, by (\ref{Itilde1}), recalling the conditions (\eqref{ass0}), and by using (\ref{RA9}), we can write
\begin{align*}
\|\tilde{I}(c,Z)\|_p
\leq&\sum_{r= 0}^{N-1}\sum_{m=r}^{N-1}\sum_{n=0}^{2(N-m-1)}\|T_{n,r}(+\infty ,%
\tilde{e}_{n,r,m}[c])\|_p\\
\leq&\sum_{r= 0}^{N-1}\sum_{m=r}^{N-1}\sum_{n= 0}^{2(N-m-1)}
C_p(4^{n+r}(n+r)!)^{1/2}(\sqrt 2\,b_p M_p)^{n+r}|\tilde e_{n,r,m}[c]|_{n,r,\infty},
\end{align*}
in which $C_p$ is a constant depending on $p$ only (which may vary in next lines) and
$$
M_p=M_p(Z,Y)=\sup_{k}\|Z_k\|_p\vee\|Y_k\|_p\leq 2\sup_k\|Z_k\|_{2p}^2=:A_p.
$$
We now split the cases $m>r+1$ and $m=r$ and we use the estimates (\ref{est1}) and (\ref{est2}). Then
\begin{align*}
\|\tilde{I}(c,Z)\|_p
\leq&C_p\sum_{r= 0}^{N-2}\sum_{m=r+1}^{N-1}\sum_{n=0}^{2(N-m-1)}
(4^{n+r}(n+r)!)^{1/2}(\sqrt 2\,b_p A_p)^{n+r}
\frac{2^n(N!)^2}{(n!)^2m!}\sqrt{n\,n!}
\frac 1{(m+1)!}\Big(\sum_{l=1}^{N}\kappa_{4, l}(c)\Big)^{1/2}\\
& + C_p
\sum_{r= 0}^{N-1}\sum_{n=0}^{2(N-r-1)}
(4^{n+r}(n+r)!)^{1/2}(\sqrt 2\,b_p A_p)^{n+r}
\frac{2^n(N!)^2}{(n!)^2r!}\sqrt{n\, n!}\Big( i_N(c)
\sum_{l=1}^{N}\delta_{l}^2(c)\Big)^{1/2}\\
\leq&C_p
\sum_{r= 0}^{N-1}\sum_{n=0}^{2(N-r-1)}
(4^{n+r}(n+r)!)^{1/2}(\sqrt 2\,b_p A_p)^{n+r}
\frac{2^n(N!)^2}{(n!)^2r!}\sqrt{n\,n!}\,\times\\
&\times
\Big[\sum_{l=1}^{N}\kappa_{4, l}(c)^{1/2}
+\sqrt{i_N(c)}\,
\sum_{l=1}^{N}\delta_{l}(c)\Big]\\
\leq
&C_p(N!)^2((2N-2)!)^{1/2}
\Big[\sum_{l=1}^{N}\kappa_{4, l}(c)^{1/2}
+\sqrt{i_N(c)}\,
\sum_{l=1}^{N}\delta_{l}(c)\Big],
\end{align*}
in which we have used the fact that $(n+r)!\leq (2N-2)!$ for $r$ and $n$ in the range of the above series.

\smallskip

\textbf{Step 2: estimate of $I_N(c,Z)-i_N(c)$.} We recall the expression (\ref{Itilde2}) for $I_{N}(c,Z)$, based on the coefficients $e_{n,r,m}[c](\eta ,\rho )$ given in (\ref{etildebis}). We notice that the term $i_N(c)$ is actually the term in the series (\ref{Itilde2}) when one takes $n=r=0$, so that
$$
I_N(c,Z)-i_N(c)
=\sum_{r\geq 0}\sum_{m\geq r}\sum_{n\geq 0\vee(1-r)}t_{n,r}(+\infty ,e_{n,r,m}[c])
$$
Following the same arguments developed in Step 1, we can say that $e_{n,r,m}[c]\equiv 0$ if $m+1>N$ and the constraints in (\ref{ass0}) hold. And by using (\ref{est-coeff-bin}), we obtain
$$
\begin{array}{ll}
\displaystyle
|e_{n,r,m}[c](\eta ,\rho )|
\leq & \displaystyle
\frac{2^n(N!)^2}{(n!)^2m!(m-r+1)}
\sum_{a=0}^{n}%
\sum_{\pi \in \Pi _{n}}|c_{(a+m+1)}\otimes
_{m-r+1}c_{(n-a+m+1)}((\mathfrak{p}_{a}(\eta _{\pi }),\rho ),(\mathfrak{q}%
_{n-a}(\eta _{\pi }),\rho ))|,
\end{array}
$$
so that
\begin{equation}\label{NL2}
\begin{array}{l}
\displaystyle
|e_{n,r,m}[c](\eta ,\rho )|^2
\smallskip\\
\leq \displaystyle
\Big(\frac{2^n(N!)^2}{(n!)^2m!(m-r+1)}\Big)^2\times n\times n!
\sum_{a=0}^{n}%
\sum_{\pi \in \Pi _{n}}|c_{(a+m+1)}\otimes
_{m-r+1}c_{(n-a+m+1)}((\mathfrak{p}_{a}(\eta _{\pi }),\rho ),(\mathfrak{q}%
_{n-a}(\eta _{\pi }),\rho ))|^2.
\end{array}
\end{equation}
We have now to split the case $r>0$ and $r=0$. For $r>0$, we use (\ref{acc9}) and, similarly as before, we obtain
$$
\begin{array}{l}
\displaystyle
|e_{n,r,m}[c]|^2_{n,r,\infty}
\leq
\displaystyle
\Big(\frac{2^n(N!)^2}{(n!)^2m!(m-r+1)}\Big)^2\times n\times n!\times
\frac 2{((m+1)!)^2}\sum_{l=m+1}^{n+m+1}\kappa_{4, l}(c).
\end{array}
$$
When $r=0$ (recall that in this case $n\geq 1$) we have a different behavior in the sum as $a=0,\ldots,n$. In fact, for $a=1,\ldots,n-1$ we use again (\ref{acc9}) and we obtain the same estimate as before. But for $a=0$ and $a=n$, we cannot use (\ref{acc9}) but we can use (\ref{acc9-1}). So, we obtain
$$
\begin{array}{ll}
\displaystyle
|e_{n,0,m}[c]|^2_{n,0,\infty}
&\leq
\displaystyle
\Big(\frac{2^n(N!)^2}{(n!)^2m!(m+1)}\Big)^2\times n\times n!\times
\frac 2{((m+1)!)^2}\sum_{l=m+2}^{n+m}\kappa_{4, l}(c)+\smallskip\\
&\displaystyle
+\Big(\frac{2^n(N!)^2}{(n!)^2m!(m+1)}\Big)^2\times n\times n!\times
2|c_{(m+1)}\otimes_{m+1}c_{(n+m+1)}|^2\smallskip\\
&
\leq
\displaystyle
\Big(\frac{2^n(N!)^2}{(n!)^2m!(m+1)}\Big)^2\times n\times n!\Big(
\frac 2{((m+1)!)^2}\sum_{l=m+2}^{n+m}\kappa_{4, l}(c)
+
\frac 2{(m+1)!}|c|_{m+1}^2\kappa_{4, n+m+1}(c)^{1/2}\Big)\smallskip\\
&
\leq
\displaystyle
2(1+i_N(c))\Big(\frac{2^n(N!)^2}{(n!)^2m!(m+1)}\Big)^2\times \frac{n\,n!}{(m+1)!}
\sum_{l=m+1}^{n+m+1}\kappa_{4, l}(c)^{1/2}
\end{array}
$$
By resuming, for $n,r\geq 0$ such that $n+r\geq 1$ we have
$$
\begin{array}{l}
\displaystyle
|e_{n,r,m}[c]|^2_{n,r,\infty}
\leq
\displaystyle
2(1+i_N(c))\Big(\frac{2^n(N!)^2}{(n!)^2m!(m-r+1)}\Big)^2
\frac 2{n\,n!((m+1)!)^2}\sum_{l=1}^{N}\kappa_{4, l}(c)^{1/2}.
\end{array}
$$
Therefore, by (\ref{Itilde2}), recalling the conditions (\ref{ass0}), and by using (\ref{RA9'}) and the estimate $M_p(Z,Y)\leq A_p$,
we can write
\begin{align*}
&\|I_N(c,Z)-i_N(c)\|_p
\leq\sum_{r= 0}^{N-1}\sum_{m=r}^{N-1}\sum_{n=0}^{2(N-m-1)}\I_{n+r\geq 1}\|t_{n,r}(+\infty ,%
e_{n,r,m}[c])\|_p\\
&\quad \leq\sum_{r= 0}^{N-1}\sum_{m=r}^{N-1}\sum_{n= 0}^{2(N-m-1)}
\I_{n+r\geq 1}\big((n+r)!\big)^{1/2}\,(
\sqrt 2\,b_{p}A_p)^{(r+n)}|e_{n,r,m}[c]|_{n,r,\infty}\\
&\quad \leq\sum_{r= 0}^{N-1}\sum_{m=r}^{N-1}\sum_{n= 0}^{2(N-m-1)}
\big((n+r)!\big)^{1/2}\,(
\sqrt 2\,b_{p}A_p)^{(r+n)}
\frac{2(1+i_N(c))^{1/2}2^n(N!)^2\sqrt{n\,n!}}{(n!)^2m!(m-r+1)(m+1)!}
\sum_{l=1}^{N}\kappa_{4, l}(c)^{1/4}\\
&\quad\leq
C_p(1+i_N(c))^{1/2}(N!)^2\sum_{l=1}^N\kappa_{4,l}(c)^{1/4}\times \mathcal{S}
\end{align*}
where $\mathcal{S}$ is a sum which has a behavior similar  to the one studied in step 1. So,
\begin{align*}
\|I_N(c,Z)-i_N(c)\|_p
\leq
C_p(N!)^2((2N-2)!)^{1/2}
\sum_{l=1}^{N}\kappa_{4, l}(c)^{1/4}.
\end{align*}
The Stirling's approximation formula now gives
$$
\exists\ \lim_{N\to\infty}\frac{((2N-2)!)^{1/2}}{2^NN^{-5/4}N!}\in(0,1)
$$
so $((2N-2)!)^{1/2}\leq 2^NN^{-5/4}N!$ and the statement finally holds.
$\square $

\end{document}